\documentclass[11pt,a4paper,twoside]{book}
\usepackage[margin=2.5cm]{geometry}
\usepackage{graphicx,url}
\usepackage{amssymb,bm}
\usepackage[english]{babel}

\def\no{\noindent}
\newtheorem{defi}{Definition}[chapter]

\newtheorem{rem}{Remark}[chapter]
\newtheorem{theo}{Theorem}[chapter]
\newtheorem{cor}{Corollary}[chapter]
\newtheorem{lem}{Lemma}[chapter]
\def\proof{\no\underline{Proof.}~}
\def\QED{\mbox{$\square$}}
\def\no{\noindent}   
\def\RR{\mathbb{R}} 
\def\CC{\mathbb{C}}

\def\DD{\mathbb{D}}
\def\pmatrix{\left(\begin{array}}
\def\endpmatrix{\end{array}\right)} 
\def\dd{{\mathrm d}}
\def\bfa{{\bf a}}
\def\bfb{{\bf b}}
\def\bfc{{\bf c}}

\def\bfe{{\bf e}}

\def\bfq{{\bf q}}
\def\bfp{{\bf p}}

\def\bfx{{\bf x}}
\def\bfy{{\bf y}}

\def\aa{{\alpha}}
\def\bb{{\beta}}
\def\bfo{{\bf 0}}
\def\ee{{\mathrm e}}

\def\eps{\varepsilon}
\def\bfeta{\bm{\eta}}
\def\bfgamma{\bm{\gamma}}

\def\II{{\cal I}}
\def\P{{\cal P}}
\def\dd{{\mathrm{d}}}
\def\rank{{\rm rank}}
\def\diag{{\rm diag}}
\def\div{{\rm div}}
\def\vol{{\rm vol}}

\title{\bf\huge \underline{LINE INTEGRAL METHODS}\\[2mm]  {\LARGE and their application to the numerical solution
of  conservative problems}}

\bigskip
\bigskip
\author{\begin{tabular}{ccc}
Luigi Brugnano  && Felice Iavernaro\\
&\qquad and \qquad\qquad &\\[-.35cm]
 University of Firenze, Italy &&
 University of Bari, Italy\end{tabular}
}
\date{~\\ \vspace{8cm} ~\\ Lecture Notes of the course\\
held at the Academy of Mathematics and Systems Science\\Chinese Academy of Sciences in Beijing\\ on December 27, 2012\,--\,January 4, 2013}

\begin{document}
\pagenumbering{roman}

\large
\maketitle

\subsection*{}

\begin{figure}
\bigskip
\bigskip
\bigskip
\bigskip
\bigskip
\bigskip
\bigskip
\bigskip

\centerline{\includegraphics[width=16cm,height=8cm]{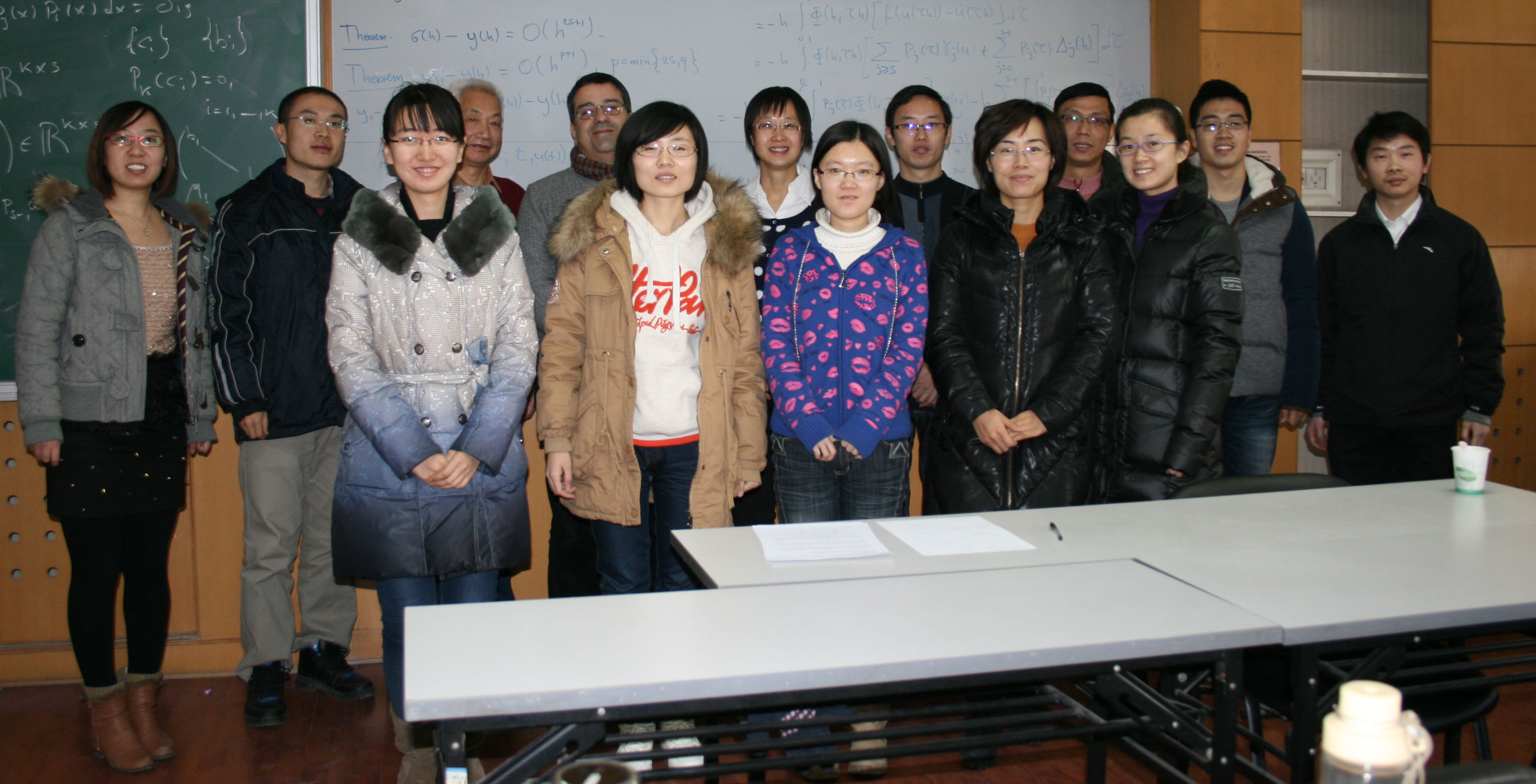}}
\end{figure}
\bigskip
\bigskip
\bigskip
\subsection*{Acknowledgements} The first author wish to thank Dr.\,Yajuan Sun for the kind invitation to deliver this course.

\bigskip
\bigskip
\bigskip
\no We devote this work to honour the memory of Professor Donato Trigiante: a brilliant scientist and a fine man, who passed down to us his love for scientific investigation.

\tableofcontents

\newpage
\pagenumbering{arabic}

\setcounter{page}{0}
\chapter{Geometric Integration} 

In this chapter, we will discuss the basic issues about Geometric Integration, giving a concrete motivation to look for {\em energy-conserving} methods, for the efficient numerical solution of Hamiltonian problems. In particular we will focus on the basic idea  Hamiltonian Boundary Value Methods (HBVMs) relies on, i.e., the definition of discrete line integrals. The material of tis chapter is based on references \cite{IP08,IT09,BIT09,BIT09_1,BI12_1}.

\section{Introduction}
The numerical solution of conservative problems is an active field of investigation dealing with the geometrical properties of the discrete vector field induced by numerical methods. The final goal is to reproduce, in the discrete setting, a number of  geometrical properties shared by the original continuous problem. Because of this reason, it has become customary to refer to this field of investigation as {\em geometric integration}, even though this concept can be led back to the early work of G.\,Dahlquist on differential equations, aimed at reproducing the asymptotic stability of equilibria for the trajectories defined by a numerical method, according to the well-known linear stability analysis (see, e.g., \cite{BT98}).

In particular, we shall deal with the numerical solution of {\em Hamiltonian problems}, which are encountered in many real-life applications, ranging from the nano-scale of molecular dynamics, to the macro-scale of celestial mechanics. Such problems have the following general form, \begin{equation}\label{ham}y' = J\nabla H(y), \qquad y(0) = y_0\in\RR^{2m},\end{equation} where $J^T=-J =J^{-1}$ is a constant, orthogonal and skew-symmetric matrix, usually given by 
\begin{equation}\label{J}
J=\pmatrix{cc} 0 &I\\ -I&0\endpmatrix
\end{equation} 
(here $I$ is the identity matrix of dimension $m$). In such a case, we speak about a problem in {\em canonical form}. The scalar function $H(y)$ is the {\em Hamiltonian} of the problem and its value is constant during the motion, namely  $$H(y(t))\equiv H(y_0), \qquad \forall t\ge 0,$$ for the solution of (\ref{ham}).  Indeed, one has:
\begin{equation}\label{Ht0}
\frac{\dd}{\dd t}H(y(t)) = \nabla H(y(t))^Ty'(t) = \nabla H(y(t))^TJ\nabla H(y(t))=0, \qquad \forall t\ge 0.
\end{equation}
Often, the Hamiltonian $H$ is also called the {\em energy}, since for isolated mechanical systems it has the physical meaning of total energy. Consequently, {\em energy  conservation} is an important feature in the simulation of such problems.
%% For the continuous dynamical system, it can be shown that energy conservation derives from the property of {\em symplecticity} of the map (see, e.g., \cite{GPS01}). 
The state vector of a Hamiltonian system splits in two $m$-length components
$$y = \pmatrix{c} q\\ p\endpmatrix,$$
where $q$ and $p$ are the vectors of generalized positions and momenta, respectively. Consequently, (\ref{ham})-(\ref{J}) becomes $$q' = \nabla_p H(q,p), \qquad p' = -\nabla_q H(q,p).$$
Depending on the case, we shall use both notations. 

Another important feature of Hamiltonian dynamical systems is that they posses a {\em symplectic} structure. To introduce this property we need a copule of ingredients:
\begin{itemize}
\item[-] The \textit{flow of the system}: it is the map acting on  the phase space $\RR^{2m}$ as $$\phi_t: y_0 \in \RR^{2m} \rightarrow y(t)\in \RR^{2m},$$
where $y(t)$ is the solution at time $t$ of (\ref{ham}) originating from the initial condition $y_0$. Differentiating both sides of (\ref{ham}) by $y_0$ and observing that $$\frac{\partial y(t)}{\partial y_0}=\frac{\partial \phi_t(y_0)}{\partial y_0}\equiv \phi'_t(y_0),$$ we see that the Jacobian matrix of the flow $\phi_t$  is the solution of the variational equation associated with (\ref{ham}), namely
\begin{equation}
\label{variational}
\frac{d}{dt} A(t)=J\nabla^2 H(y(t)) A(t), \qquad A(0)=I,
\end{equation}
where $\nabla^2 H(y)$ is the Hessian matrix of $H(y)$.
\item[-] The definition of  a \textit{symplectic transformation}: a map $u= (q,p)\in \RR^{2m} \mapsto  u(q,p) \RR^{2m}$ is said symplectic  if its Jacobian matrix $u'(q,p) \in \RR^{2m \times 2m}$ is a symplectic matrix, that is
$$
u'(q,p)^TJ u'(q,p) =J, \quad \mbox{for all }\, q,p\in\RR^m.
$$
\end{itemize}
That said, it is not difficult to prove that, under regularity assumptions on $H(q,p)$, the flow associated to a Hamiltonian system is symplectic. Indeed, setting $$A(t)=\frac{\partial \phi_t}{\partial y_0},$$ and considering (\ref{variational}), on has that
$$
\begin{array}{rl}
\displaystyle \frac{d}{dt}\left( A(t)^T J  A(t) \right) & \displaystyle =
\left( \frac{d}{dt} A(t)^T J  A(t) \right)+
\left( A(t)^T J  \frac{d}{dt}  A(t) \right) \\[.4cm]
\displaystyle  &= \displaystyle \left( A(t)^T \nabla^2 H(y(t))J^T J  A(t) \right) +
\left( A(t)^T J J \nabla^2 H(y(t)) A(t) \right)=0.
\end{array}
$$
Therefore $$A(t)^T J  A(t)\equiv A(0)^T J  A(0)=J.$$ The converse of the above property is also true: if the flow associated with a dynamical system $\dot y=f(y)$ defined on $\RR^{2m}$  is symplectic then necessarily $f(y)=J\nabla H(y)$ for a suitable scalar function $H(y)$. Consequently, conservation of $H(y)$ follows, by virtue of (\ref{Ht0}).

Symplecticity has relevant implications on the dynamics of Hamiltonian systems. Among the most important are:
\begin{itemize}
\item[(i)] \textit{Canonical transformations}.  A change of variables $z=\psi(y)$ is \textit{canonical}, namely it preserve the structure of (\ref{ham}), if and only if it is symplectic. Canonical transformations were known from Jacobi and used to recast (\ref{ham}) in simpler form.
\item[(ii)] \textit{Volume preservation}.  The flow $\phi_t$ of a Hamiltonian system  is volume preserving in phase space. Recall that if $V$ is a (suitable) domain of $\RR^{2m}$, we have:
    $$
    \vol (V) =\int_V \dd y, \qquad \vol (\phi_t(V)) = \int_{\phi_t(V)} \dd y = \int_V \left| \det \frac{ \partial \phi_t(y)}{\partial y} \right|  \dd y.
    $$
However, since  $\frac{ \partial \phi_t(y)}{\partial y}\equiv A(t)$ is a symplectic matrix, from $A(t)^T J A(t)=J$ it follows that $\det(A(t))^2=1$ for any $t$ and, 
%, since $\det A(0)=\det I=1$ we have, in particular, $\det A(t)=1$ and 
hence, $ \vol (\phi_t(V))= \vol(V)$.

More in general, volume  preservation is a characteristic feature of divergence-free vector fields.
Recall that the divergence of a vector field $f:\RR^n\rightarrow \RR^n$ is the trace of its Jacobian matrix:
    $$
    \div f(y)= \frac{\partial f_1}{\partial y_1}+\frac{\partial f_2}{\partial y_2}+\dots+\frac{\partial f_n}{\partial y_n}.
    $$
    The vector field $J\nabla H$ associated with a Hamiltonian system has zero divergence. In fact, considering that $J\nabla H =[\frac{\partial H}{\partial p_1},\dots,\frac{\partial H}{\partial p_m},-\frac{\partial H}{\partial q_1},\dots,-\frac{\partial H}{\partial q_m}]^T$ we obtain
    $$
    \div \, \nabla H= \frac{\partial^2 H}{\partial q_1 \partial p_1}+\dots+\frac{\partial^2 H}{\partial q_m \partial p_m}-\frac{\partial^2 H}{\partial p_1 \partial q_1}- \dots - \frac{\partial^2 H}{\partial p_m \partial q_m}=0
    $$
    since the partial derivatives commute.
An important consequence of the previous property is Liouville's theorem, which states that the flow $\phi_t$ associated with a divergence-free vector field $f:\RR^n \rightarrow \RR^n$ is volume preserving.
\end{itemize}

The above properties and the fact that symplecticity is a characterizing property of Hamiltonian systems somehow reinforces the search of symplectic methods for their numerical integration. A one-step method $$y_1=\Phi_h(y_0)$$ is per se a transformation of the phase space. Therefore the method is symplectic if $\Phi_h$ is a symplectic map. An important consequence of symplecticity in Runge-Kutta methods is  the conservation of all quadratic first integral of a Hamiltonian system. 
 
 A first integral for system (\ref{ham}) is a scalar function $I(y)$ which remains constant if evaluated  along any solution $y(t)$ of (\ref{ham}): $I(y(t))=I(y_0)$ or, equivalently, $$\nabla I(y)^T J \nabla H(y)=0,\qquad \mbox{for any}\qquad y.$$ A quadratic first integral takes the form $I(y)=y^TCy$, with $C$ a symmetric matrix.

As previously seen, the most noticeable first integral of a Hamiltonian system is the Hamiltonian function itself. It is worth noticing that while in the continuous setting energy conservation derives from the property of {\em symplecticity} of the flow (see, e.g., \cite{GPS01}), as sketched above, the same is no longer true in the discrete setting: a symplectic integrator is not able to yield energy conservation in general. Consequently, devising {\em energy  conservation} methods form an important branch of the geometric integration.

%Concerning the numerical integration of (\ref{ham}), two main lines of investigation have been the definition and the study of {\em symplectic} methods and {\em energy-conserving} methods, respectively. This reflects the fact that the symplecticity of the flow associated with (\ref{ham}) and energy conservation are the most relevant features characterizing a Hamiltonian system.

Symplectic methods can be found in early work of Gr\"obner (see, e.g., \cite{G75}). Symplectic Runge-Kutta methods have been then studied by Feng Kang \cite{FK85}, Sanz Serna \cite{SS88}, and Suris \cite {Suris89}. Such methods are obtained by imposing that the discrete map, associated with a given numerical method,  is symplectic, as is the continuous one. In particular, in \cite{SS88} an easy criterion for simplecticity is provided, for an $s$-stage Runge-Kutta method with tableau given by
\begin{equation}\label{RK}
\begin{array}{c|c}
\bfc &A \\ \hline & \bfb^T\end{array}
\end{equation}
where, as usual, $\bfc=(c_i)\in\RR^s$ is the vector of the abscissae, $\bfb=(b_i)\in\RR^s$ is the vector of the weights, and $A=(a_{ij})\in\RR^{s\times s}$ is the corresponding Butcher matrix. 
\begin{theo}\label{SSth}(\ref{RK}) is symplectic if and only if, by setting $B=\diag(\bfb)$, one has:
\begin{equation}\label{RKsymp}
BA+A^TB = \bfb\bfb^T.
\end{equation}
\end{theo}
Since  for the continuous map symplecticity implies energy-conservation, though this is no more true for the discrete one, then one expects that at least {\em something similar} happens for the discrete map as well. As a matter of fact, under suitable assumptions, it can be proved that, when a symplectic method is used with a constant step-size, the numerical solution satisfies a perturbed Hamiltonian problem, thus providing a quasi-conservation property over ``exponentially long times'' \cite{BG94}. Even though this is an interesting feature, nonetheless, it constitutes a somewhat weak stability result since, in general, it does not extend to infinite intervals.

\begin{figure}
\centerline{\includegraphics[width=12cm,height=8cm]{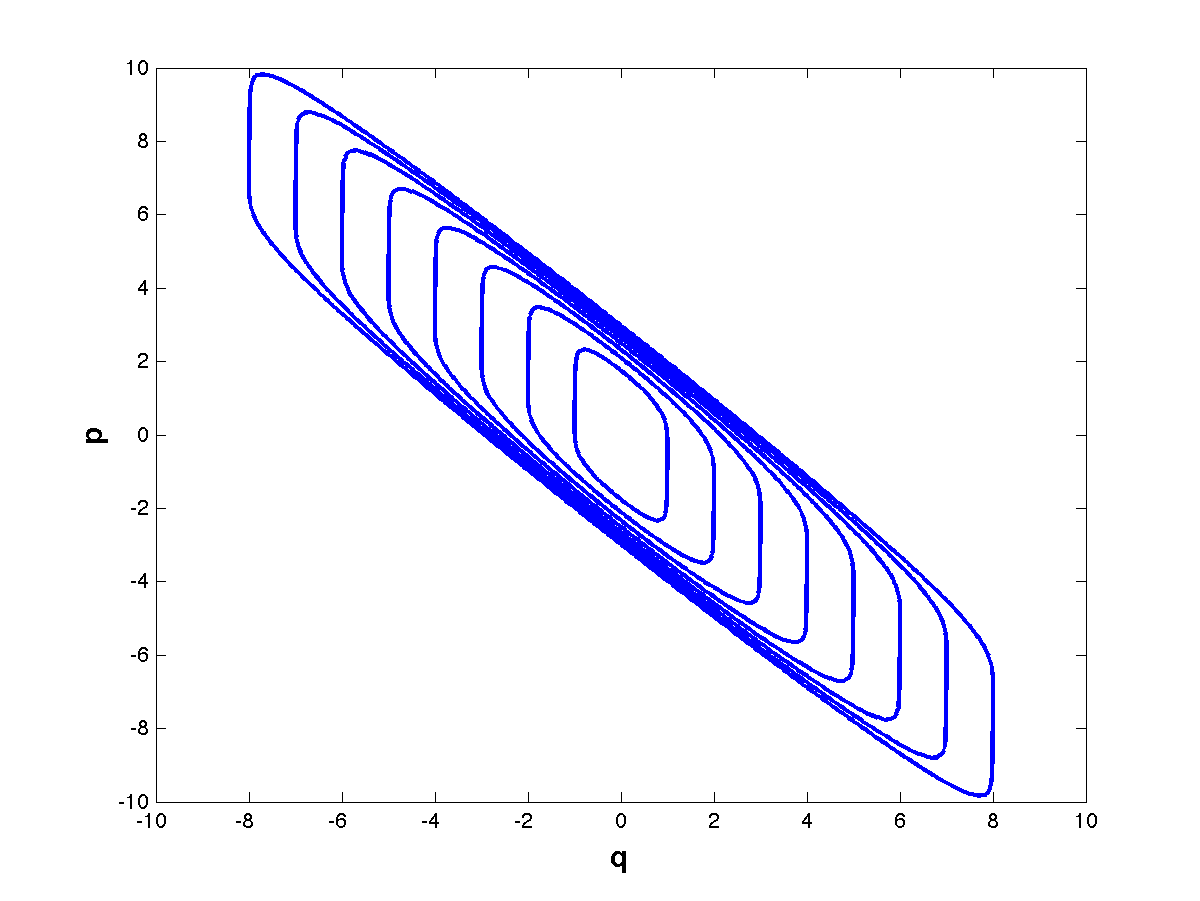}}
\caption{Level curves for problem (\ref{Hex})--(\ref{initex}).}
\label{level}
\end{figure}
\begin{figure}
\centerline{\includegraphics[width=12cm,height=8cm]{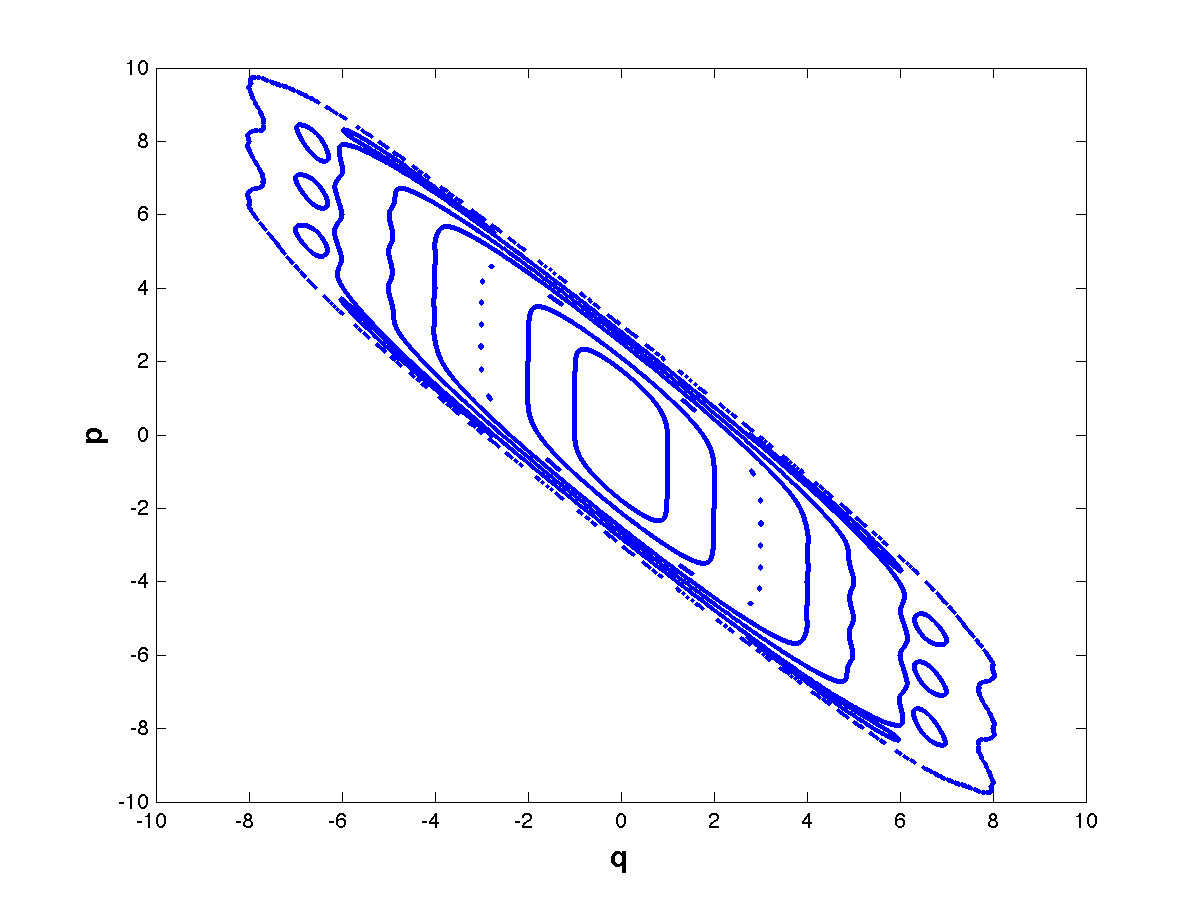}}
\caption{2-stage Gauss method, $h=10^{-3}$, for approximating problem (\ref{Hex})--(\ref{initex}). }
\label{gauss2}
\end{figure}

Moreover, the perturbed dynamical system could be not ``so close'' to the original one, meaning that,  if the stepsize $h$ is not small enough, the perturbed Hamiltonian could not correctly approximate the exact one.  As an example, consider the problem defined by the Hamiltonian \cite{BIT11}
\begin{equation}\label{Hex}H(q,p) = p^2 + (\beta q)^2 +\alpha (q+p)^{2n}.\end{equation}
The corresponding dynamical system has exactly one (marginally stable) equilibrium at the origin. 
Let us select the following parameters
\begin{equation}\label{parex}
\beta = 10, \qquad \alpha = 1, \qquad n = 4,
\end{equation}
and suppose we are interested in approximating the level curves of the Hamiltonian (shown in Figure~\ref{level}) passing from the points
\begin{equation}\label{initex} (q_0,p_0) = (i,-i), \qquad i = 1,\dots,8.\end{equation}
This can be done by integrating the trajectories starting at such initial points, for the corresponding Hamiltonian system but, if we use the {\em symplectic} 2-stage Gauss method, with stepsize $h=10^{-3}$, we obtain the phase portrait depicted in Figure~\ref{gauss2} which is clearly wrong.\footnote{Additional examples may be found in reference \cite{BIS10}.}

A way to get rid of this problem is to directly look for {\em energy-conserving} methods, able to provide an exact conservation of the Hamiltonian function along the numerical trajectory.

The very first attempts to face this problem were based on projection techniques coupled with standard non conservative numerical methods. However, it is well-known that this approach suffers from many drawbacks, in that this is usually not enough to correctly reproduce the dynamics (see, e.g., \cite[p.\,111]{HLW06}).

A completely new approach is represented by {\em discrete gradient methods}, which are based upon the definition of a discrete counterpart of the gradient operator, so that energy conservation of the numerical solution is guaranteed at each step and for any choice of the integration step-size \cite{Go96,McLQR99}.

A different approach is based on the concept of {\em time finite element methods} \cite{Hu72}, where one finds local Galerkin approximations on each subinterval of a given mesh of size $h$ for the equation (\ref{ham}). This, in turn, has led to the definition of energy-conserving Runge-Kutta methods \cite{BS00,Bo97,TC07,TS12}.

A partially related approach is given by {\em discrete line integral methods} \cite{IP07,IP08,IT09}, where the key idea is to exploit the relation between the method itself and the
{\em discrete line integral}, i.e., the discrete counterpart of the line integral in conservative vector fields. This tool  
 yields exact conservation for polynomial Hamiltonians of arbitrarily high-degree, and  results in the class of methods later named {\em Hamiltonian Boundary Value Methods (HBVMs)}, which have been developed in a series of papers \cite{BIT09,BIT09_1,BIS10,BIT10,BIT10_1,BIT11,BIT12,BIT12_1}.

Another approach, strictly related to the latter one, is given by the {\em Averaged Vector Field} method  \cite{QMcL08,CMcLMcLOQW09} and its generalizations \cite{Ha10}, which have been also analysed in the framework of B-series \cite{Bseries}  (i.e., methods admitting a Taylor expansion with respect to the step-size), see e.g., \cite{HZ12}.

Further generalizations of HBVMs can be also found in \cite{BIT12_2,BIT12_3,BCMR12,BI12,BB12}.

\section{Discrete line integral methods}
The basic idea  which such methods rely on is straightforward. We shall first sketch it in the simplest case, as was done in \cite{IP07}, and then the argument will be generalized. Assume that, in problem (\ref{ham}), the Hamiltonian is a polynomial of degree $\nu$. Moreover, starting from the initial condition $y_0$ we want to produce a new approximation at $t=h$, say $y_1$, such that the Hamiltonian is conserved. By considering the simplest possible path joining $y_0$ and $y_1$, i.e., the segment
\begin{equation}\label{sigma_1}
\sigma(c h) = c y_1 + (1-c) y_0, \qquad c\in[0,1],
\end{equation}
one obtains:
\begin{eqnarray*}
H(y_1) - H(y_0) &=& H(\sigma(h)) -H(\sigma(0)) \\
&=& \int_0^h \nabla H(\sigma(t))^T\sigma'(t)\dd t \\
&=& h\int_0^1 \nabla H( \sigma(ch) )^T\sigma'(ch)\dd c\\
&=&h\int_0^1 \nabla H(cy_1+(1-c)y_0)^T(y_1-y_0)\dd c\\
&=& h \left[\int_0^1 \nabla H(cy_1+(1-c)y_0)\dd c\right]^T(y_1-y_0)~=~0,
\end{eqnarray*} 
provided that
\begin{equation}\label{y11}
y_1=y_0+hJ\int_0^1\nabla H(cy_1+(1-c)y_0)\dd c.
\end{equation}
In fact, due to the fact that $J$ is skew symmetric, one obtains:
\begin{eqnarray*}
\lefteqn{h^{-1}\left[\int_0^1 \nabla H(cy_1+(1-c)y_0)\dd c\right]^T(y_1-y_0)}\\
&=&\left[\int_0^1 \nabla H(cy_1+(1-c)y_0)\dd c\right]^TJ\left[\int_0^1 \nabla H(cy_1+(1-c)y_0)\dd c\right]~=~0.
\end{eqnarray*}
If $H\in\Pi_{\nu}$, then the integrand at the right-hand side in (\ref{y11}) has degree $\nu-1$ and, therefore, can be exactly computed by using, say, a Newton-Cotes formula based at $\nu$ equidistant abscissae in $[0,1]$. By setting, hereafter,
\begin{equation}\label{f}
f(\cdot) = J\nabla H(\cdot),
\end{equation}
one then obtains
\begin{equation}\label{y11d}
y_1 = y_0 + h \sum_{i=1}^\nu \bb_i f(c_iy_1+(1-c_i)y_0) \equiv y_0 + h \sum_{i=1}^\nu b_i f(Y_i)
\end{equation}
where 
\begin{equation}\label{y11d2}
c_i = \frac{i-1}{\nu-1}, \qquad Y_i = \sigma(c_i)\equiv c_iy_1+(1-c_i)y_0, \qquad i=1,\dots,\nu,
\end{equation}
and the $\{b_i\}$ are the quadrature weights:
$$b_i = \int_0^1 \prod_{j=1,\,j\ne i}^\nu\frac{t-c_j}{c_i-c_j}\dd t, \qquad i=1,\dots,\nu.$$

\subsection*{Some examples:}
\begin{itemize}
\item when $\nu=2$, one obtains the usual {\em trapezoidal rule}, $$y_1 = y_0+\frac{h}2\left( f(y_0) + f(y_1)\right)$$

\item when $\nu=3$, one obtains the following fomula:
$$y_1 = y_0+\frac{h}6\left( f(y_0) +4f\left(\frac{y_0+y_1}2\right) + f(y_1)\right)$$

\item when $\nu=5$, one obtains the formula:
$$y_1 = y_0+\frac{h}{90}\left( 7f(y_0) + 32f\left( \frac{3y_0+y_1}4\right)+12f\left( \frac{y_0+y_1}2\right)+32f\left( \frac{y_0+3y_1}4\right)+7f(y_1) \right).$$

\end{itemize}
The above fromulae were named $s$-stages trapezoidal rules in \cite{IP07}.  They provide exact conservation for polynomial Hamiltonian functions of degree no larger than $2\lceil\frac{\nu}2\rceil$, for all $\nu\ge1$. Their order of accuracy can be easily determined by recasting (\ref{y11d})--(\ref{y11d2}) as a $\nu$-stage Runge-Kutta method:
\begin{equation}\label{nustage}
\begin{array}{c|c}
\bfc & \bfc\bfb^T\\ \hline &\bfb^T
\end{array} \qquad\mbox{with}\quad \bfc = (c_1,\dots,c_\nu)^T \quad\mbox{and}\quad \bfb=(b_1,\dots,b_\nu)^T,
\end{equation}
which satisfies some of the usual {\em simplifying assumptions }(see, e.g., \cite[p.\,71]{HW96}) for an $s$-stage Runge-Kutta method (see (\ref{RK}) with coefficients $b_i,c_i,a_{ij}$, $i,j=1,\dots,s$:
\begin{description}
\item{$B(p)$:\quad}   $\sum_{i=1}^s b_i c_i^{q-1} = \frac{1}q, \qquad q=1,\dots,p$,

\item{$C(\eta)$:\quad} $\sum_{j=1}^s a_{ij} c_i^{q-1} = \frac{c_i^q}q, \qquad q=1,\dots,\eta, \quad i=1,\dots,s$,

\item{$D(\zeta)$:\quad} $\sum_{i=1}^s b_i c_i^{q-1}a_{ij} = \frac{b_j}q(1-c_j^q), \qquad q=1,\dots,\zeta, \quad j=1,\dots,s$.
\end{description}
In such a case, in fact, the following result holds true.

\begin{theo}[Butcher, 1964] If a Runge-Kutta method satisfies conditions $B(p)$, $C(\eta)$, and $D(\zeta)$, with
$$p\le\min\{ \eta+\zeta+1, ~2(\eta+1)\},$$ then it has order $p$.
\end{theo}
As a matter of fact, $B(2)$ and $C(1)$ turn out to be satisfied, for (\ref{nustage}), thus resulting in a second order method. In more details:
\begin{itemize} 
\item the quadrature is exact for polynomials of degree 1, so that $B(2)$ holds true;

\item moreover, by setting $\bfe=(1,\dots,1)^T$, one has
$$\bfc \bfb^T\bfe = \bfc \qquad\Leftrightarrow\qquad C(1).$$ 

\end{itemize}
\begin{rem} It is worth mentioning that, even though (\ref{nustage}) is formally a $\nu$-stage implicit Runge-Kutta method, nevertheless the actual size of the generated discrete problem consists of only {\em one} nonlinear equation, in the unknown $y_1$, as the above examples clearly show. The mono-implicit character of these methods comes from the fact that their coefficient matrix has rank one. 
\end{rem}

\section{Generalizing the approach}
The next step is to generalize the above approach, where we have assumed that  the  path $\sigma(ch)$, defined in (\ref{sigma_1}), is a linear function. Now we consider a polynomial path $\sigma$ of degree $s\ge1$.
Having fixed a suitable basis $\{P_0,\dots,P_{s-1}\}$ for $\Pi_{s-1}$, one can expand the derivative of $\sigma$ as
\begin{equation}\label{sigma1_s}
\sigma'(ch) = \sum_{j=0}^{s-1} P_j(c) \gamma_j, \qquad c\in[0,1],
\end{equation}
for certain set of coefficients $\{\gamma_j\}$ to be determined. By imposing the initial condition $$\sigma(0)=y_0,$$ one then formally obtains
\begin{equation}\label{sigma_s}
\sigma(ch) = y_0 + h\sum_{j=0}^{s-1} \int_0^c P_j(x)\dd x\,\gamma_j, \qquad c\in[0,1],
\end{equation} 
with the new approximation given by $y_1=\sigma(h)$. Energy conservation may be obtained  by following a similar computation as before, namely
\begin{eqnarray*}
H(y_1) - H(y_0) &=& H(\sigma(h)) -H(\sigma(0)) \\
&=& \int_0^h \nabla H(\sigma(t))^T\sigma'(t)\dd t \\
&=& h\int_0^1 \nabla H( \sigma(ch) )^T\sigma'(ch)\dd c\\
&=&h\int_0^1 \nabla H(\sigma(c h))^T\sum_{j=0}^{s-1} P_j(c) \gamma_j \dd c\\
&=&h\sum_{j=0}^{s-1} \left[\int_0^1 \nabla H(\sigma(c h))P_j(c)\dd c\right]^T \gamma_j~=~0,
\end{eqnarray*} 
provided that the unknown coefficients $\{\gamma_j\}$ satisfy
\begin{equation}\label{gammajs}
\gamma_j = \eta_j J\int_0^1\nabla H(\sigma(ch))P_j(c)\dd c, \qquad j=0,\dots,s,
\end{equation}
for a suitable set of nonzero scalars $\eta_0,\dots,\eta_{s-1}$.
The new approximation is then given by plugging (\ref{gammajs}) into (\ref{sigma_s}):
\begin{equation}\label{y1s}
y_1 = \sigma(h) = y_0 + h\sum_{j=0}^{s-1} \eta_j\int_0^1 P_j(x)\dd x \int_0^1 P_j(\tau)f(\sigma(\tau h))\dd\tau.
\end{equation}
As before, assume $H\in\Pi_\nu$. Then, the integrands in (\ref{gammajs}) and (\ref{y1s}) have at most degree $(\nu-1)s+\nu-1 \equiv \nu s-1$. Therefore, by fixing a suitable set of $k$ abscissae $0\le c_1<\dots<c_k\le1$, and corresponding quadrature weights $\{b_1,\dots,b_k\}$, such that the resulting quadrature formula is exact for polynomials of degree $\nu s - 1$, the integrals in (\ref{gammajs}) and (\ref{y1s}) may be replaced by the corresponding quadrature formulae, which yields
\begin{equation}\label{gammajsd}
\gamma_j = \eta_j \sum_{i=1}^kb_if(\sigma(c_ih))P_j(c_i), \qquad j=0,\dots,s,
\end{equation}
and
\begin{equation}\label{y1sd}
y_1 \equiv \sigma(h) = y_0 + h\sum_{j=0}^{s-1} \eta_j\int_0^1 P_j(x)\dd x \sum_{i=1}^k b_i P_j(c_i)f(\sigma(c_i h)),
\end{equation}
respectively. By setting, as before, $$Y_i=\sigma(c_ih), \qquad i=i,\dots,k,$$ one then obtains:
\begin{eqnarray}\label{aij}
Y_i &=& y_0 + h\sum_{j=1}^k \overbrace{\left[ b_j\sum_{\ell=0}^{s-1} \eta_\ell P_\ell(c_j)\int_0^{c_i} P_\ell(x)\dd x\right]}^{=\,a_{ij}} f(Y_j)\\ &\equiv& y_0 + h\sum_{j=1}^k a_{ij} f(Y_i), \qquad i=1,\dots,k, \nonumber \\
y_1&=&y_0+h\sum_{i=1}^k\underbrace{\left[b_i  \sum_{\ell=0}^{s-1} \eta_\ell P_\ell(c_i) \int_0^1 P_\ell(x)\dd x\right]}_{=\,\hat{b}_i} f(Y_i)\label{hbi} \\ \nonumber &\equiv& y_0+h\sum_{i=1}^k\hat{b}_i f(Y_i).
\end{eqnarray} We are then speaking of the following $k$-stage Runge-Kutta method:
\begin{equation}\label{kstage}
\begin{array}{c|c}
\bfc & A \equiv (a_{ij}) \in \RR^{k\times k}\\ \hline &\hat\bfb^T
\end{array} \qquad\mbox{with}\quad \bfc = (c_1,\dots,c_k)^T,\qquad \hat\bfb=(\hat{b}_1,\dots,\hat{b}_k)^T,
\end{equation}
with $a_{ij}$, $\hat{b}_i$ defined according to (\ref{aij}) and (\ref{hbi}), respectively.

In so doing, energy conservation can always be achieved, provided that the quadrature has a suitable high order. For example, we can place the $k$ abscissae $\{c_i\}$ at the $k$ Gauss-Legendre nodes on $[0,1]$ thus obtaining maximum order $2k$. In such a case, energy conservation is guaranteed for polynomial Hamiltonians of degree no larger that \begin{equation}\label{nu} \nu \le \frac{2k}s.\end{equation} 

However, it is quite difficult to discuss the order of accuracy and the properties of the $k$-stage Runge-Kutta method (\ref{kstage}), when a generic polynomial basis is considered. As matter of fact, different choices of the basis could provide different methods, having different orders. As an example, fourth-order energy-conserving Runge-Kutta methods were derived in \cite{IP08,IT09}, by using the Newton polynomial basis, defined at the abscissae $\{c_i\}$.
We shall see that things will greatly simplify, by choosing an orthonormal polynomial basis.

\begin{rem}\label{RKform0}
It is worth noticing that we can cast in matrix form the Butcher tableau of the $k$-stage Runge-Kutta method (\ref{kstage}) by introducing the matrices
$$\P_s =\pmatrix{ccc}
P_0(c_1) & \dots & P_{s-1}(c_1)\\
\vdots   &           &\vdots \\
P_0(c_k) &  \dots &P_{s-1}(c_k)\endpmatrix \in\RR^{k\times s},
$$
$$
\II_s = \pmatrix{ccc}
\int_0^{c_1}P_0(x)\dd x & \dots &\int_0^{c_1}P_{s-1}(x)\dd x\\
\vdots   &           &\vdots \\
\int_0^{c_k}P_0(x)\dd x &  \dots &\int_0^{c_k}P_{s-1}(x)\dd x\endpmatrix \in\RR^{k\times s},
$$
$$\Lambda_s = \pmatrix{ccc} \eta_1\\ &\ddots \\ &&\eta_s\endpmatrix\in\RR^{s\times s},$$
$$\Omega = \pmatrix{ccc} b_1\\ &\ddots \\ &&b_k\endpmatrix\in\RR^{k\times k},$$
and the row vector
\begin{equation}\label{Is1}
\II_s^1 = \pmatrix{ccc}
\int_0^1P_0(x)\dd x & \dots &\int_0^1P_{s-1}(x)\dd x\endpmatrix,
\end{equation}
as follows:
$$
\begin{array}{c|c}
\bfc & \II_s\Lambda_s\P_s^T\Omega\\
\hline
       & \II_s^1\Lambda_s \P_s^T\Omega
\end{array} 
$$
which will be further studied later.
\end{rem}

\chapter{Background results}

In this chapter, we state a few preliminary results concerning Legendre polynomials and differential equations, for later reference.
This chapter is based on references \cite{BIT09,BIT10_1}.

\section{Legendre polynomials}\label{legpol}
The following polynomials, denoted by $P_i$, are the {\em Legendre} polynomials shifted on the interval $[0,1]$, and scaled in order to be orthonormal:
\begin{equation}\label{orto}
\deg P_i = i, \qquad \int_0^1 P_i(x)P_j(x)\dd x = \delta_{ij}, \qquad \forall i,j\ge 0,
\end{equation}
where $\delta_{ij}$ is the Kronecker symbol (its value is 1, when $i=j$, and 0, otherwise). As any family of orthogonal polynomials, they satisfy a 3-terms recurrence, which is given, in this specific case, by:
\begin{eqnarray}\nonumber
&&P_0(x)\equiv 1, \qquad P_1(x)=\sqrt{3}(2x-1),\\ \label{Legendre}\\
&&P_{i+1}(x) = (2x-1)\frac{2i+1}{i+1}\sqrt{\frac{2i+3}{2i+1}}P_i(x)-\frac{i}{i+1}\sqrt{\frac{2i+3}{2i-1}}P_{i-1}(x), \qquad i\ge1.\nonumber
\end{eqnarray}
We recall that the roots $\{c_1,\dots,c_k\}$ of $P_k(x)$ are all distinct and belong to the interval $(0,1)$. Thus they may be identified via the following conditions:  
\begin{equation}\label{legnodes}
P_k(c_i)=0, \quad i=1,\dots,k, \qquad\mbox{with}\qquad 0<c_1<\dots<c_k<1.
\end{equation}
It is known that they are symmetrically distributed in the interval [0,1]:
\begin{equation}\label{symci}
c_i = 1-c_{k-i+1}, \qquad i=1,\dots,k.
\end{equation}
They are referred to as the Gauss-Legendre abscissae on $[0,1]$, and generate the Gauss-Legendre quadrature formula of order $2k$, namely an interpolating quadrature formula which is exact for polynomials of degree no larger than $2k-1$. In fact, if $p(x)\in\Pi_{2k-1}$, then it can be written as
$$p(x) = q(x)P_k(x)+r(x), \qquad q(x),r(x)\in\Pi_{k-1}.$$ Consequently,
\begin{eqnarray*}
\int_0^1p(x)\dd x &=& \int_0^1 \left[q(x)P_k(x)+r(x)\right]\dd x\\
&=& \underbrace{\int_0^1 q(x)P_k(x)\dd x}_{=\,0}+\int_0^1 r(x)\dd x ~=~ \int_0^1 r(x)\dd x,
\end{eqnarray*}
since $P_k(x)$ is orthogonal to polynomials of degree less than $k$ (see (\ref{orto})). On the other hand, for the quadrature formula $(c_i,b_i)$, with the quadrature weights given by
$$b_i = \int_0^1 \prod_{\stackrel{j=1}{j\not=i}}^{k} \frac{x-c_j}{c_i-c_j}\dd x, \qquad i=1,\dots,k,$$
one obtains:
$$
\sum_{i=1}^k b_ip(c_i) ~=~ \sum_{i=1}^k b_i\left[q(c_i)\overbrace{P_k(c_i)}^{=\,0}+r(c_i)\right]
~=~ \sum_{i=1}^k b_ir(c_i) ~=~ \int_0^1 r(x)\dd x,
$$
due to the fact that any quadrature based at $k$ distinct abscissae is exact for polynomials of degree no larger than $k-1$. As a matter of fact, for such a quadrature formula, for any function $f\in C^{2k}([0,1])$, one has
\begin{equation}\label{ek}
\int_0^1 f(x)\dd x = \sum_{i=1}^k b_i f(c_i) + \Delta_k, \qquad \Delta_k = \rho_k f^{(2k)}(\zeta),
\end{equation}
for a suitable $\zeta\in(0,1)$, and with $\rho_k$ independent of $f$. More in general, if the quadrature would have order $q\le2k$,  one would obtain
\begin{equation}\label{eq}
\int_0^1 f(x)\dd x = \sum_{i=1}^k b_i f(c_i) + \Delta_k, \qquad \Delta_k = \rho_k f^{(q)}(\zeta),
\end{equation}
with $\zeta$ and $\rho_k$ defined similarly as above, thus showing that the formula is exact for polynomials of degree no larger than $q-1$. 

In particular, in the sequel, we shall need to discuss the case where the integrand in (\ref{eq}) has the following form,
\begin{equation}\label{Pjf}
f(\tau) = P_j(\tau)f(\tau h), \qquad \tau\in[0,1],
\end{equation}
with $P_j$ the $j$th Legendre polynomial. The following result then holds true.
\begin{lem}\label{erroq}
Let $f\in C^{(q)}$, $q$ being the order of the given quadrature formula $(c_i,b_i)$ over the interval $[0,1]$. Then $$\int_0^1P_j(\tau)f(\tau h)\dd\tau-\sum_{i=1}^kb_i P_j(c_i)f(c_ih)= O(h^{q-j}), \qquad j=0,\dots,q.$$
\end{lem}
\proof
The thesis follows from (\ref{eq}), by considering that
\begin{eqnarray*}
\frac{\dd^q}{\dd \tau^q} P_j(\tau)f(\tau h) &\equiv& \left[ P_j(\tau) f(\tau h)\right]^{(q)} ~=~ \sum_{i=0}^q 
\pmatrix{c} q\\ i\endpmatrix P_j^{(i)}(\tau) f^{(q-i)}(\tau h) h^{q-i}\\
&=& \sum_{i=0}^j \pmatrix{c} q\\ i\endpmatrix P_j^{(i)}(\tau) f^{(q-i)}(\tau h) h^{q-i} ~=~ O(h^{q-j}),
\end{eqnarray*}
since $P_j^{(i)}(\tau)\equiv0$, for $i>j$.\,\QED
\bigskip

We also need a further result concerning integrals with integrands in the form (\ref{Pjf}), which is stated below.

\begin{lem}\label{hj}
Let $G:[0,h]\rightarrow V$, with $V$ a suitable vector space, a function which admits a Taylor expansion at 0. Then $$\int_0^1 P_j(\tau)G(\tau h)\dd\tau = O(h^j),\qquad j\ge0.$$
\end{lem}

\proof One obtains, by expanding $G(\tau h)$ at $\tau=0$:
\begin{eqnarray*}
\int_0^1P_j(t)G(\tau h)\dd\tau &=& \int_0^1P_j(t)\sum_{k\ge0}\frac{G^{(k)}(0)}{k!}(\tau h)^k\dd\tau ~=~\sum_{k\ge0}\frac{G^{(k)}(0)}{k!}h^k\int_0^1P_j(\tau)\tau^k\dd\tau\\
&=& \sum_{k\ge j}\frac{G^{(k)}(0)}{k!}h^k\int_0^1P_j(\tau)\tau^k\dd\tau ~=~O(h^j),
\end{eqnarray*}
where last but one equality follows from the fact that
$$\int_0^1P_j(\tau)\tau^k\dd\tau = 0, \qquad\mbox{for}\qquad k<j.\,\QED$$

\section{Matrices defined by the Legendre polynomials}
The integrals of the Legendre polynomials are related to the polynomial themselves as follows. For all $c\in[0,1]$:
\begin{eqnarray}\nonumber
\int_0^c P_0(x)\dd x &=& \xi_1 P_1(c)+\frac{1}2 P_0(c), \\ \label{xi} 
\int_0^c P_i(x)\dd x &=& \xi_{i+1} P_{i+1}(c)-\xi_i P_{i-1}(c), \qquad i\ge 1,\\
\xi_i &=& \frac{1}{2\sqrt{4i^2-1}}. \nonumber
\end{eqnarray}

\begin{rem}\label{remleg1} From the orthogonal conditions (\ref{orto}), and taking into account that $P_0(x)\equiv 1$, one obtains:
\begin{equation}\label{intPj}
\int_0^1 P_0(x)\dd x = 1, \qquad \int_0^1 P_j(x)\dd x = 0, \qquad \forall j\ge1.
\end{equation}
\end{rem}
Legendre polynomials possess the following symmetry property:
\begin{equation}\label{symPj}
P_j(c) = (-1)^j P_j(1-c), \qquad c\in[0,1],\qquad j\ge0.
\end{equation}
Consequently, their integrals share a similar symmetry:
\begin{equation}\label{symIj}
\int_{c_1}^{c_2}P_j(x)\dd x = (-1)^j \int_{1-c_2}^{1-c_1}P_j(x)\dd x, \qquad c_1,c_2\in[0,1],\qquad j\ge0.
\end{equation}

In the sequel, we shall use the following matrices, defined by means of the Legendre polynomials evaluated at the $k\ge s$ abscissae (\ref{legnodes}):\footnote{They have been formally introduced, for a generic polynomial basis, at the end of the previous chapter>}
\begin{equation}\label{Ps}
\P_s =\pmatrix{ccc}
P_0(c_1) & \dots & P_{s-1}(c_1)\\
\vdots   &           &\vdots \\
P_0(c_k) &  \dots &P_{s-1}(c_k)\endpmatrix \in\RR^{k\times s},
\end{equation}
and
\begin{equation}\label{Is}
\II_s = \pmatrix{ccc}
\int_0^{c_1}P_0(x)\dd x & \dots &\int_0^{c_1}P_{s-1}(x)\dd x\\
\vdots   &           &\vdots \\
\int_0^{c_k}P_0(x)\dd x &  \dots &\int_0^{c_k}P_{s-1}(x)\dd x\endpmatrix \in\RR^{k\times s}.
\end{equation}
Because of the (\ref{xi}), they are related by the following relation:
\begin{equation}\label{Xs}
\II_s = \P_{s+1}\hat{X}_s, \qquad \hat{X}_s = \pmatrix{cccc}
\frac{1}2 &-\xi_1\\
\xi_1      &0   &\ddots\\
              &\ddots &\ddots &-\xi_{s-1}\\
              &           &\xi_{s-1} &0\\
                            \hline
              &           &               &\xi_s\endpmatrix \equiv \pmatrix{c} X_s\\\hline  0\,\dots\,0~\xi_s\endpmatrix.
\end{equation}
We also set
\begin{equation}\label{Omega}
\Omega = \pmatrix{ccc} b_1\\ &\ddots \\ &&b_k\endpmatrix\in\RR^{k\times k}
\end{equation}
the diagonal matrix with the corresponding Gauss-Legendre weights.  
The following simple properties then hold true.

\begin{theo}\label{Pkeqs} Matrices (\ref{Ps}) and (\ref{Is}) have full column rank, for all $s=1,\dots,k$. Moreover,
\begin{equation}\label{orto1}\P_s^T\Omega \P_{s+1} =\pmatrix{cc} I_s &\bfo\endpmatrix.\end{equation}
\end{theo}

\proof By considering any set of $s\le k$ rows of $\P_s$,  the resulting sub-matrix is the Gram matrix of the $s$ linearly independent polynomials $P_0,\dots,P_{s-1}$ defined at the corresponding $s$ (distinct) abscissae. It is, therefore, nonsingular and, then, $\P_s$ has full column rank. Moreover, when $s=k+1$, one has
\begin{equation}\label{kpiu1}
\P_{k+1} = \pmatrix{cc} \P_k &\bfo\endpmatrix,
\end{equation}
since the last column contains $P_k(c_i)=0$,  $i=1,\dots,k$.
As a consequence, because of (\ref{Xs}), for matrix $\II_s$ one obtains:
\begin{itemize}
\item when $s<k$, then both $\P_{s+1}$ and $\hat{X}_s$ have full column rank and so has $\II_s$;
\item when $s=k$, then from (\ref{kpiu1}) it follows that $$\II_k =  \P_{k+1} \hat{X}_k= \P_kX_k,$$ and both $\P_k$ and $X_k$ are nonsingular.\footnote{Indeed, it can be proved that $X_k$ is nonsingular for all $k\ge1$.}
\end{itemize}
Concerning (\ref{orto1}), one has, by considering that the quadrature formula $(c_i,b_i)$ is exact for polynomials of degree no larger that $2k-1\ge 2s-1$, and setting $\bfe_i\in\RR^s$ and $\hat\bfe_j\in\RR^{s+1}$ the $i$th and $j$th unit vectors:
\begin{eqnarray*}
\bfe_i^T \P_s^T\Omega \P_{s+1} \hat\bfe_j &=& \sum_{\ell=1}^k b_\ell P_{i-1}(c_\ell)P_{j-1}(c_\ell) 
~=~ \int_0^1 P_{i-1}(x)P_{j-1}(x)\dd x ~=~ \delta_{ij},\\
&& \qquad\qquad\qquad \forall i=1,\dots,s, ~\mbox{and}~j=1,\dots,s+1.\,\QED\end{eqnarray*}

From the previous theorem, the following result easily follows.
\begin{cor}\label{Psmeno1}
When $k=s$, then $\P_s^{-1}=\P_s^T\Omega$.
\end{cor}

\section{Additional preliminary results}
In order to deal with the analysis of the methods, we need the following perturbation result concerning the solution of the initial value problem for ordinary differential equations:
\begin{equation}\label{yt0y0}
y'(t) = f(y(t)), \qquad t\ge t_0, \qquad y(t_0) = y_0,
\end{equation}
whose solution will be denoted by $y(t;t_0,y_0)$, in order to emphasize the dependence on the initial condition, set at $(t_0,y_0)$.  

Associated with this problem is the corresponding {\em fundamental matrix}, $\Phi(t,t_0)$, satisfying the {\em variational problem} (see also (\ref{variational}))
$$\Phi'(t,t_0) = J_f(y(t;t_0,y_0))\Phi(t,t_0), \qquad t\ge t_0, \qquad \Phi(t_0,t_0)=I,$$
where the derivative (i.e., $'$) is with respect to $t$, and $J_f(y)$ is the Jacobian matrix of $f(y)$.
The following result then holds true.

\begin{lem}\label{initper} With reference to the solution $y(t;t_0,y_0)$ of problem (\ref{yt0y0}), one has:
$$(i)\quad \frac{\partial}{\partial y_0} y(t;t_0,y_0) = \Phi(t,t_0);\qquad\qquad (ii)\quad\frac{\partial}{\partial t_0} y(t;t_0,y_0) = -\Phi(t,t_0)f(y_0).$$
\end{lem}

\proof Let us consider a perturbation $\delta y_0$ of the initial condition, and let $y(t;t_0,y_0+\delta y_0)$ be the corresponding solution. Consequently, 
\begin{eqnarray*}
y'(t;t_0,y_0+\delta y_0) &=& f(y(t;t_0,y_0+\delta y_0)) \\[2mm]
&=& \underbrace{f(y(t;t_0,y_0))}_{=\,y'(t;t_0,y_0)} + J_f(y(t;t_0,y_0))\left[ y(t;t_0,y_0+\delta y_0) - y(t;t_0,y_0)\right]\\[-4mm]
&&\qquad\qquad\qquad\qquad +~ O\left(\left\| y(t;t_0,y_0+\delta y_0) - y(t;t_0,y_0)\right\|^2\right).
\end{eqnarray*}
Therefore, by setting $$z(t)=y(t;t_0,y_0+\delta y_0) - y(t;t_0,y_0),$$ one obtains that, at first order (as is the case, when we let $\delta y_0\rightarrow0$),
$$z'(t) = J_f(y(t;t_0,y_0)) z(t), \qquad z(t_0) = \delta y_0.$$
The solution of this linear problem is easily seen to be 
$$z(t) = \Phi(t,t_0)\delta y_0$$ and, consequently, 
$$\frac{\partial}{\partial y_0} y(t;t_0,y_0) = \frac{\partial}{\partial (\delta y_0)} z(t) = \Phi(t,t_0),$$ i.e., the part $(i)$ of the thesis. 

Concerning the part $(ii)$, let consider a scalar $\eps\approx 0$ and observe that, by setting,  $y(t)=y(t;t_0,y_0)$, then $$y(t;t_0+\eps,y_0)\equiv y(t-\eps).$$ Consequently, at first order, the solution of the perturbed problem
\begin{equation}\label{yt0y0e}
y'(t) = f(y(t)), \qquad t\ge t_0+\eps, \qquad y(t_0+\eps) = y_0,
\end{equation}
coincides with that of the problem
\begin{equation}\label{yt0y0e1}
y'(t) = f(y(t)), \qquad t\ge t_0, \qquad y(t_0) = y_0(\eps)\equiv y_0-\eps f(y_0).
\end{equation}
Letting $\eps\rightarrow0$, one then obtains:
$$\frac{\partial}{\partial t_0} y(t;t_0,y_0) ~=~ \underbrace{\frac{\partial}{\partial y_0} y(t;t_0,y_0)}_{=\,\Phi(t,t_0)}\,\overbrace{ \frac{\partial}{\partial \eps} y_0(\eps)}^{=\,-f(y_0)} ~=~ -\Phi(t,t_0)f(y_0).$$ This concludes the proof.\,\QED\bigskip

\chapter{A Framework for HBVMs}
In this chapter, we provide a novel framework for discussing the order and conservation properties of  HBVMs, based on a local Fourier expansion of the vector field defining the dynamical system. In particular, this approach allows us to easily discuss the linear stability properties of the methods.
The material in this chapter is based on \cite{BIT10_1,BIT12_1,BIT11}. It is worth noticing that an interesting extension of this approach has been recently proposed in \cite{TS12}.

\section{Local Fourier expansion}
Legendre polynomials, previously introduced, constitute an {\em orthonomal basis} for the functions defined on the interval $[0,1]$. Therefore, we can formally expand the second member of (\ref{ham}) over the interval $[0,h]$ as follows (we use the notation (\ref{f})):
\begin{equation}\label{four1}
f(y(ch)) = \sum_{j\ge0} P_j(c) \gamma_j(y), \qquad c\in[0,1],
\end{equation}
where  
\begin{equation}\label{gjy}
\gamma_j(y) = \int_0^1 P_j(\tau)f(y(\tau h))\dd\tau, \qquad j\ge0.
\end{equation}
The expansion (\ref{four1})-(\ref{gjy}) is known as the {\em Neumann expansion} of an analytic function,\footnote{E.T.\,Whittaker, G.N.\,Watson, {\em A Course in Modern Analysis, Fourth edition}, Cambridge University Press, 1950, page 322.} and converges uniformly, provided that the function $g(c) = f(y(ch))$ has continuous second derivative:\footnote{E.\,Isaacson, H.B.\,Keller, {\em Analysis of Numerical Methods}, Wiley \& Sons, 1966, page 206.} for sake of simplicity, hereafter we shall assume $g(t)$ to be analytic.

In so doing, we are transforming the initial value problem
\begin{equation}\label{ivp1}
y'(t) = f(y(t)), \qquad t\in[0,h], \qquad y(0)=y_0,
\end{equation}
into the equivalent {\em integro-differential} problem
\begin{equation}\label{idp1}
y'(ch) = \sum_{j\ge0} P_j(c) \int_0^1 P_j(\tau)f(y(\tau h))\dd\tau, \qquad c\in[0,1], \qquad y(0)=y_0.
\end{equation}
In order to obtain a polynomial approximation of degree $s$ to (\ref{ivp1})-(\ref{idp1}), we just truncate the infinite series to a finite sum. The resulting initial value problem is
\begin{eqnarray}\nonumber
\sigma'(ch) &=& \sum_{j=0}^{s-1} P_j(c) \int_0^1 P_j(\tau)f(\sigma(\tau h))\dd\tau ~\equiv~ \sum_{j=0}^{s-1} P_j(c) \gamma_j(\sigma), \qquad c\in[0,1], \\
\sigma(0) &=&y_0, \label{idp2}
\end{eqnarray}
and its solution evidently defines a polynomial  $\sigma\in \Pi_{s}$. Integrating both sides of (\ref{idp2}) yields the equivalent formulation  
\begin{equation}\label{sigmac}
\sigma(ch) = y_0 +h\sum_{j=0}^{s-1} \int_0^cP_j(x)\dd x\, \gamma_j(\sigma), \qquad c\in[0,1],
\end{equation}
where the notation (\ref{gjy}) has been used again. One easily recognizes that (\ref{idp2}) defines the very same expansion 
(\ref{sigma1_s})--(\ref{gammajs}) seen before (with all $\eta_j=1$). Consequently, such a method is energy-conserving, if we are able to exactly compute the integrals providing the coefficients $\gamma_j(\sigma)$ at the right-hand side in (\ref{idp2}). From Remark~\ref{remleg1}, one obtains that
\begin{equation}\label{sigmah}
\sigma(h) = y_0 + \int_0^hf(\sigma(\tau h))\dd\tau.
\end{equation}

Let now discuss the order of the approximation $\sigma(h)\approx y(h)$.  
\begin{lem}\label{gamj} Let $\gamma_j(\sigma)$ be defined according to (\ref{gjy}). Then $\gamma_j(\sigma)=O(h^j)$.
\end{lem}
\proof The proof follows immediately from (\ref{gjy}), by vrtue of Lemma~\ref{hj}.\,\QED\bigskip

We are now able to prove the following result.
\begin{theo}\label{ord1} $\sigma(h)-y(h)=O(h^{2s+1})$.\end{theo}

\proof Denoting by $y(t;t_0,y_0)$ the solution of problem (\ref{yt0y0}) and considering that $\sigma(0)=y_0$, by virtue of Lemmas~\ref{hj} and \ref{gamj}  one has:
\begin{eqnarray*}
\sigma(h)-y(h)&=& y(h;h,\sigma(h))-y(h;0,y_0) \\ 
&\equiv& y(h;h,\sigma(h))-y(h;0,\sigma(0)) ~=~ \int_0^h \frac{\dd}{\dd t} y(h;t,\sigma(t))\dd t\\
&=& \int_0^h \left(\frac{\partial}{\partial \theta} y(h;\theta,\sigma(t))\Big|_{\theta=t} + \frac{\partial}{\partial \omega}y(h;t,\omega)\Big|_{\omega=\sigma(t)}\sigma'(t)\right)\dd t\\
&=&\int_0^h \left[-\Phi(h,t)f(\sigma(t)) +\Phi(t,t_0)\sigma'(t)\right]\dd t \\
&=& \int_0^h \Phi(h,t)[-f(\sigma(t)) +\sigma'(t)]\dd t \\
&=& h\int_0^1 \Phi(h,\tau h)[-f(\sigma(\tau h)) +\sigma'(\tau h)]\dd\tau \\
&=&-h\int_0^1 \Phi(h,\tau h)\left[ \sum_{j\ge0}P_j(\tau)\gamma_j(\sigma) - \sum_{j=0}^{s-1}P_j(\tau)\gamma_j(\sigma)\right]\dd\tau\\
&=&-h\int_0^1 \Phi(h,\tau h)\sum_{j\ge s}P_j(\tau)\gamma_j(\sigma)\dd\tau\\
&=&
-h\sum_{j\ge s}\underbrace{\left[\int_0^1 \overbrace{\Phi(h,\tau h)}^{\equiv G(\tau h)}\, P_j(\tau)\dd\tau\right]}_{=\,O(h^j)}\overbrace{\gamma_j(\sigma)}^{=\,O(h^j)}\\
&=& h\sum_{j\ge s} O(h^{2j})~=~O(h^{2s+1}).\,\QED
\end{eqnarray*}
We observe that, unless we can compute exactly the integrals defining the $\{\gamma_j(\sigma)\}$ in (\ref{idp2}) (which is the case, for example, when $f$ is a polynomial or in very special situations), (\ref{idp2}) is {\em not yet} an operative method, but rather a {\em formula}. In order to obtain a numerical approximation procedure, we need to approximate those integrals by means of a suitable quadrature formula, which we define at $k\ge s$ Gauss abscissae in $[0,1]$ defined in (\ref{legnodes}).
As was seen in Section~\ref{legpol}, the corresponding quadrature formula $(c_i,b_i)$ has order $q=2k$, that is it is exact for polynomials of degree no larger that $2k-1$.  By recalling Lemma~\ref{erroq}, let then approximate the integrals in (\ref{idp2}) by means of a quadrature $(c_i,b_i)$ over $k$ distinct abscissae. Consequently, in place of $\sigma$ defined by (\ref{idp2}) or (\ref{sigmac}), we shall compute the new polynomial $u\in\Pi_s$ such that:
\begin{eqnarray}\nonumber
u'(ch) &=& \sum_{j=0}^{s-1} P_j(c) \sum_{\ell=1}^k b_\ell P_j(c_\ell)f(u(c_\ell h)), \qquad c\in[0,1], \\
u(0) &=&y_0, \label{idp3}
\end{eqnarray}
that is,
\begin{equation}\label{uc}
u(ch) = y_0 +h\sum_{j=0}^{s-1}\int_0^c P_j(x)\dd x \sum_{\ell=1}^k b_\ell P_j(c_\ell)f(u(c_\ell h)), \qquad c\in[0,1],
\end{equation}
with the new approximation given by:
\begin{equation}\label{y1equh}
y_1 \equiv u(h) = y_0 + h\sum_{i=1}^k b_i f(u(c_ih)).
\end{equation}
If the quadrature formula $(c_i,b_i)$ has order $q$, then, by virtue of Lemma~\ref{erroq} and taking into account (\ref{gjy}), one obtains
\begin{eqnarray}\label{gju}
\gamma_j(u) &\equiv& \int_0^1 P_j(\tau)f(u(\tau h))\dd\tau ~=~ \sum_{\ell=1}^k b_\ell P_j(c_\ell)f(u(c_\ell h)) ~-~\Delta_j(h),\\ \nonumber
\Delta_j(h)&=&O(h^{q-j}), \quad j=0,\dots,q.\end{eqnarray}
Consequently, we can rewrite the first equation in (\ref{idp3}) in the following equivalent form:
\begin{equation}\label{idp4}
u'(ch) = \sum_{j=0}^{s-1} P_j(c) \left[ \gamma_j(u)-\Delta_j(h)\right], \qquad c\in[0,1].
\end{equation}
This allows us to derive the following result, which we state for a generic quadrature of order $q$.

\begin{theo}\label{ord2} $y_1-y(h)=O(h^{p+1})$, where $p=\min\{q,2s\}$.\end{theo}

\proof The proof proceeds on the same line as that of Theorem~\ref{ord1}:
\begin{eqnarray*}
y_1-y(h) &=&u(h)-y(h) ~=~ y(h;h,u(h))-y(h;0,u(0))\\
&=& \int_0^h \frac{\dd}{\dd t} y(h;t,u(t))\dd t ~=~ \int_0^h \left(\frac{\partial}{\partial \theta}y(h;\theta,u(t))\Big|_{\theta=t} + \frac{\partial}{\partial \omega}y(h;t,\omega)\Big|_{\omega=u(t)}u'(t)\right)\dd t\\
&=& \int_0^h \Phi(h,t)[-f(u(t)) +u'(t)]\dd t ~=~ h\int_0^1 \Phi(h,\tau h)[-f(u(\tau h)) +u'(\tau h)]\dd\tau \\
&=&-h\int_0^1 \Phi(h,\tau h)\left[ \sum_{j\ge0}P_j(\tau)\gamma_j(u) - \sum_{j=0}^{s-1}P_j(\tau)\left(\gamma_j(u)-\Delta_j(h)\right)\right]\dd\tau\\
&=&h\int_0^1 \Phi(h,\tau h)\sum_{j=0}^{s-1}P_j(\tau)\Delta_j(u)\dd\tau-h\int_0^1 \Phi(h,\tau h)\sum_{j\ge s}P_j(\tau)\gamma_j(u)\dd\tau\\
&=&
h\sum_{j=0}^{s-1}\underbrace{\left[\int_0^1 \overbrace{\Phi(h,\tau h)}^{\equiv G(\tau h)}\, P_j(\tau)\dd\tau\right]}_{=\,O(h^j)}\overbrace{\Delta_j(u)}^{=\,O(h^{q-j})}
-h\sum_{j\ge s}\underbrace{\left[\int_0^1 \overbrace{\Phi(h,\tau h)}^{\equiv G(\tau h)}\, P_j(\tau)\dd\tau\right]}_{=\,O(h^j)}\overbrace{\gamma_j(u)}^{=\,O(h^j)}\\
&=& O(h^{q+1}) + h\sum_{j\ge s} O(h^{2j})~=~O(h^{p+1}), \qquad p=\min\{q,2s\}.\,\QED
\end{eqnarray*}

\begin{defi}
The method (\ref{idp4})--(\ref{y1equh}) defines a {\em Hamiltonian Boundary Value Method (HBVM) with $k$ stages and degree $s$}, in short HBVM$(k,s)$.\end{defi} As a consequence, by setting the abscissae at the $k$ Gauss points (\ref{legnodes}), the following result holds true.

\begin{cor}\label{ordhbvm} By choosing the $k$ abscissae $\{c_i\}$ as in (\ref{legnodes}), a HBVM$(k,s)$ method has order $2s$, for all $k\ge s$.\end{cor}

\section{Runge-Kutta form of HBVM$(k,s)$} Before studying the conservation properties of the methods, let us derive the Runge-Kutta formulation of HBVM$(k,s)$. The basic fact is that, at the right-hand sides of equations (\ref{uc})--(\ref{y1equh}), one only requires to know the value of the polynomial $u$ at the abscissae $\{c_ih\}$. Consequently, by setting $$Y_i = u(c_ih), \qquad i=1,\dots,k,$$ one obtains:
\begin{eqnarray}\nonumber
Y_i &=& y_0 + h\sum_{j=1}^k \overbrace{\left[ b_j\sum_{\ell=0}^{s-1} P_\ell(c_j)\int_0^{c_i} P_\ell(x)\dd x \right]}^{\equiv\, a_{ij}}f(Y_j) \\ &\equiv&y_0+h\sum_{j=1}^k a_{ij}f(Y_j), \qquad i=1,\dots,k,\label{RKform}\\
y_1 &=& y_0+h\sum_{i=1}^k b_i f(Y_i). \nonumber
\end{eqnarray}
In other words, we have defined the following $k$-stage Runge-Kutta method:
\begin{equation}\label{kstage}
\begin{array}{c|c}
\bfc & A\equiv(a_{ij}) \\ \hline &\bfb^T
\end{array}
\end{equation}
with (see (\ref{RKform})),
$$\bfc = (c_1,\dots,c_k)^T, \quad \bfb=(b_1,\dots,b_k)^T, \mbox{\quad and \quad}A=(a_{ij})\in\RR^{k\times k}.$$
The Butcher tableau (\ref{kstage}) defines the {\em Runge-Kutta shape of a HBVM$(k,s)$ method}.
We can easily derive a more compact form for the Butcher array $A$ in (\ref{kstage}).

\begin{theo}\label{AeqPsIsO} $A=\II_s\P_s^T\Omega$, with the matrices $\II_s,\P_s,\Omega$ defined according to (\ref{Ps})--(\ref{Omega}).\end{theo}

\proof By setting $\bfe_i,\bfe_j\in\RR^k$ the $i$-th and $j$-th unit vectors, one obtains:
\begin{eqnarray*}
\bfe_i^T\II_s \P_s^T\Omega\bfe_j &=& \pmatrix{ccc} \int_0^{c_i} P_0(x)\dd x&\dots \int_0^{c_i} P_{s-1}(x)\dd x\endpmatrix \pmatrix{c}P_0(c_j)\\ \vdots \\ P_{s-1}(c_j)  \endpmatrix b_j\\
&=& b_j \sum_{\ell=0}^{s-1} P_\ell(c_j)\int_0^{c_i} P_\ell(x)\dd x ~\equiv~a_{ij} ~=~\bfe_i^TA\bfe_j,
\end{eqnarray*}
according to (\ref{RKform}).\,\QED
\bigskip

\no Consequently, the Buther tableau (\ref{kstage}) can be casted as:
\begin{equation}\label{kstage1}
\begin{array}{c|c}
\bfc & \II_s\P_s^T\Omega \\ \hline &\bfb^T
\end{array}
\end{equation}
or, equivalently, by taking into account (\ref{Xs}),
\begin{equation}\label{kstage2}
\begin{array}{c|c}
\bfc & \P_{s+1} \hat{X}_s\P_s^T\Omega \\ \hline &\bfb^T
\end{array}\,.
\end{equation}

\begin{rem} We observe that the Runge-Kutta form (\ref{kstage1}) of a HBVM$(k,s)$ method is simplified, with respect to that sketched in Remark~\ref{RKform0} for a discrete line-integral methods defined by using a general polynomial basis. In particular, the diagonal matrix $\Lambda_s$ is now authomatically fixed, in order to maximize the accuracy of the method, and the vector of the quadrature coincide with that used for approximating the integrals involved in the coefficients of the polynomial $u$.
\end{rem}

\subsection{HBVM$(s,s)$}\label{HBVMss}
In the case $k=s$, the matrices $\II_s,\P_s,\Omega\in\RR^{s\times s}$. Moreover, the following results follows immediately from 
Theorem~\ref{Pkeqs} and Corollary~\ref{Psmeno1}:
$$\II_s = \P_sX_s, \qquad \P_s^T\Omega = \P_s^{-1}.$$
Consequently, in such a case, we can write the Butcher tableau (\ref{kstage2}) as that of the following $s$-stage method,
\begin{equation}\label{sstage}
\begin{array}{c|c}
\bfc & \P_s X_s\P_s^{-1} \\ \hline &\bfb^T
\end{array}
\end{equation}
which is the {\em $W$-transform} defining the $s$-stage Gauss-Legendre Runge-Kutta collocation method \cite[p.\,79]{HW96}, which has order $2s$. In this sense, in the case $k\ge s$, HBVM$(k,s)$ can be ragarded as {\em low-rank generalizations} of the $s$-stage Gauss method. Indeed, the following result holds true.

\begin{theo}\label{ranks}
For all $k\ge s$ the rank of the matrix $A=\P_{s+1}\hat{X}_s\P_s^T\Omega$ is $s$. Moreover, the nonzero eigenvalues coincides with those of the basic $s$-stage Gauss method.
\end{theo}

\proof The rank of the matrix $\P_{s+1}$ is $s$ or $s+1$ (when $k>s$), whereas that of matrices $\hat{X_s},\P_s$ is $s$, and $\Omega$ is nonsingular. Therefore, the rank of $A$ cannot exceed $s$. Moreover, from Theorem~\ref{Pkeqs}, one has 
$$
\P_s^T\Omega A \P_s ~=~\P_s^T\Omega \P_{s+1}\hat{X}_s\P_s^T\Omega\P_s
~=~ (I_s~\bfo)\hat{X}_s I_s ~=~X_s\in\RR^{s\times s},
$$
which is known to be nonsingular. Consequently, $\rank(A)=s$. Moreover,
$$\P_s^T\Omega A = \P_s^T\Omega \P_{s+1}\hat{X}_s\P_s^T\Omega = (I_s~\bfo)\hat{X}_s \P_s^T\Omega = X_s\P_s^T\Omega.$$
This means that the columns of $\Omega\P_s$ span an $s$-dimensional left invariant subspace of $A$. Therefore, the eigenvalues of $X_s$ will coincide with the nonzero eigenvalues of $A$. On the other hand, from (\ref{sstage}) one obtains immediately that the eigenvalues of $X_s$ are the eigenvalues of the Butcher matrix of the $s$-stage Gauss method.\,\QED\bigskip

\no This property has been named {\em isospectrality of HBVMs}, in \cite{BIT12}. It will be used for the efficient implementation of HBVM$(k,s)$ methods.

\section{Energy conservation} We now consider the issue of energy conservation for HBVM$(k,s)$ methods. From (\ref{idp3})--(\ref{y1equh}) with $f=J\nabla H$, we obtain:
\begin{eqnarray*}
H(y_1)-H(y_0)&=& H(u(h))-H(u(0)) ~=~\int_0^h \nabla H(u(t))^Tu'(t)\dd t\\
&=& h\int_0^1 \nabla H(u(\tau h))^Tu'(\tau h)\dd\tau \\
&=& h\int_0^1 \nabla H(u(\tau h))^T\sum_{j=0}^{s-1}P_j(\tau)\sum_{i=1}^k b_i P_j(c_i)J\nabla H(c_i h)\dd\tau\\
&=& h\sum_{j=0}^{s-1}\left[\underbrace{\int_0^1 P_j(\tau)J\nabla H(u(\tau h))\dd\tau}_{=\,O(h^j)}\right]^TJ\left[\sum_{i=1}^k b_i P_j(c_i)J\nabla H(c_i h)\right]\\
&\equiv& E_H
\end{eqnarray*}
Now, two possibilities may occur:
\begin{itemize}

\item $\displaystyle \int_0^1 P_j(\tau)J\nabla H(u(\tau h))\dd\tau=\sum_{i=1}^k b_i P_j(c_i)J\nabla H(c_i h)$\,: in such case, $E_H=0$, so that energy is {\em exactly conserved}. This is the case of a polynomial Hamiltonian of degree $\nu$ no larger that $2k/s$;

\item $\displaystyle \int_0^1 P_j(\tau)J\nabla H(u(\tau h))\dd\tau=\sum_{i=1}^k b_i P_j(c_i)J\nabla H(c_i h) - \Delta_j(h)$\,: in such a case, by taking into account (\ref{gju}), one obtains that $E_H=O(h^{2k+1})$, provided that the Hamiltonian is suitably regular, as we have assumed.
\end{itemize}
We have then proved the following result.

\begin{theo} HBVM$(k,s)$ is energy-conserving for all polynomial Hamiltonian of degree 
\begin{equation}\label{nule2ks}
\nu\le \frac{2k}s.
\end{equation}
In any other case, $H(y_1)-H(y_0)=O(h^{2k+1})$, even though the method has order $s$.\end{theo}

\begin{figure}
\centerline{\includegraphics[width=12cm,height=8cm]{./figure/k2s2.png}}
\caption{Numerical level curves for problem (\ref{initex1})--(\ref{parex1}), 2-stage Gauss method, $h=10^{-3}$.}
\label{k2s2}

\bigskip
\bigskip
\centerline{\includegraphics[width=12cm,height=8cm]{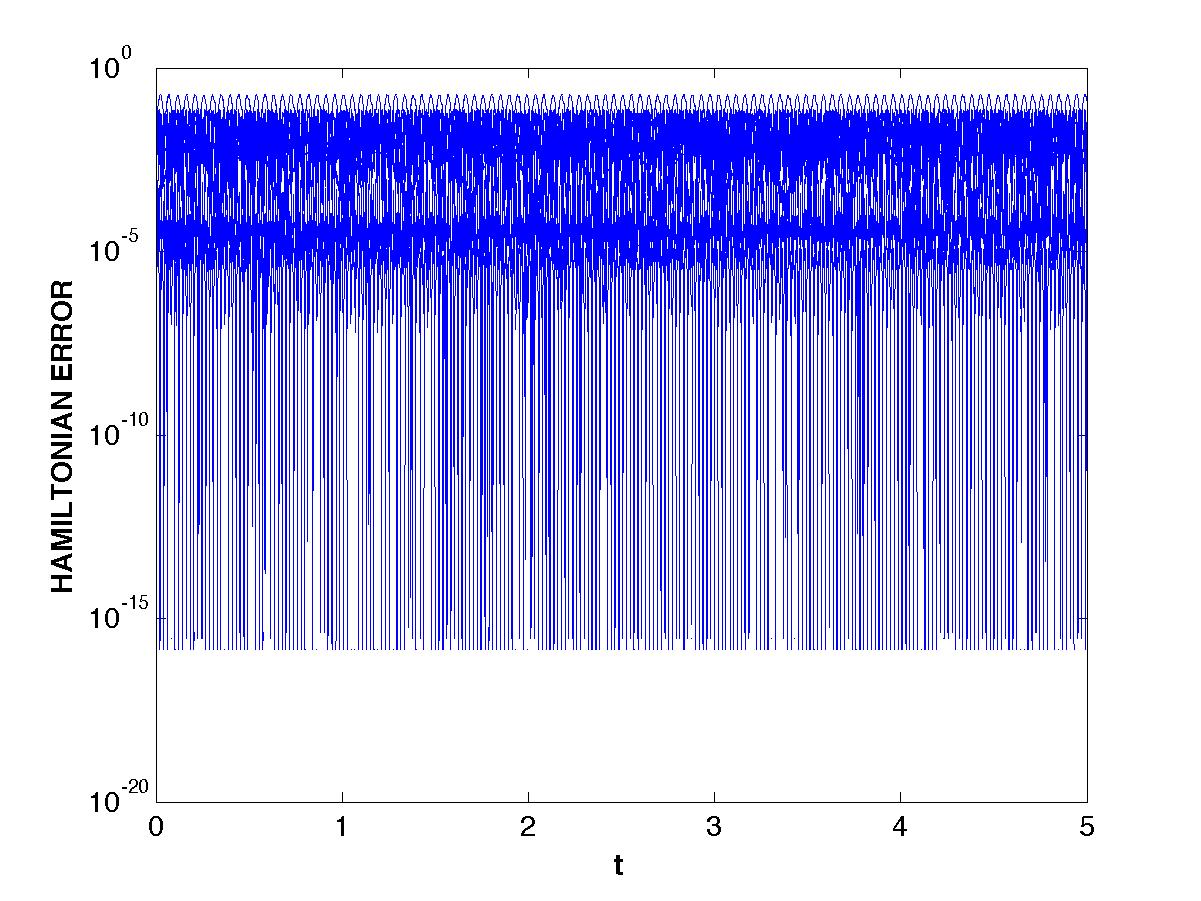}}
\caption{Hamiltonian error for problem (\ref{initex1})--(\ref{parex1}), 2-stage Gauss method, $h=10^{-3}$.}
\label{k2s2H}
\end{figure}

\begin{rem}
We observe that:
\begin{itemize}
\item for polynomial Hamiltonians, energy conservation can be {\em always} obtained, by choosing $k$ large enough, by virtue of (\ref{nule2ks});

\item even in the case of non polynomial Hamiltonians, energy conservation can be {\em practically} gained by choosing $k$ large enough, provided that $|E_{H}|$, which is $O(h^{2k+1})$, is within roundoff errors.
\end{itemize}
\end{rem}
As an example, in Figure~\ref{level} we plotted the level curves passing at
\begin{equation}\label{initex1} (q_0,p_0) = (i,-i), \qquad i = 1,\dots,8,\end{equation}
 for the Hamiltonian problem with Hamiltonian
\begin{equation}\label{Hex1}H(q,p) = p^2 + (\beta q)^2 +\alpha (q+p)^{2n},\end{equation}
with parameters:
\begin{equation}\label{parex1}
\beta = 10, \qquad \alpha = 1, \qquad n = 4.
\end{equation}
By using the 2-stage Gauss method (fourth-order), with step size $h=10^{-3}$, the obtained phase portrait is wrong, as is shown in Figure~\ref{k2s2}, due to the error in the numerical Hamiltonian, which is shown in Figure~\ref{k2s2H}. Indeed, even though no drift in the Hamiltonian occurs, nevertheless it is not negligible, for the problem at hand.

\begin{figure}
\centerline{\includegraphics[width=12cm,height=8cm]{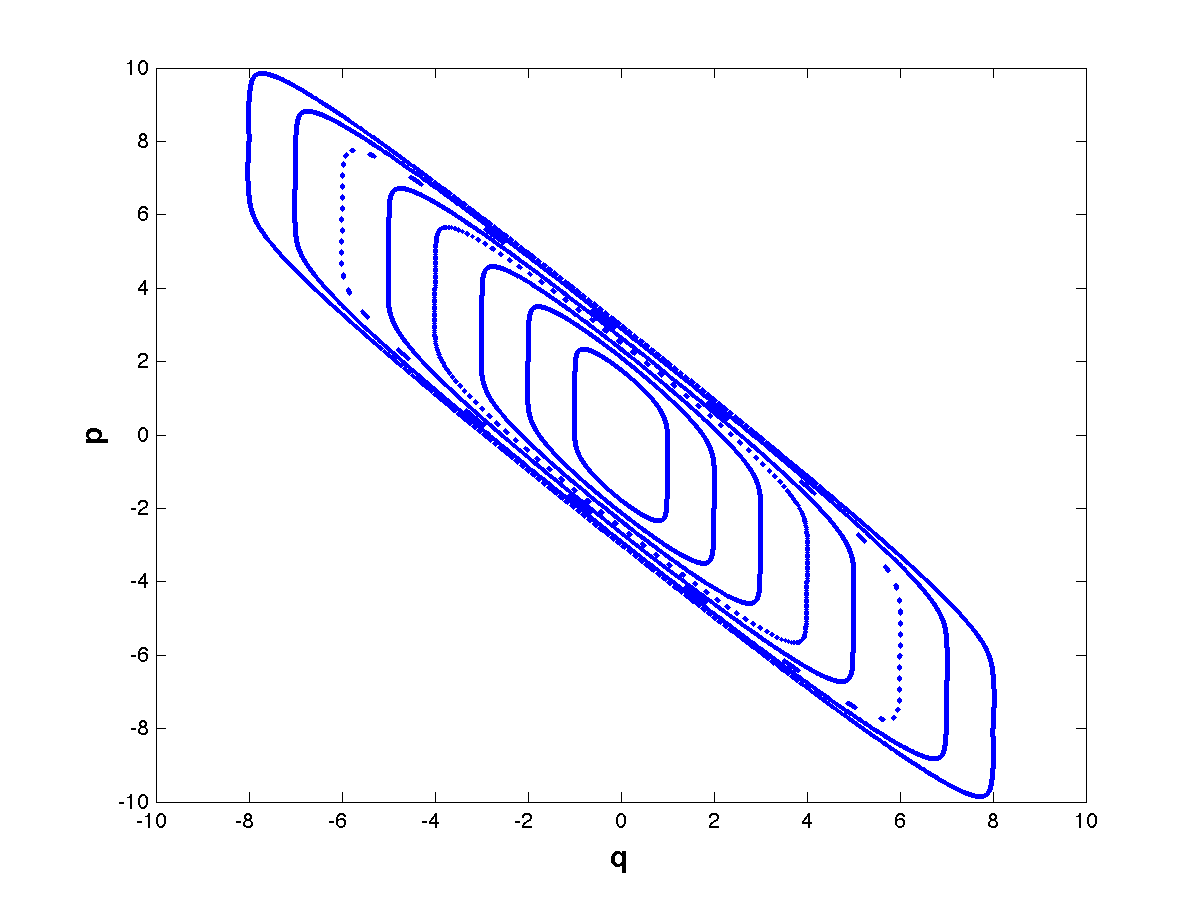}}
\caption{Numerical level curves for problem (\ref{initex1})--(\ref{parex1}), HBVM(3,2) method, $h=10^{-3}$.}
\label{k3s2}

\bigskip
\bigskip
\centerline{\includegraphics[width=12cm,height=8cm]{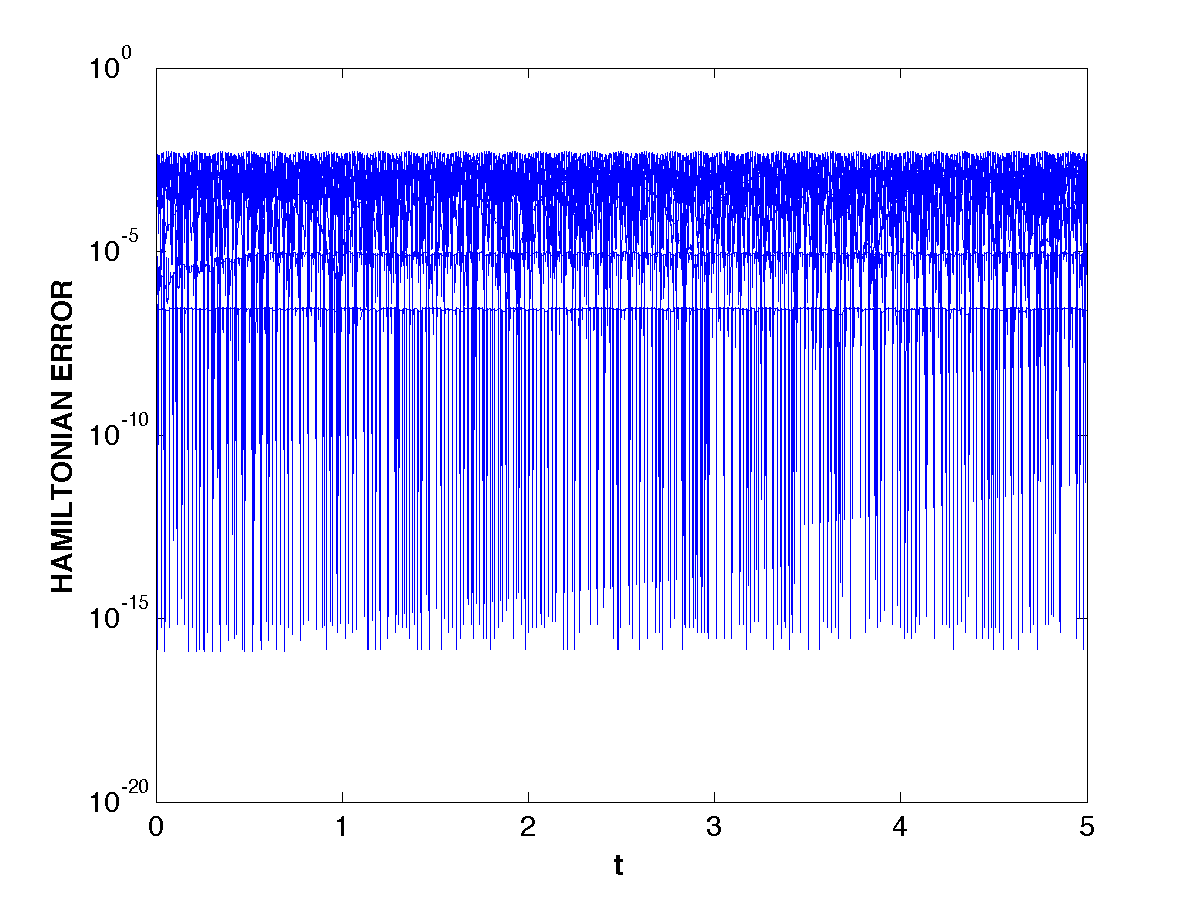}}
\caption{Hamiltonian error for problem (\ref{initex1})--(\ref{parex1}), HBVM(3,2) method, $h=10^{-3}$.}
\label{k3s2H}
\end{figure}

However, if we use HBVM$(3,2)$ with the same step-size, the error in the Hamiltonian is of sixth-order: this is enough to have a smaller error in the numerical Hamiltonian, as is shown in Figure~\ref{k3s2H}, resulting in a correct phase portrait, as is shown in Figure~\ref{k3s2}.

\begin{figure}
\centerline{\includegraphics[width=12cm,height=8cm]{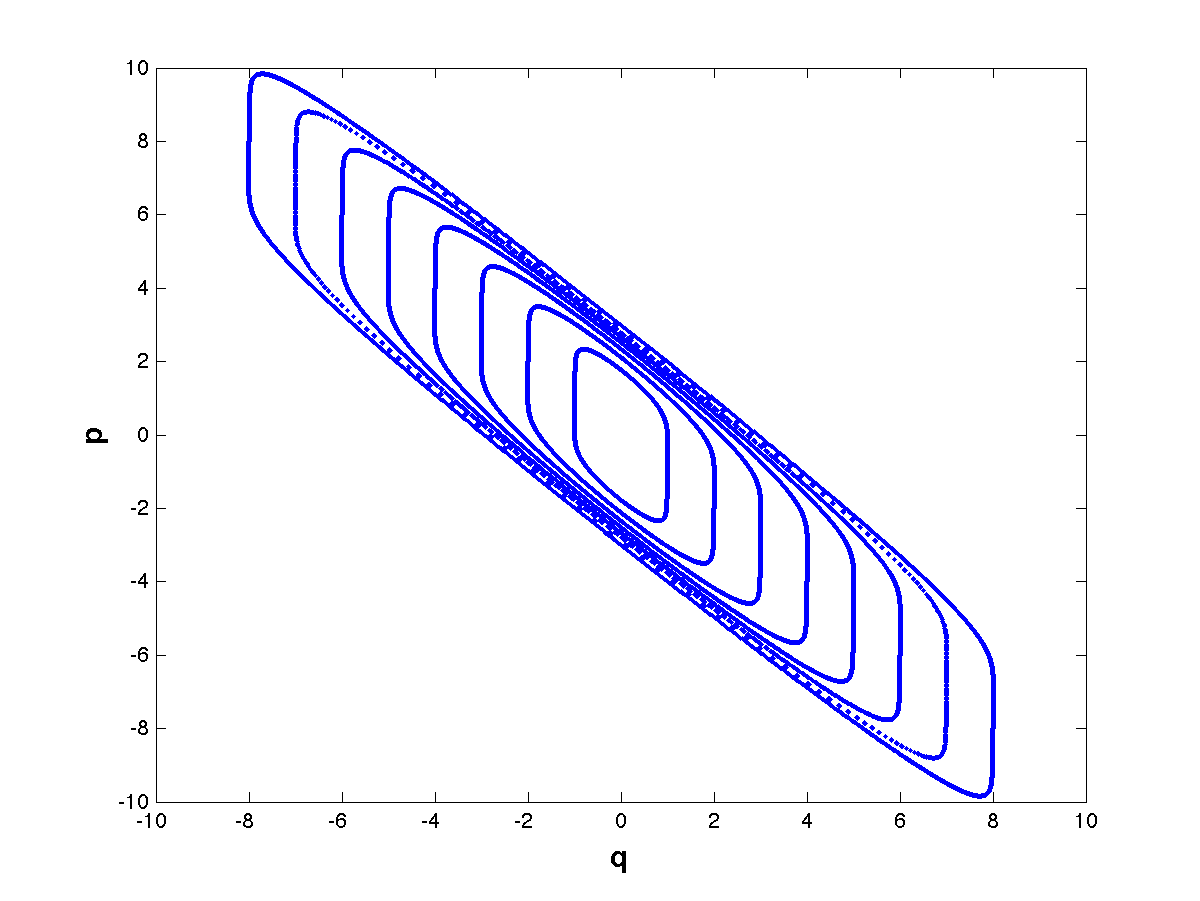}}
\caption{Numerical level curves for problem (\ref{initex1})--(\ref{parex1}), HBVM(8,2) method, $h=10^{-3}$.}
\label{k8s2}

\bigskip
\bigskip
\centerline{\includegraphics[width=12cm,height=8cm]{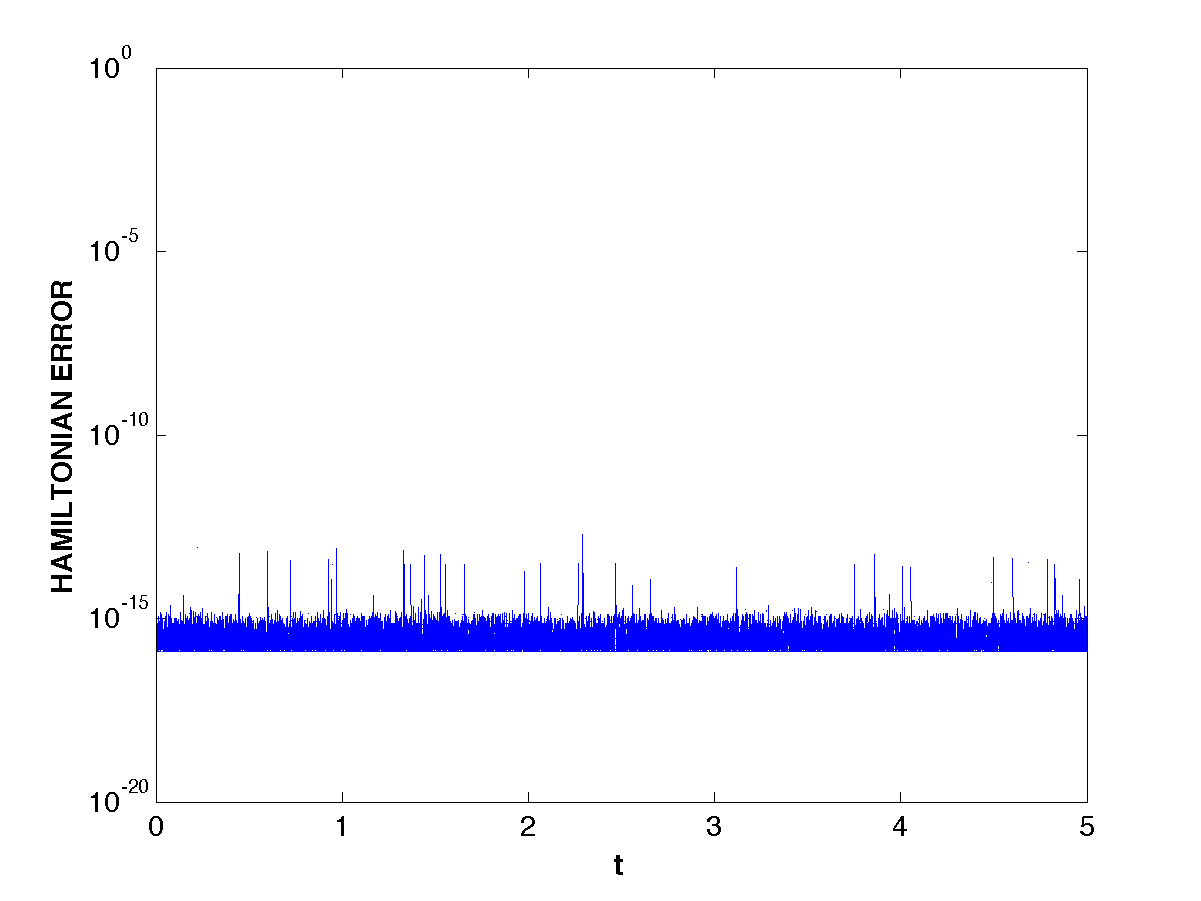}}
\caption{Hamiltonian error for problem (\ref{initex1})--(\ref{parex1}), HBVM(8,2) method, $h=10^{-3}$.}
\label{k8s2H}
\end{figure}

At last, by using HBVM$(8,2)$ with the same step size, the Hamiltonian error is of the order of roundoff errors, as is shown in Figure~\ref{k8s2H}, thus allowing a perfect reconstruction of the phase portrait, depicted in Figure~\ref{k8s2}. Indeed, since the Hamiltonian (\ref{Hex1}) has degree eight, the quadrature is exact, in this case, according to (\ref{nule2ks}).

\begin{figure}
\centerline{\includegraphics[width=12cm,height=8cm]{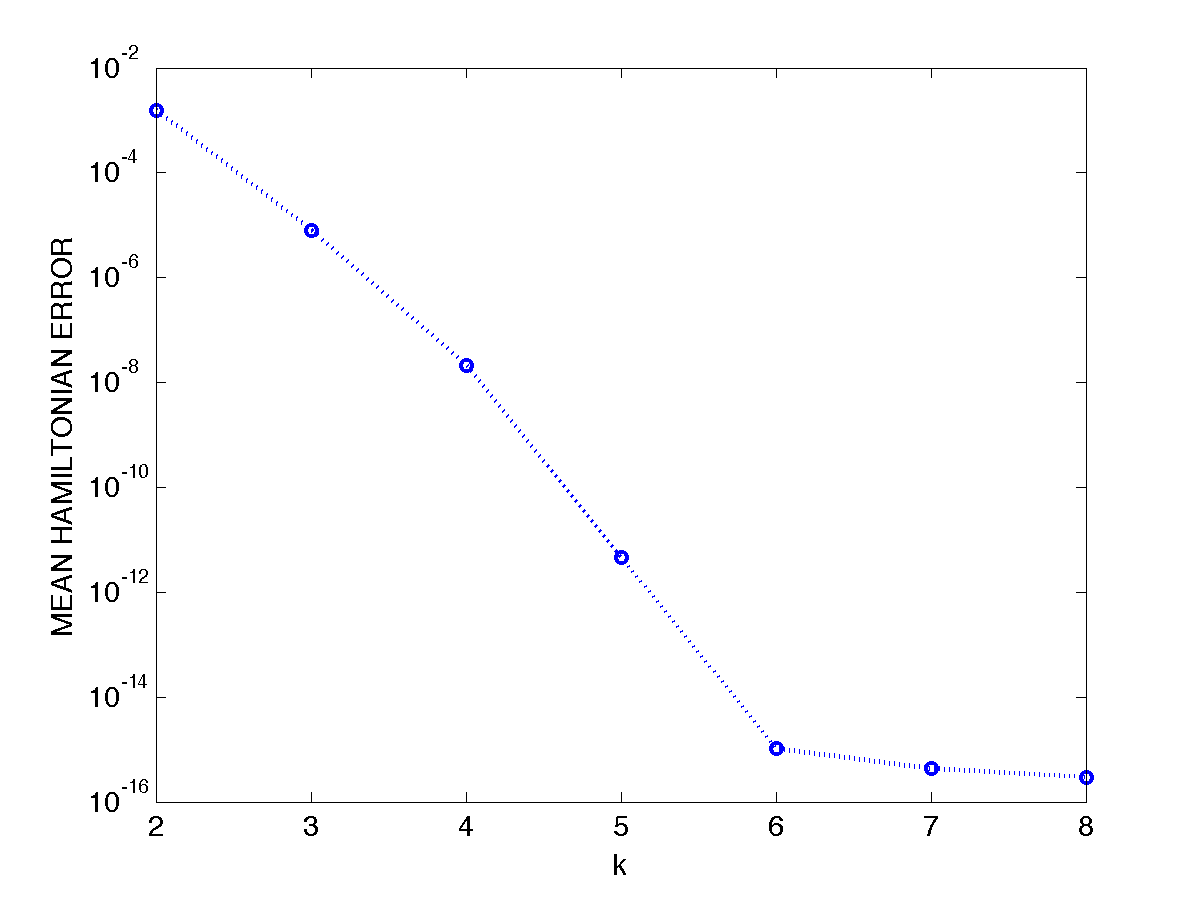}}
\caption{Mean Hamiltonian error for problem (\ref{initex1})--(\ref{parex1}), HBVM($k$,2) method, $k=2,\dots,8$, by using a stepsize $h=10^{-3}$.}
\label{erroH}
\end{figure}

For sake of completeness, in Fugure~\ref{erroH} we also plot the mean error in the numerical Hamiltonian, for HBVM$(k,2)$ methods, used with the stepsize $h=10^{-3}$, for $k=2,\dots,8$. As one can see, for the largest values of $k$ the error is essentially due to roundoff.


\section{Symmetry}\label{symmet}
We here prove that, provided that the abscissae $\{c_i\}$ are symmetrically distributed in the interval $[0,1]$, as is the case of the Gauss-Legendre nodes (see (\ref{symci})), a HBVM$(k,s)$ method is symmetric. In more detail, if applied to the initial value problem
$$y' = f(y), \qquad y(0)=y_0,$$
it provides the approximation $y_1\approx y(h)$, then it will provide the same discrete solution, as well as the same internal stages, though in reversed order, if applied to the
initial value problem
\begin{equation}\label{ivpr}
z' = -f(z), \qquad z(0) = y_1.
\end{equation}
For proving this, let us define the following matrices:
$$%\begin{equation}\label{Jr}
J_r = \pmatrix{cccc} &&&1\\ &&\cdot&\\ &\cdot&&\\ 1\endpmatrix\in\RR^{r\times r}, \qquad r=k,k+1,k+2,
$$%\end{equation}
$$%\begin{equation}\label{LandD}
L = \pmatrix{cccc} 1\\ -1&1\\ &\ddots &\ddots\\ &&-1&1\endpmatrix\in\RR^{k+1\times k+1},\qquad
D = \pmatrix{cccc} 1\\ &-1\\ &&\ddots\\ &&&(-1)^{s-1}\endpmatrix\in\RR^{s\times s},
$$%\end{equation}
and, by recalling the vector $\II_s^1$ defined at (\ref{Is1}),
$$%\begin{equation}\label{hIIs}
\hat\II_s = \pmatrix{c} \II_s\\ \II_s^1\endpmatrix\in\RR^{k+1\times s}.
$$%\end{equation}
Moreover, by setting 
\begin{equation}\label{ci01}
0\equiv c_0<c_1<\cdots<c_k<c_{k+1}\equiv1,
\end{equation}
we need to define matrix $$L\,\hat\II_s\equiv \Delta\II_s = \left( \int_{c_{i-1}}^{c_i} P_{j-1}(x)\dd x\right)_{\small\begin{array}{l}i=1,\dots,k+1\\ j=1,\dots,s\end{array}}.$$ 
The following properties then hold true, provided that the abscissae are symmetrically distributed in the interval [0,1], i.e., by taking into account (\ref{ci01}), $c_i=1-c_{k-i}$, $i=0,\dots,k+1$:
\begin{itemize}
\item $J_r^T=J_r^{-1}=J_r$;

\item $J_k\Omega J_k = \Omega\quad\Rightarrow\quad \Omega J_k = J_k\Omega$;

\item $J_{k+1} \Delta\II_s = \Delta\II_s D$;

\item $J_k\P_s=\P_sD$;
\end{itemize}
where the last two properties follow from (\ref{symIj}) and (\ref{symPj}), respectively.
The discrete solution generated by a HBVM$(k,s)$ method can then be cast in vector form as
$$\pmatrix{cc} -\hat\bfe &I_{k+1}\endpmatrix\otimes I\, \hat{Y} = h \hat\II_s\P_s^T\Omega\otimes I\, f(\hat{Y}),$$
where $\hat\bfe\in\RR^{k+1}$ is the unit vector, and
$$\hat{Y} = \pmatrix{c} y_0\\[2mm] Y\\ y_1\endpmatrix, \qquad Y = \pmatrix{c} Y_1\\ \vdots\\ Y_k\endpmatrix.$$
Left multiplication by $L\otimes I$ then gives
\begin{equation}\label{hAhB}
\hat{A}\otimes I\,\hat{Y} = h\hat{B}\otimes I\,f(\hat{Y}),
\end{equation}
with 
$$\hat{A} = \pmatrix{cccc} -1&1\\ &\ddots &\ddots\\ &&-1&1\endpmatrix,\quad
\hat{B} = \pmatrix{ccc} \bfo &\Delta\II_s\P_s^T\Omega &\bfo\endpmatrix\quad\in\RR^{k+1\times k+2}.$$ Since one easily realizes that
$$J_{k+1} \hat{A} J_{k+2} = -\hat{A},$$
the method would be symmetric provided that
$$J_{k+1} \hat{B} J_{k+2} = \hat{B}.$$ In fact, by observing that
$$\hat{Z} = J_{k+2}\otimes I\,\hat{Y} = \pmatrix{c} y_1\\[2mm] J_k\otimes I\,Y\\ y_0\endpmatrix$$
is the reversed-time discrete solution, left multiplication of (\ref{hAhB}) by $J_{k+1}\otimes I$ then gives:
\begin{eqnarray*}
\bfo&=& J_{k+1}\hat{A}\otimes I\,\hat{Y}-hJ_{k+1}\hat{B}\otimes I\,f(\hat{Y})\\
      &=& J_{k+1}\hat{A}J_{k+2}^2\otimes I\,\hat{Y}-hJ_{k+1}\hat{B}J_{k+2}^2\otimes I\,f(\hat{Y})\\
      &=& -\hat{A}\otimes I\hat{Z}-h\hat{B}\otimes I\,f(\hat{Z}).
\end{eqnarray*}
That is, the reversed-time vector satisfies the equation
$$\hat{A}\otimes I\,\hat{Z} = -h\hat{B}\otimes I\,f(\hat{Z}),$$
which consists in applying the HBVM$(k,s)$ method to problem (\ref{ivpr}), thus providing the approximation $z_1=y_0$, and using stages $Z=J_k\otimes I\,Y$. As matter of fact, one has:
\begin{eqnarray*}
J_{k+1}\hat{B}J_{k+2} &=& \pmatrix{ccc} \bfo &J_{k+1}\Delta\II_s \P_s^T\Omega J_k&\bfo\endpmatrix\\
&=&   \pmatrix{ccc} \bfo &\Delta\II_s D\P_s^TJ_k\Omega&\bfo\endpmatrix\\
&=&   \pmatrix{ccc} \bfo &\Delta\II_s D (J_k\P_s)^T\Omega&\bfo\endpmatrix\\
&=&   \pmatrix{ccc} \bfo &\Delta\II_s D^2\P_s^T\Omega&\bfo\endpmatrix ~=~ \hat{B},
\end{eqnarray*}
and the symmetry of the method follows.

\section{Linear stability analysis}
We now consider the linear stability analysis of HBVM$(k,s)$: indeed, such methods can be defined independently from the problem of energy conservation, by considering a general function $f$ in (\ref{idp3}). As matter of fact, we have seen, in Section~\ref{HBVMss}, that HBVM$(k,s)$ methods, with $k>s$ can be regarded as a low-rank generalization of the basic $s$-stage Gauss-Legendre method. 

Then, let us apply one such a method to the celebrated test equation $$y'=\lambda y, \qquad y(0)=y_0\ne0, \qquad \Re(\lambda)<0.$$
Setting $$\lambda= \aa+i\bb,\qquad y=x_1+i x_2,$$ the test equation becomes:
\begin{equation}\label{testRI}
\bfx' \equiv \pmatrix{c} x_1\\ x_2\endpmatrix' = \pmatrix{cc}\aa &-\bb\\ \bb &\aa\endpmatrix\pmatrix{c} x_1\\ x_2\endpmatrix\equiv A\bfx, \qquad \bfx(0)=\bfx_0\ne\bfo.
\end{equation}
Defining the scalar function 
\begin{equation}\label{lyapV}
V(\bfx) = \frac{1}2 \bfx^T\bfx \equiv \frac{1}2\|\bfx\|_2^2,
\end{equation}
the application of a HBVM$(k,s)$ method for solving (\ref{testRI}) defines the  polynomial $\sigma$ such that $\sigma(0)=\bfx_0$ and, moreover
\begin{eqnarray}\nonumber
\sigma'(ch) &=& \sum_{j=0}^{s-1} P_j(c) \sum_{i=1}^k b_i P_j(c_i) A\sigma(c_ih) 
                  ~\equiv~ \sum_{j=0}^{s-1} P_j(c) \int_0^1 P_j(\tau) A\sigma(\tau h)\dd\tau\\
                  &\equiv&A\sum_{j=0}^{s-1} P_j(c) \int_0^1 P_j(\tau) \nabla V(\sigma(\tau h))\dd\tau. 
\label{sig1}
\end{eqnarray}
By considering that $\sigma(0)=\bfx_0$, and the new approximation is defined by $$\bfx_1\equiv \sigma(h),$$ one obtains:
\begin{eqnarray}\nonumber
\Delta V(\bfx_0) &=& V(\bfx_1)-V(\bfx_0) ~=~ V(\sigma(h))-V(\sigma(0))\\
\nonumber
&=& \int_0^h \frac{\dd}{\dd t} V(\sigma(t))\dd t ~=~ \int_0^h \nabla V(\sigma(t))^T\sigma'(t)\dd t\\
\nonumber
&=& h\int_0^1 \nabla V(\sigma(\tau h))^TA\sum_{j=1}^{s-1} P_j(\tau) \left[\int_0^1 P_j(c)\nabla V(\sigma(c h))\dd c\right]\dd\tau\\
\nonumber
&=& \aa h \sum_{j=0}^{s-1} \left\| \int_0^1 P_j(\tau)\nabla V(\sigma(\tau h))\dd\tau\right\|_2^2\\ 
&=& \aa h \sum_{j=0}^{s-1} \left\| \int_0^1 P_j(\tau)\sigma(\tau h)\dd\tau\right\|_2^2 
~\equiv~\aa h \Gamma^2.\label{uno}
\end{eqnarray}
Moreover, the following result holds true.

\begin{lem}\label{lem0} \quad 
$\Gamma^2=0\quad\Rightarrow\quad \bfx_0=\bfo$.\end{lem}
\proof Indeed, one has:
$$\Gamma^2=0\qquad \Rightarrow \qquad \sigma'(ch)\equiv \bfo \mbox{\quad and\quad} \int_0^1\overbrace{P_0(ch)}^{\equiv\,1}\sigma(ch)\dd c=\int_0^1\sigma(ch)\dd c=\bfo.$$ From the first equality one obtains $\sigma(ch)\equiv \bfx_0$ and, therefore, from the second equality one derives $\bfx_0=\bfo$.\,\QED\bigskip

\no From (\ref{uno}) and Lemma~\ref{lem0}, the following result then easily follows.

\begin{theo}\label{pAstab} For all $k\ge s$, and for any choice of the nodes, HBVM$(k,s)$ is {\em perfectly $A$-stable}, i.e., its stability region coincides with the negative-real complex plane, $\CC^-$\end{theo}

\proof From (\ref{uno}) and Lemma~\ref{lem0}, one has, by considerdering (lyapV) and that $\aa=\Re(\lambda)$:
$$\|\bfx_1\|_2^2 = \|\bfx_0\|_2^2 +\aa h\Gamma ^2< \|\bfx_0\|_2^2 \qquad \Leftrightarrow\qquad
\Re(\lambda)<0.$$ Consequently, a HBVM$(k,s)$ method turns out to be perfectly $A$-stable, since its absolute stability region coincides with $\CC^-$, for all $k\ge s\ge1$.\,\QED

\chapter{Implementation of the methods}

In this chapter, we discuss the efficient implementation of HBVM$(k,s)$ methods. In particular, it is clearly shown that their computational cost depends essentially on $s$,  in the sense that, for all $k\ge s$, the discrete problem turns out to have always block-dimension $s$. A nonlinear iteration procedure, based on the {\em blended implementation} of the methods is also sketched. The material in this chapter is based on \cite{IT09,BIT09,BIT10_0,BIT11,Br00,BM02,BM04,BMM06,BM09}.

\section{Fundamental and silent stages}
From (\ref{kstage1})-(\ref{kstage2}), we see that a HBVM$(k,s)$ method, with $k>s$, is defined by a Butcher matrix of rank $s$. Consequently, $k-s$ of the stages of the method can be expressed as a linear combination of the remaining $s$ stages: we shall, therefore, name {\em fundamental stages} the latter ones, and {\em silent stages} the former ones. For this purpose, let us partition the stage vector $Y$ as
$$Y = \pmatrix{c} Y^{(1)}\\ Y^{(2)} \endpmatrix$$
where, by supposing for sake of brevity that the fundamental satges are the first $s$-ones,\footnote{Indeed, this can be always achieved, by using a suitable premutation of the abscissae.} 
$$Y^{(1)} = \pmatrix{c} Y_1\\ \vdots\\ Y_s\endpmatrix, \qquad Y^{(2)} = \pmatrix{c} Y_{s+1}\\ \vdots\\ Y_k\endpmatrix.$$
Similarly, we partition matrices $\II_s$ and $\P_s$, respectively, as
$$\II_s = \pmatrix{c} \II_s^{(1)}\\ \II_s^{(2)}\endpmatrix,
\qquad \P_s = \pmatrix{c} \P_s^{(1)}\\ \P_s^{(2)}\endpmatrix,$$ 
containing the corresponding rows as those of $Y^{(1)}$ and $Y^{(2)}$, respectively. 
Moreover, we also consider the partition
$$\Omega = \pmatrix{cc} \Omega_1\\ &\Omega_2\endpmatrix, \qquad \Omega_1\in\RR^{s\times s},
\qquad \Omega_2\in\RR^{k-s\times k-s}.$$ 
Consequently, by setting $\bfe^{(1)}$ end $\bfe^{(2)}$ the unit vectors of length $s$ and $k-s$, respectively, one obtains:
\begin{eqnarray}\label{Y1}
Y^{(1)} &=& \bfe^{(1)} \otimes y_0 + h \II_s^{(1)}\P_s^T\Omega \otimes I\, \pmatrix{c} f(Y^{(1)})\\ f(Y^{(2)})\endpmatrix,\\
Y^{(2)} &=& \bfe^{(2)} \otimes y_0 + h \II_s^{(2)}\P_s^T\Omega \otimes I\, \pmatrix{c} f(Y^{(1)})\\ f(Y^{(2)})\endpmatrix.\label{Y2}
\end{eqnarray}
From (\ref{Y1}), one then obtains that
$$\P_s^T\Omega \otimes I\, \pmatrix{c} f(Y^{(1)})\\ f(Y^{(2)})\endpmatrix = \left( h \II_s^{(1)}\right)^{-1}\left[ Y^{(1)}-\bfe^{(1)}\otimes y_0\right],$$ which substituted into (\ref{Y2}) gives:
 \begin{eqnarray*}
 Y^{(2)} &=& \bfe^{(2)} \otimes y_0 + \II_s^{(2)} \left(\II_s^{(1)}\right)^{-1}\left[ Y^{(1)}-\bfe^{(1)}\otimes y_0\right]\\
 [2mm]
 &=&\underbrace{\left[ \bfe^{(2)}- \II_s^{(2)} \left(\II_s^{(1)}\right)^{-1} \bfe^{(1)}\right]}_{=\,\bfa}\otimes y_0 +\II_s^{(2)} \left(\II_s^{(1)}\right)^{-1}Y^{(1)}\\
 &\equiv& \bfa \otimes y_0 +\II_s^{(2)} \left(\II_s^{(1)}\right)^{-1}Y^{(1)}.
 \end{eqnarray*}
 Consequently, we can rewrite (\ref{Y1})-(\ref{Y2}) as:
 \begin{eqnarray}\nonumber
 Y^{(1)} &=& \bfe^{(1)} \otimes y_0 + h \II_s^{(1)}\P_s^T\Omega \otimes I\,
 \pmatrix{c}
  f(Y^{(1)}) \\ 
  f\left(\bfa \otimes y_0 +\II_s^{(2)} (\II_s^{(1)})^{-1}Y^{(1)}\right)\endpmatrix\\ \nonumber
  &\equiv& \bfe^{(1)} \otimes y_0 + h \II_s^{(1)}\left[(\P_s^{(1)})^T\Omega_1 \otimes I\,f(Y^{(1)}) ~+\right.\\
  &&\left.(\P_s^{(2)})^T\Omega_2 \otimes I\, f\left(\bfa \otimes y_0 +\II_s^{(2)} (\II_s^{(1)})^{-1}Y^{(1)}\right)\right],
\label{Y1only}
\end{eqnarray}
involving only the fundamental stages, thus confirming that the actual discrete problem, to be solved at each time step, amounts to a set of $s$ (generally) nonlinear equations, each having the same size as that of the continuous problem. For solving such a problem, one could use, e.g., a {\em fixed-point iteration},
 \begin{equation}\label{fixit}
 Y_{\ell+1}^{(1)} = \bfe^{(1)} \otimes y_0 + h \II_s^{(1)}\P_s^T\Omega \otimes I\,
 \pmatrix{c}
  f(Y_\ell^{(1)})\\
 f\left(\bfa \otimes y_0 +\II_s^{(2)} (\II_s^{(1)})^{-1}Y_\ell^{(1)}\right)
 \endpmatrix,
\qquad \ell=0,1,\dots,
\end{equation}
or, if the case, a {\em simplified-Newton iteration}. In more details, setting
\begin{eqnarray*}
F(Y^{(1)}) &=& Y^{(1)} - \bfe^{(1)} \otimes y_0 - h \II_s^{(1)}\left[(\P_s^{(1)})^T\Omega_1 \otimes I\, f(Y^{(1)})~+\right.\\
&&\left. (\P_s^{(2)})^T\Omega_2 \otimes I\, f\left(\bfa \otimes y_0 +\II_s^{(2)} (\II_s^{(1)})^{-1}Y^{(1)}\right)\right],
\end{eqnarray*}
one then solves,
\begin{equation}\label{NewtC}
\left[I-hC\otimes J_0\right] \Delta_\ell = -F(Y_\ell^{(1)}), \qquad 
Y_{\ell+1}^{(1)} =  Y_\ell^{(1)} + \Delta_\ell, \qquad \ell=0,1,\dots,\end{equation} 
where  $J_0=J_f(y_0)$ and matrix $C$ is defined as follows:
\begin{equation}\label{C}
C = \II_s^{(1)}\left[(\P_s^{(1)})^T\Omega_1 + (\P_s^{(2)})^T\Omega_1\II_s^{(2)} (\II_s^{(1)})^{-1}\right]
\end{equation}
The following result holds true.

\begin{figure}
\centerline{\includegraphics[width=12cm,height=8cm]{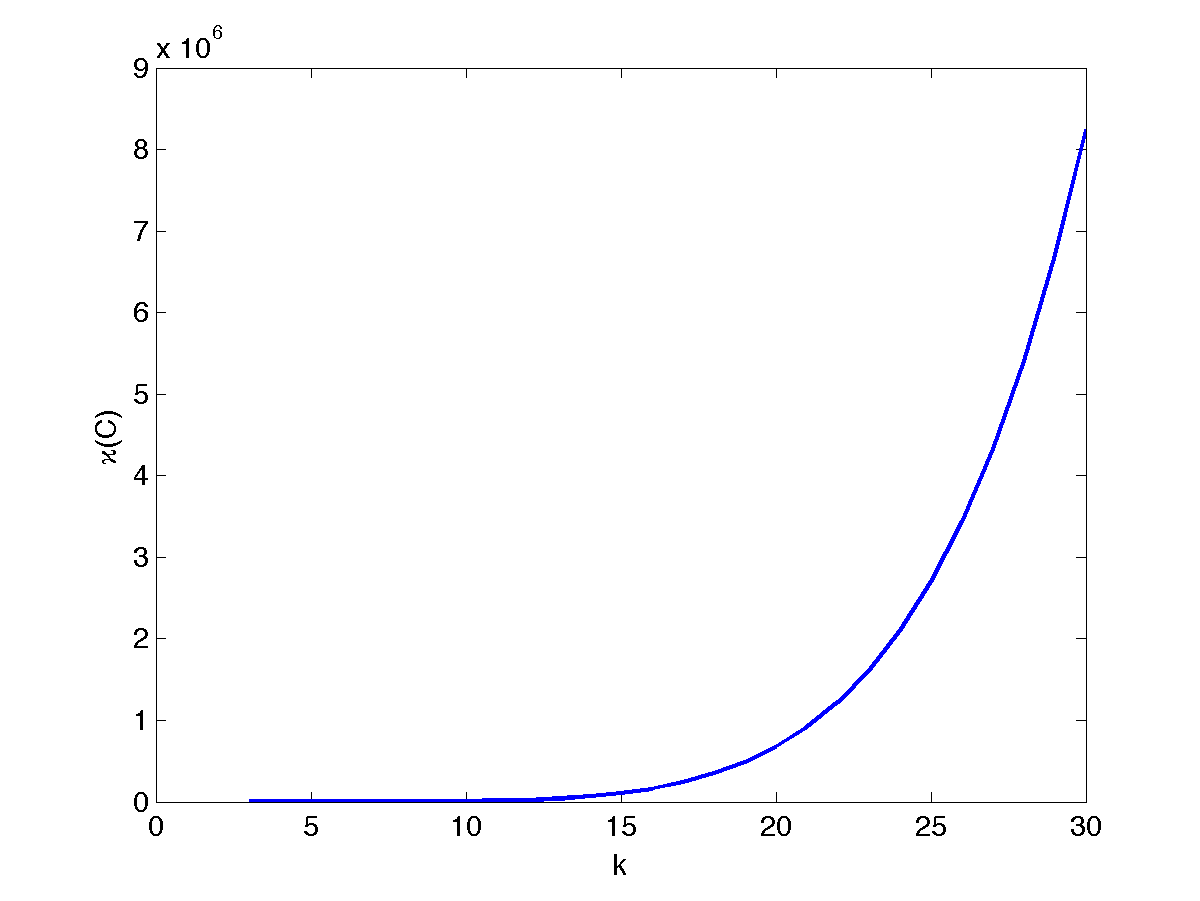}}
\caption{Condition number of matrix (\ref{C}), fundamental abscissae fixed at the first $s$ ones.}
\label{condC1}

\bigskip
\bigskip
\centerline{\includegraphics[width=12cm,height=8cm]{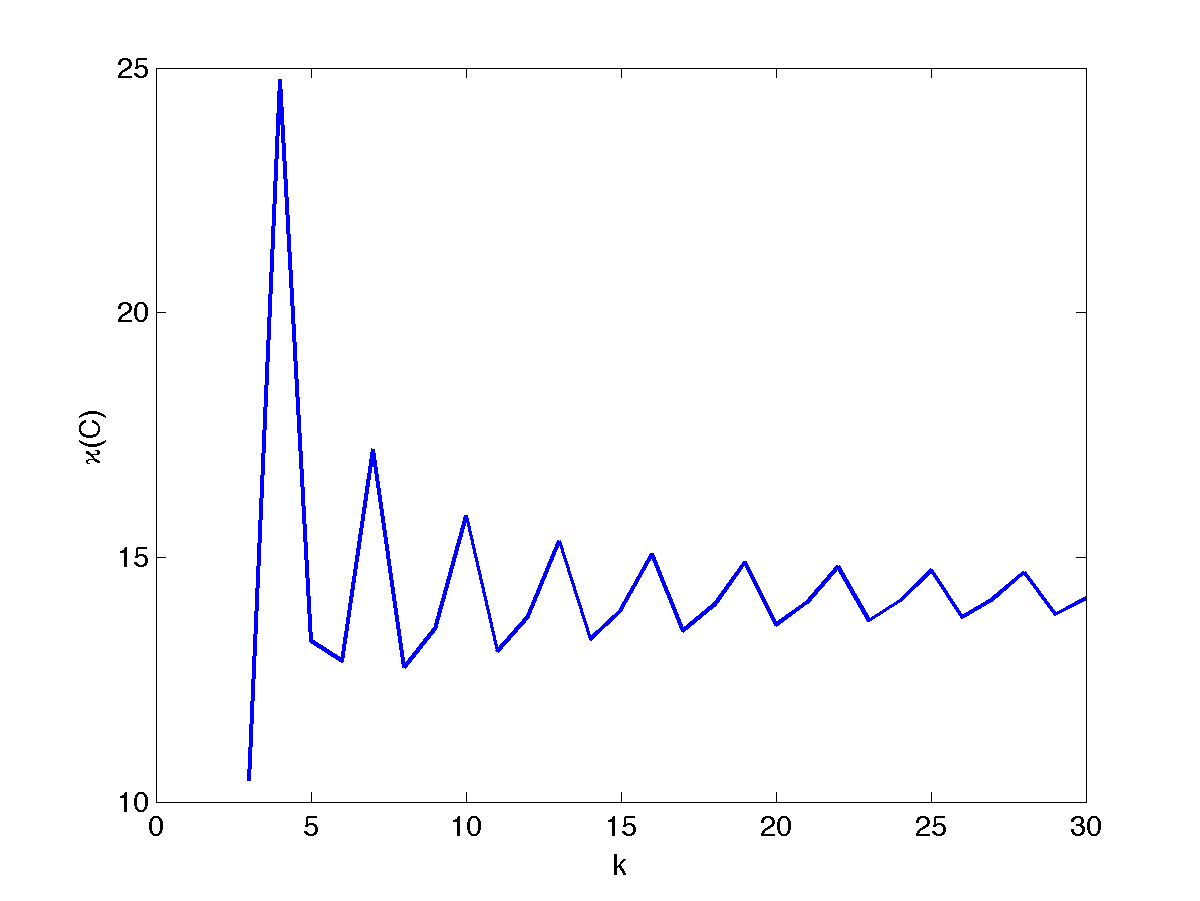}}
\caption{Condition number of matrix (\ref{C}), fundamental abscissae almost evenly spaced.}
\label{condC2}
\end{figure}

\begin{theo}\label{isoC}
The eigenvalues of matrix $C$, as defined in (\ref{C}), coincide with those of matrix $X_s$ defined in (\ref{Xs}), that is the eigenvalues of the Butcher matrix of the $s$-stage Gauss method.
\end{theo}
\proof One has:
\begin{eqnarray*}
C &=& \II_s^{(1)}\left[(\P_s^{(1)})^T\Omega_1 + (\P_s^{(2)})^T\Omega_2\II_s^{(2)} (\II_s^{(1)})^{-1}\right]\\
&=& \II_s^{(1)}\left[(\P_s^{(1)})^T\Omega_1\II_s^{(1)} + (\P_s^{(2)})^T\Omega_2\II_s^{(2)} \right](\II_s^{(1)})^{-1}\\
&=& \II_s^{(1)}\left[\P_s^T\Omega \II_s\right](\II_s^{(1)})^{-1}\\
&\sim& \P_s^T\Omega \II_s\\
&=& \P_s^T\Omega\P_{s+1}\hat{X}_s\\
&=& [I_s~\bfo]\hat{X_s} ~=~ X_s.\,\QED
\end{eqnarray*}
Consequently, matrix $C$ has {\em always} the same spectrum, independently of the choice of the {\em fundamental} and {\em silent abscissae}.\footnote{I.e., the abscissae corresponding to the fundamental and silent stages, respectively.} This, in turn, coincides with the nonzero eigenvalues of the corresponding Butcher array (see Theorem~\ref{ranks}). Nevertheless, its condition number is greatly affected from this choice. Clearly, a badly conditioned matrix $C$ would affect the convergence of both the iterations (\ref{fixit}) and (\ref{NewtC}). 
As an example, in Figures~\ref{condC1} and \ref{condC2} we plot the condition number of matrix $C$ corresponding to the following choices of the fundamental abscissae, in the case $k\ge s=3$:
\begin{itemize}
\item the first $s$ abscissae of the $k$ ones (Figure~\ref{condC1});

\item $s$ almost evenly spaced abscissae among the $k$ ones  (Figure~\ref{condC2}).
\end{itemize}
As one may see, in the first case $\kappa(C)$ grows exponentially with $k$, whereas it is uniformly bounded in the second case. Because of this reason, we shall consider a more favorable formulation of the discrete problem itself, which will be independent of the choice of the fundamental abscissae.

\section{Alternative formulation of the discrete problem}
In order to overcome the previous drawback, the basic idea is to reformulate the discrete problem by considering as unknowns the coefficients, say 
$$\hat\gamma_j = \sum_{\ell=1}^k b_\ell P_{j}(c_\ell) f(u(c_\ell h)), \qquad j=0,\dots,s-1,$$
of the polynomial approximation defining a HBVM$(k,s)$ method (see (\ref{uc})). In more details, recalling that
$$Y_i \equiv u(c_i h) = y_0 +h\sum_{j=0}^{s-1} \hat\gamma_j \int_0^{c_i} P_j(x)\dd x, \qquad i=1,\dots,k,$$
one may cast the disctere problem as follows:
\begin{equation}\label{hgamma}
\hat{\bfgamma}\equiv\pmatrix{c} \hat\gamma_0\\ \vdots\\ \hat\gamma_{s-1}\endpmatrix =
\P_s^T\Omega\otimes I \, f\left( \bfe\otimes y_0+h\II_s\otimes I\,\hat{\bfgamma}\right),
\end{equation} with the new approximation given by $$y_1 = y_0+h\hat\gamma_0.$$
We observe that (\ref{hgamma}) has always (block) dimension $s$, whatever is the value of $k$ considered. For solving such a problem, one can still use a {\em fixed-point iteration},
$$\hat{\bfgamma}^{\ell+1} =
\P_s^T\Omega\otimes I \, f\left( \bfe\otimes y_0+h\II_s\otimes I\,\hat{\bfgamma}^{\ell}\right),\qquad \ell=0,1,\dots,$$ whose implementation is straightforward. One can also consider a {\em simplified-Newton iteration}. Setting
\begin{equation}\label{Fgamma}
F(\hat{\bfgamma}) = \hat{\bfgamma}-
\P_s^T\Omega\otimes I \, f\left( \bfe\otimes y_0+h\II_s\otimes I\,\hat{\bfgamma}\right),
\end{equation}
and, as before, $J_0=J_f(y_0)$, it takes the form
\begin{equation}\label{NewtC1}
\left[I-hC\otimes J_0\right] \Delta^\ell = -F(\hat{\bfgamma}^\ell), \qquad 
\hat{\bfgamma}^{\ell+1} = \hat{\bfgamma}^\ell + \Delta^\ell, \qquad \ell=0,1,\dots,\end{equation} 
where matrix $C$ is now defined as follows:
\begin{equation}\label{Cnew}
C = \P_s^T\Omega \II_s = \P_s^T\Omega \P_{s+1}\hat{X}_s = (I_s~\bfo)\hat{X}_s = X_s.
\end{equation}
Consequently, the iteration (\ref{NewtC1}) becomes:
\begin{equation}\label{NewtC2}
\left[I-hX_s\otimes J_0\right] \Delta^\ell = -F(\hat{\bfgamma}^\ell), \qquad 
\hat{\bfgamma}^{\ell+1} = \hat{\bfgamma}^\ell + \Delta^\ell, \qquad \ell=0,1,\dots.\end{equation}
\begin{rem} It is worth noticing that (\ref{NewtC2}) holds independently of the choice of the $k$ abscissae $\{c_i\}$, the only requirement being the order $2s$ of the quadrature, so that the property $\P_s^T\Omega \P_{s+1}\hat{X}_s = (I_s~\bfo)$ holds true. 
\end{rem}
\begin{rem}
We observe that both matrices (\ref{C}) and (\ref{Cnew}) share the same eigenvalues which, in turn, are the nonzero eigenvalues of the Butcher array of the given HBVM$(k,s)$ method (see Theorem~\ref{ranks}).
\end{rem}
We observe that, remarkably enough, at each step of the simplified-Newton iteration we have to solve a linear system of dimension $sm\times sm$ of the form
\begin{equation}\label{sm}\left[ I-hX_s\otimes J_0\right]\bfx = \bfeta,\end{equation}
whose coefficient matrix is thus {\em independent} of $k$ and of the choice of the abscissae. Its cost is then approximately given by $$\frac{2}3 (sm)^3 \qquad flops,$$ due to the cost of the corresponding $LU$ factorization. In the next section, we shall consider an alternative, iterative, procedure for solving (\ref{sm}), able to reduce the cost for the factrization to $$\frac{2}3 m^3 \qquad flops.$$

\section{{\em Blended} HBVMs} 
We now introduce an iterative procedure for solving (\ref{sm}), 
which has been already succesfully implemented in the computational codes {\tt BiM} \cite{BM04} and {\tt BiMD} \cite{BMM06} for the numerical solution of stiff ODE-IVPs and linearly implicit DAEs up to order 3.

For this iterative procedure a linear analysis of convergence is provided. To this purpose, let us consider the ``usual'' test equation,
\begin{equation}\label{testeq}
y' = \lambda y, \qquad \Re(\lambda)<0.
\end{equation}
In such a case, by setting as usual $q=h\lambda$, problem (\ref{sm}) becomes the linear system, of dimension $s$,
\begin{equation}\label{b1}
(I-qX_s)\bfx = \bfeta.\end{equation}
The solution of this linear system is not affected by left-multiplication by $\zeta X_s^{-1}$, where $\zeta>0$ is a free parameter to be chosen later. Thus, we obtain the following equivalent formulation of (\ref{b1}):
\begin{equation}\label{b2}
\zeta(X_s^{-1} -qI)\bfx = \zeta X_s^{-1}\bfeta \equiv \bfeta_1.
\end{equation}
Let us define the {\em weighting function}
\begin{equation}\label{teta}
\theta(q) = I(1-\zeta q)^{-1}.
\end{equation}
It satisfies the following properties:
\begin{itemize}
\item $\theta(q)$ is well defined for all $q\in\CC^-$, since $\zeta>0$;

\item $\theta(0)=I$;

\item $\theta(q)\rightarrow O$, as $q\rightarrow\infty$.

\end{itemize}
We can derive a further equivalent formulation of problem (\ref{b1}), as the {\em blending}, with weights $\theta(q)$ and $I-\theta(q)$ of the two equivalent formulations (\ref{b1}) and (\ref{b2}), thus obtaining
\begin{equation}\label{blend}
M(q)\bfx = \bfeta(q),\end{equation}
with:
\begin{eqnarray}\nonumber
M(q)       &=& \theta(q) (I-qX_s) + \zeta(I-\theta(q))(X_s^{-1}-qI),\\[-2mm]
\label{Meta}\\
\bfeta(q) &=& \theta(q)\bfeta + \zeta(I-\theta(q))X_s^{-1}\bfeta.\nonumber
\end{eqnarray}
Equations (\ref{blend})-(\ref{Meta}) define the {\em blended formulation} of the original problem (\ref{b1}). The next step is now to devise an iterative procedure, defined by a suitable splitting, for solving (\ref{blend}). To this end we observe that, due to the properties of the weighting function $\theta(q)$ defined in (\ref{teta}), one has:
\begin{eqnarray*}
M(q) &\approx& I, \qquad q\approx 0,\\
M(q) &\approx& -\zeta qI, \qquad |q|\gg1.
\end{eqnarray*}
Consequently, $N(q)=I(1-\zeta q) \approx M(q)$, both for $q\approx0$, and $|q|\gg1$. It is then natural to define the following iterative procedure, for solving (\ref{blend}):
$$N(q)\bfx_{r+1} = (N(q)-M(q))\bfx_r + \bfeta(q), \qquad r=0,1,\dots.$$ That is, observing that $N(q)^{-1}=\theta(q)$:
\begin{equation}\label{blit}
\bfx_{r+1} = (I-\theta(q)M(q))\bfx_r +\theta(q)\bfeta(q), \qquad r=0,1,\dots.
\end{equation} Equation (\ref{blit}) defines the {\em blended iteration} associated with the blended formulation (\ref{blend}) of the problem.
By considering that the solution, say $\bfx^*$, of (\ref{blend}) satisfies also (\ref{blit}), by setting
$$\bfe_r=\bfx_r-\bfx^*$$ the error at the $r$th iteration, one then obtains the {\em error equation}
$$\bfe_{r+1} = (I-\theta(q)M(q))\bfe_r, \qquad r=0,1,\dots.$$
Consequently, the iteration (\ref{blit}) will converge to the solution $\bfx^*$ of the problem iff the spectral radius of the {\em iteration matrix},
\begin{equation}\label{roqZq}
\rho(q) =\max_{\xi\in\sigma(Z(q))}|\xi|, \qquad\mbox{with}\qquad Z(q) =I-\theta(q)M(q),
\end{equation} is less than 1, where $\sigma(\cdot)$ denotes the spectrum of the matrix in argument.
The set $$\DD = \left\{ q\in\CC\,:\,\rho(q)<1\right\}$$ is the {\em region of convergence} of the iteration (\ref{blit}). The iteration will be said to be:
\begin{itemize}
\item $A$-convergent if $\CC^-\subseteq \DD$;
\item $L$-convergent if, in addition, $\rho(q)\rightarrow0$, as $q\rightarrow\infty$.
\end{itemize}  
\begin{rem}
$A$-convergent iterations are then appropriate when the underlying method is $A$-stable. Similarly, $L$-convergent iterations are appropriate in the case of  $L$-stable methods.
\end{rem}
We observe that:
\begin{itemize}
\item $Z(0)=O ~\quad \Rightarrow\quad \rho(0)=0$;
\item $Z(q)\rightarrow O \quad\Rightarrow\quad \rho(q)\rightarrow0$, \quad as \quad $q\rightarrow\infty$;
\item $Z(q)$ is well-defined for all $q\in\CC^-$, since $\zeta>0$.
\end{itemize}
Consequently, for the blended iteration (\ref{blit}) $A$-convergence and $L$-convergence are equivalent to each other. From the maximum modulus theorem, in turn, it follows that this is equivalent to requiring that the {\em maximum amplification factor},
$$\rho^*=\sup_{\Re(q)=0} \rho(q) = \sup_{x\in\RR} \rho(ix),$$ satisfies $$\rho^*\le1.$$
For the blended iteration, due to the fact that $\rho(q)\rightarrow0$, as $q\rightarrow\infty$, and since the matrix $X_s$ is real, so that $\rho(\bar{q})=\rho(q)$, one has actually to prove that
\begin{equation}\label{Lconv}
\rho^* = \max_{x>0} \rho(ix)\le 1.
\end{equation} We shall choose the free positive parameter $\zeta$, in order to minimize $\rho^*$, so that (\ref{Lconv}) turns out to be fulfilled for all $s\ge1$ (see (\ref{b1}) and (\ref{Xs})). The following result holds true.

\begin{theo}\label{muXs}
$\mu\in\sigma(X_s) \quad\Leftrightarrow \quad \frac{\displaystyle q(\mu-\zeta)^2}{\displaystyle \mu(1-q\zeta)^2} \in\sigma(Z(q)).$%%\footnote{In this setting, $\sigma(\cdot)$ denotes the spectrum of the matrix in argument.}
\end{theo}
\proof From (\ref{roqZq}), (\ref{Meta}), and (\ref{teta}), one obtains:
\begin{eqnarray*}
Z(q) &=& I-\theta(q)M(q)\\
        &=& I-\theta(q)\left[ \theta(q)(I-qX_s) +\zeta(I-\theta(q))(X_s^{-1}-qI)\right]\\
        &=& I-\theta(q)^2\left[ (I-qX_s) -\zeta^2q (X_s^{-1}-qI)\right]\\
        &=& \theta(q)^2\left[ (1+\zeta^2q^2-2\zeta q)I -I+qX_s +\zeta^2qX_s^{-1}-\zeta^2q^2I\right]\\
        &=& q\theta(q)^2X_s^{-1}\left[ X_s^2 -2\zeta X_s+\zeta^2I\right] \\
        &=&q\theta(q)^2X_s^{-1}(X_s-\zeta I)^2\\
        &\equiv& q(X_s-\zeta I)^2 \left[ X_s(1-\zeta q)^2I\right]^{-1},
\end{eqnarray*}
from which the thesis easily follows.\,\QED
\bigskip

\no As a consequence, one obtains the following result.
\begin{cor}\label{corrho} The maximum amplification factor (\ref{Lconv}) of the blended iteration (\ref{blit}) is given by: 
$$\rho^* = \max_{\mu\in\sigma(X_s)} \frac{|\mu-\zeta|^2}{2\zeta|\mu|}.$$
\end{cor}
\proof One has:
$$
\rho^* ~=~ \max_{x>0} \max_{\mu\in\sigma(X_s)}~ \frac{x|\mu-\zeta|^2}{|\mu|\,|1-ix\zeta|^2}
~=~ \max_{x>0} \frac{x}{1+\zeta^2x^2} ~ \max_{\mu\in\sigma(X_s)} \frac{|\mu-\zeta|^2}{|\mu|}.
$$
The thesis then follows immediately, by considering that
$$\max_{x>0} \frac{x}{1+\zeta^2x^2} ~=~\frac{1}{2\zeta},$$
which is obtained at $x=\zeta^{-1}$.\,\QED\bigskip

\begin{figure}
\centerline{\includegraphics[width=12cm,height=8cm]{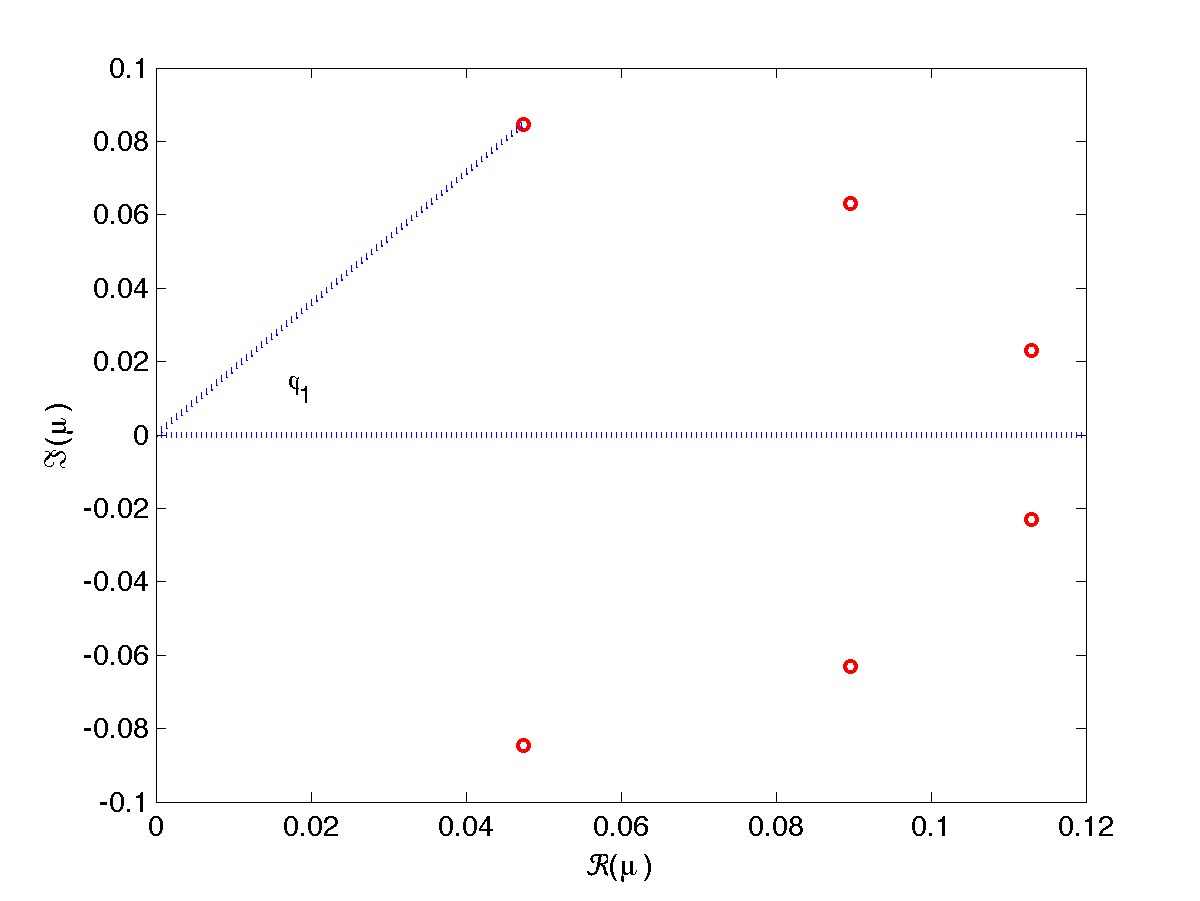}}
\caption{Eigenvalues of matrix $X_6$.}
\label{eigX6}

\bigskip
\bigskip
\centerline{\includegraphics[width=12cm,height=8cm]{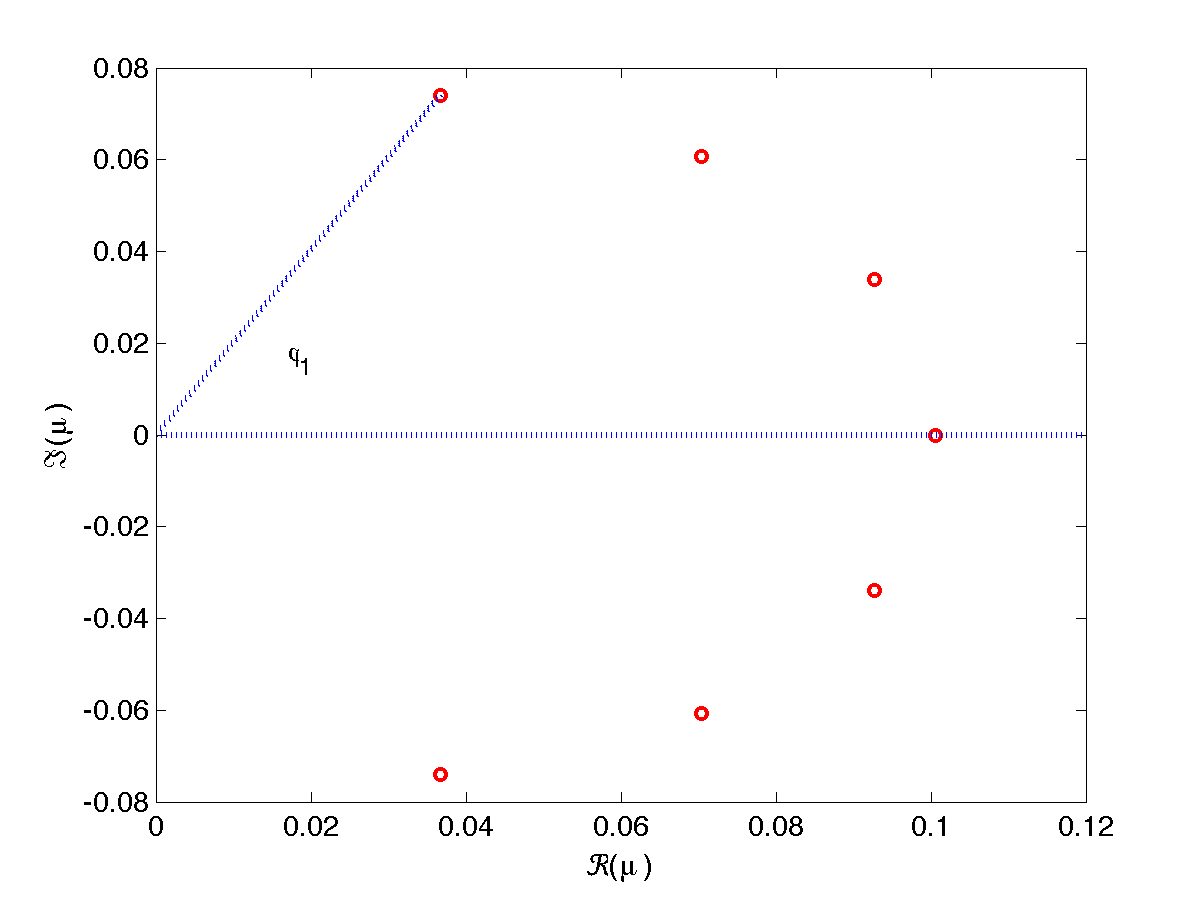}}
\caption{Eigenvalues of matrix $X_7$.}
\label{eigX7}
\end{figure}

\no We are now in the position to choose the positive parameter $\zeta$ in order for $\rho^*$ to be minimized. 
This clearly will depend on the eigenvalues of matrix $X_s$. Since this matrix is real, the complex ones occur as complex-conjugate pairs. Consequently, if we set
$$\mu_j = |\mu_j|\ee^{i\phi_j}, \qquad j=1,\dots,s,$$ we can sort them by decreasing arguments:
$$\frac{\pi}2> \phi_1 >\phi_2>\dots >\phi_s>-\frac{\pi}2,$$ due to the fact that $$\Re(\mu_j)>0, \qquad j=1,\dots,s.$$ Moreover, we can neglect the complex conjugate ones, thus obtaining:
$$\frac{\pi}2>\phi_1>\dots \phi_\ell\ge0, \qquad \ell =\lceil\frac{s}2\rceil.$$
In addition to this, it turns out that the eigenvalues of matrix $X_s$ also satisfy:
$$0<|\mu_1|<\dots<|\mu_\ell|,$$
as is shown in Figures~\ref{eigX6} and \ref{eigX7}, in the cases $s=6$ and $s=7$, respectively. In such a case, the following result holds true.

\begin{theo}\label{gammastar}
$\rho^*$ is minimized by choosing 
\begin{equation}\label{gammabest}
\zeta = |\mu_1| \equiv \min_{\mu\in\sigma(X_s)} |\mu|,
\end{equation}
resulting in 
\begin{equation}\label{rostarmin}
\rho^*= \left.\frac{1}{2\zeta}\frac{|\mu_1-\zeta|^2}{|\mu_1|}\right|_{\zeta=|\mu_1|}.
\end{equation}
In such a case, one obtains: 
\begin{equation}\label{rostarmin1}
\rho^*= 1-\cos\phi_1 <1.
\end{equation}
\end{theo}
\proof For (\ref{gammabest})-(\ref{rostarmin}), see \cite{BM02}. Concerning (\ref{rostarmin1}), one has:
\begin{eqnarray*}
\rho^* &=&\frac{1}{2|\mu_1|}\frac{|\mu_1-|\mu_1||^2}{|\mu_1|} ~=~ \frac{|\mu_1|^2\left[(1-\cos\phi_1)^2 +(\sin\phi_1)^2\right]}{2|\mu_1|^2}\\ &=&\frac{1+(\cos\phi_1)^2 +(\sin\phi_1)^2-2\cos\phi_1}2 ~=~\frac{2-2\cos\phi_1}2\\
&=& 1-\cos\phi_1.\,\QED
\end{eqnarray*}
Consequently, the blended implementation of HBVM$(k,s)$ methods is {\em always} $A$-convergent and, therefore, $L$-convergent. We can also characterize the speed of convergence when $q\approx0$, by considering that, from Theorem~\ref{muXs}, Corollary~\ref{corrho}, and Theorem~\ref{gammastar}, it follows that
$$\rho(q) = \frac{|q|\left|\mu_1-|\mu_1|\right|^2}{|\mu_1|\,\left|1-q|\mu_1|\right|^2} = \frac{\left|\mu_1-|\mu_1|\right|^2}{|\mu_1|}|q| +O(|q|^2) \approx \tilde\rho |q|,$$
where the parameter $$\tilde\rho = \frac{\left|\mu_1-|\mu_1|\right|^2}{|\mu_1|}$$ is called  the {\em non-stiff} amplification factor.
In Table~\ref{fact} we list the relevant information for the iteration of HBVM$(k,s)$ methods.

\begin{table}
\caption{Bended iteration of HBVM$(k,s)$ methods.}
\label{fact}
\begin{center}
\begin{tabular}{|c|c|c|c|}
\hline
\hline
$s$ &$\zeta$ & $\rho^*$ & $\tilde\rho$\\
\hline
2 & 0.2887  &  0.1340  &  0.0774\\
3 & 0.1967  &  0.2765  &  0.1088\\
4 & 0.1475  &  0.3793  &  0.1119\\
5 & 0.1173  &  0.4544  &  0.1066\\
6 & 0.0971  &  0.5114  &  0.0993\\
7 & 0.0827  &  0.5561  &  0.0919\\
\hline
\hline
\end{tabular}
\end{center}
\end{table}

\section{Actual blended implementation}
Let us now sketch the blended implementation of HBVMs, when applied to a general, nonlinear system, also analyzing its complexity. In the case of the initial value problem $$y'=f(y), \qquad y(0)=y_0\in\RR^m,$$ the previous aguments can be generalized in a straightworfard way, by considering that now the weighting function becomes
$$\theta = I_s\otimes \Omega^{-1}, \qquad \Omega = I_m-h\zeta J_0,$$ where $h$ is the stepsize, $\zeta$ is the optimal parameter specified in the second column in Table~\ref{fact}, and $J_0$ is the Jacobian of $f$ evaluated at $y_0$ (clearly, we are speaking about the first step in the numerical integration).

 From (\ref{Fgamma}) and (\ref{NewtC2}), we have to solve the {\em outer-inner} iteration described in Table~\ref{oii}. Let us analyze its computational complexity (let $m$ denote the dimension of the continuous problem and $e\in\RR^s$ be the unit vector), by considering, for each item, only the leading term in the complexity and denoting, as 1 {\em flop}, an elementary (binary) algebraic {\em flo}ating-point {\em op}eration. One obtains:
\begin{itemize}

\item $\Omega$: 1 Jacobian evaluation;

\item $\theta$: $\frac{2}3m^3$ flops for the $LU$ factorization of $\Omega$;
\item $y^\ell$: $2ksm$ flops;

\item $f^\ell$: $k$ function evaluations;

\item $\eta^\ell$: $2ksm$ flops;

\item $z^{\ell,r}$: $2sm^2$ flops;

\item $u^{\ell,r}$: $2s^2m$ flops;

\item $w^{\ell,r}$: $2s^2m$ flops;

\item $\Delta^{\ell,r+1}$: $4sm^2$ flops;

\item $\hat\gamma^{\ell+1}$: $sm$ flops.

\end{itemize}
Consequenly, this algorithm has a fixed computational cost of 1 Jacobian evaluation and $\frac{2}3m^3$ flops, plus, assuming that $\nu$ {\em inner} iterations are performed, a cost of $k$ function evaluations and $4ksm +\nu(6sm^2 +4s^2m) + sm$ flops per {\em outer} iteration.

A simplified (and sometimes more efficient) procedure is that of performing a {\em nonlinear} iteration, obtained by performing exactly 1 inner iteration (i.e., that with $r=0$) in the above procedure, thus obtaining the algorithm depicted in Table~\ref{nli}. In such a case, the resulting computational cost is obtained as follows:
\begin{itemize}

\item $\Omega$: 1 Jacobian evaluation;

\item $\theta$: $\frac{2}3m^3$ flops for the $LU$ factorization of $\Omega$;

\item $y^\ell$: $2ksm$ flops;

\item $f^\ell$: $k$ function evaluations;

\item $\eta^\ell$: $2ksm$ flops;

\item $u^\ell$: $2s^2m$ flops;

\item $\Delta^\ell$: $4sm^2$ flops;

\item $\hat\gamma^{\ell+1}$: $sm$ flops.

\end{itemize}
Consequenly, this latter algorithm has a fixed computational cost of 1 Jacobian evaluation and $\frac{2}3m^3$ flops, plus a cost of $s$ function evaluations and $4sm^2 +4ksm+ sm$ flops per iteration.

\begin{table}[t]
\caption{Outer-inner iteration for the blended implementation of HBVMs.}
\label{oii}
\begin{eqnarray*}
\Omega &=& I-(h\zeta) J_0\\
\theta &=& I_s\otimes\Omega^{-1} \qquad\qquad \% \mbox{~actually}, \Omega \mbox{~is factored} ~LU\\
\hat\gamma^0&&given \qquad\qquad \% \mbox{~e.g.,~} \hat\gamma^0=0\\
for&&\ell~ =~0,1,\dots\\
&&y^\ell = e\otimes y_0 + (h\II_s)\otimes I\, \hat\gamma^\ell\\
&&f^\ell = f(y^\ell)\\
&&\eta^\ell =   -\hat\gamma^\ell + (\P_s^T\Omega)\otimes I\, f^\ell \qquad\qquad \%~ F(\hat\gamma^\ell)\\
&&\Delta^{\ell,0} = 0\\
&&for \quad r~=~0,1,\dots\\
&&\qquad if \quad r>0\\
&&\qquad\quad z^{\ell,r} = [I_s\otimes J_0]\Delta^{\ell,r}\\
&&\qquad\quad u^{\ell,r} = [(\zeta X_s^{-1})\otimes I](\Delta^{\ell,r}+\eta^\ell) -(h\zeta) z^{\ell,r}\\
&&\qquad\quad w^{\ell,r} =  \Delta^{\ell,r} +\eta^\ell -[(hX_s)\otimes I]z^{\ell,r}\\
&&\qquad else\\
&&\qquad\quad u^{\ell,0} = [(\zeta X_s^{-1})\otimes I]\eta^\ell\\
&&\qquad\quad w^{\ell,0} =  \eta^\ell\\
&&\qquad end\\
&&\qquad \Delta^{\ell,r+1} = \Delta^{\ell,r} - \theta\left[ u^{\ell,r}+\theta (w^{\ell,r}-u^{\ell,r}) \right]\\
&&end \qquad \Rightarrow \qquad \mbox{returns}\quad \Delta^\ell\\
&&\hat\gamma^{\ell+1} = \hat\gamma^\ell+\Delta^\ell\\
end
\end{eqnarray*}
\bigskip

\caption{Nonlinear iteration for the blended implementation of HBVMs.}
\label{nli}
\begin{eqnarray*}
\Omega &=& I-(h\zeta) J_0\\
\theta &=& I_s\otimes\Omega^{-1} \qquad\qquad \% \mbox{~actually}, \Omega \mbox{~is factored}~LU\\
\hat\gamma^0&&given \qquad\qquad \% \mbox{~e.g.,~} \hat\gamma^0=0\\
for&&\ell~ =~0,1,\dots\\
&&y^\ell = e\otimes y_0 + (h\II_s)\otimes I\, \hat\gamma^\ell\\
&&f^\ell = f(y^\ell)\\
&&\eta^\ell =   -\hat\gamma^\ell + (\P_s^T\Omega)\otimes I\, f^\ell \qquad\qquad \%~ -F(\hat\gamma^\ell)\\
&&u^\ell = [(\zeta X_s^{-1})\otimes I] \eta^\ell\\
&&\Delta^\ell = \theta\left[\theta (u^\ell-\eta^\ell)-u^\ell \right]\\
&&\hat\gamma^{\ell+1} = \hat\gamma^\ell+\Delta^\ell\\
end
\end{eqnarray*}
\end{table}

\chapter{Line Integral Methods}
Sometimes conservative problems are not in Hamiltonian form and/or they posses multiple constants of motions, which are functionally independent. In certain cases, it is crucial to be able to preserve all of them, in order to obtain a faithfully simulation of the underlying dynamical system. For this reason, we extend the polynomial methods studied above, in order to cope with these, more general, conservative problems. The material of this chapter is based on \cite{BI12,BCMR12}.

\section{Introduction}
In the previous chapters, we have studied polynomial methods for approximately solving, on the interval $[0,h]$, the initial value problem
\begin{eqnarray}\label{ivpode}
y'(ch) &=& f(y(ch)) ~\equiv~ \sum_{j\ge 0} \gamma_j(y) P_j(c), \qquad c\in[0,1],\qquad y(0)=y_0\in\RR^{2m},\\
\gamma_j(y) &=& \int_0^1P_j(\tau)f(y(\tau h))\dd\tau, \qquad j\ge0. \label{gamj1}
\end{eqnarray}
The methods that we have considered are characterized by a suitable polynomial $\sigma\in\Pi_s$ such that
\begin{equation}\label{polimet}
\sigma'(ch) = \sum_{j=0}^{s-1} \gamma_j(\sigma) P_j(c),\qquad c\in[0,1], \qquad \sigma(0)=y_0,
\end{equation}
then approximating the integrals $\gamma_j(\sigma)$ by a suitable quadrature formula. When the problem is Hamiltonian, that is, $f(\cdot)=J\nabla H(\cdot)$, with $J^T=-J=J^{-1}$, energy is conserved, for the discrete-time dynamical system defined by (\ref{polimet}).
Indeed,
\begin{eqnarray*}
H(\sigma(h))-H(\sigma(0)) &=&\int_0^h \nabla H(\sigma(t))^T\sigma'(t)\dd t ~=~h\int_0^1 \nabla H(\sigma(\tau h))^T\sigma'(\tau h)\dd\tau\\
&=&h\int_0^1 \nabla H(\sigma(\tau h))^T\sum_{j=0}^{s-1}\gamma_j(\sigma)P_j(\tau)\dd\tau\\
&=&h\sum_{j=0}^{s-1}\left[\int_0^1 \nabla H(\sigma(\tau h))P_j(\tau)\dd\tau\right]^T\gamma_j(\sigma)\\
&=&h\sum_{j=0}^{s-1} \gamma_j(\sigma)^TJ\gamma_j(\sigma)~=~0,
\end{eqnarray*}
due to the fact that $J$ is skew-symmetric.

Let now consider the case where problem (\ref{ivpode}) is a general {\em conservative} (not necessarily Hamiltonian) problem, whose dimension will be denoted by $m$, for sake of brevity. Assume that it possess a set of $\nu$ smooth, functionally independent  (clearly, $\nu<m$), invariants. That is, there exists $$L:\,\RR^m\rightarrow\RR^\nu$$ such that
\begin{equation}\label{nablaL}
\nabla L(y)^Tf(y) = 0, \qquad \forall y\in\RR^m,
\end{equation}
where $\nabla L(y)^T$ denotes the Jacobian matrix of $L$. In such a case, along the solution of (\ref{ivpode}), one has:
\begin{equation}\label{Lyt}
\frac{\dd}{\dd t} L(y(t)) ~=~ \nabla L(y(t))^Ty'(t) ~=~ \nabla L(y(t))^Tf(y(t)) ~=~ 0\in\RR^\nu,
\end{equation}
because of (\ref{nablaL}). Consequently,
$$L(y(t))\equiv L(y_0),\qquad \forall t\ge0.$$
 From (\ref{ivpode}) and (\ref{Lyt}) one obtains:
\begin{eqnarray}\nonumber
L(y(h))-L(y(0)) &=&\int_0^h\nabla L(y(t))^Tf(y(t))\dd t \\  \nonumber
   &=& h\int_0^1 \nabla L(y(ch))^T\sum_{j\ge0} P_j(c)\gamma_j(y)\dd c\\
   &=& h\sum_{j\ge0} \phi_j(y)^T\gamma_j(y) ~=~0, \label{orto}
\end{eqnarray}
where $\gamma_j(y)$ is formally still defined by (\ref{gamj1}) and, moreover, we have set
\begin{equation}\label{phij}
\phi_j(y) = \int_0^1 P_j(c) \nabla L(y(ch))\dd c, \qquad j\ge0.
\end{equation}
\begin{rem} It is important noticing that, because of (\ref{nablaL}), (\ref{orto}) continues to hold, if we replace $y(t)$ by any other path $\sigma(t)$.\end{rem}

We observe that, because of the result of Lemma~\ref{hj}, one has (assuming, for sake of simplicity, that both $L(y(t))$ and $f(y(t))$ can be expanded in Taylor series at $t=0$):
\begin{equation}\label{phijh}
\phi_j(y) = O(h^j), \qquad j\ge0,
\end{equation}
However, if we replace the infinite series in (\ref{ivpode}) with the finite sum in (\ref{polimet}), one obtains, because of (\ref{nablaL}) and repeating similar steps as above:
\begin{equation}\label{h2s1}
L(\sigma(h)) - L(\sigma(0)) = h\sum_{j=0}^{s-1} \phi_j(\sigma)^T\gamma_j(\sigma)
                                          = -h\sum_{j\ge s} \phi_j(\sigma)^T\gamma_j(\sigma)
                                          = O(h^{2s+1}),
\end{equation}
since (see (\ref{orto})) $$\sum_{j\ge0} \phi_j(\sigma)^T\gamma_j(\sigma) ~=~0.$$ In order to get conservation for the polynomial dynamical system, we perturb (\ref{polimet}) as follows:
\begin{equation}\label{polimet1}
\sigma'(ch) = \sum_{j=0}^{s-1} \gamma_j(\sigma) P_j(c) ~-~\phi_0(\sigma)\aa,\qquad c\in[0,1], \qquad \sigma(0)=y_0,
\end{equation}
where $\phi_0(\sigma)$ is defined according to (\ref{phij}), and $\aa\in\RR^\nu$ is determined in order to enforce the conservation of the invariants. By repeating similar steps as above, we obtain:
\begin{eqnarray}\label{Lcons}
0 &=& L(\sigma(h))-L(\sigma(0)) ~=~h\int_0^1 \nabla L(\sigma(ch))^T\sigma'(ch)\dd c\\ \nonumber
   &=& h\sum_{j=0}^{s-1}\phi_j(\sigma)^T\gamma_j(\sigma) ~-~h\left[\phi_0(\sigma)^T\phi_0(\sigma)\right]\aa.
\end{eqnarray}
Consequently, conservation is gained, provided that
\begin{equation}\label{aa}
\left[\phi_0(\sigma)^T\phi_0(\sigma)\right]\aa = \sum_{j=0}^{s-1}\phi_j(\sigma)^T\gamma_j
\end{equation}
The following result holds true.

\begin{theo}\label{aatheo} The vector $\aa$ exists and is unique, for all sufficiently small step sizes $h$ and, moreover,
\begin{equation}\label{aah}
\aa = O(h^{2s}).
\end{equation}
\end{theo}
\proof
If the invariants are functionally independent, $\sigma(0)$ is a regular point for the constraints, so that $\nabla L(\sigma(0))$ has full column rank (i.e., $\nu$). Considering that
$$\phi_0(\sigma) = \int_0^1 \nabla L(\sigma(ch))\dd c\rightarrow \nabla L(y_0), \qquad\mbox{as}\qquad h\rightarrow0,$$
one has that matrix
$$M_0\equiv \left[\phi_0(\sigma)^T\phi_0(\sigma)\right]$$ is symmetric and positive definite and, therefore, nonsingular. The existence and uniqueness of $\aa$ then follows from the Implicit Function Theorem. Moreover, since (see (\ref{phijh}))
$$M_0 = O(h^0),$$ then
$$\aa = M_0^{-1}\sum_{j=0}^{s-1}\phi_j(\sigma)^T\gamma_j(\sigma) = -M_0^{-1}\sum_{j\ge s}\phi_j(\sigma)^T\gamma_j(\sigma) = O(h^{2s}).$$This completes the proof.\,\QED\bigskip

We now consider the following question: i.e., the polynomial $\sigma$ as defined in (\ref{polimet}) doesn't satisfy, in general, $L(\sigma(h))=L(y_0)$, even though $\sigma(h)-y(h)=O(h^{2s+1})$. Conversely, for the polynomial $\sigma$ defined by (\ref{polimet1})-(\ref{aa})  one has $L(\sigma(h))=L(y_0)$. Moreover, next theorem states that its order of convergence remains the same.

\begin{theo}\label{ordpolimet1} Let $\sigma$ be defined according to (\ref{polimet1})-(\ref{aa}). Then $\sigma(h)-y(h)=O(h^{2s+1})$.
\end{theo}
\proof By using similar steps as those used in the proof of Theorem~\ref{ord1}, one has:
\begin{eqnarray*}
\sigma(h)-y(h)&=& y(h;h,\sigma(h))-y(h;0,\sigma(0)) ~=~ \int_0^h \frac{\dd}{\dd t} y(h;t,\sigma(t))\dd t\\
&=& \int_0^h \left(\frac{\partial}{\partial \theta} y(h;\theta,\sigma(t))\Big|_{\theta=t} + \frac{\partial}{\partial \omega}y(h;t,\omega)\Big|_{\omega=\sigma(t)}\sigma'(t)\right)\dd t\\
&=& \int_0^h \Phi(h,t)[-f(\sigma(t)) +\sigma'(t)]\dd t \\
&=& h\int_0^1 \Phi(h,\tau h)[-f(\sigma(\tau h)) +\sigma'(\tau h)]\dd\tau \\
&=&-h\int_0^1 \Phi(h,\tau h)\left[ \sum_{j\ge0}P_j(\tau)\gamma_j(\sigma) - \sum_{j=0}^{s-1}P_j(\tau)\gamma_j(\sigma) + \phi_0(\sigma)\aa\right]\dd\tau\\
&=&-h\int_0^1 \Phi(h,\tau h)\sum_{j\ge s}P_j(\tau)\gamma_j(\sigma)\dd\tau -h\int_0^1\Phi(h,\tau h)\phi_0(\sigma)\aa\dd\tau \\
&=&
-h\sum_{j\ge s}\underbrace{\left[\int_0^1 \overbrace{\Phi(h,\tau h)}^{\equiv G(\tau h)}\, P_j(\tau)\dd\tau\right]}_{=\,O(h^j)}\overbrace{\gamma_j(\sigma)}^{=\,O(h^j)} -h\underbrace{\int_0^1\Phi(h,\tau h)\dd\tau \phi_0(\sigma)}_{=\,O(1)}\overbrace{\aa}^{=\,O(h^{2s})}\\
&=& h\sum_{j\ge s} O(h^{2j})~+O(h^{2s+1})~=~O(h^{2s+1}).\,\QED
\end{eqnarray*}

\section{Discretization}
As was previously observed,  (\ref{polimet1})-(\ref{aa}) doesn't yet define a method, but a conservative {\em formula}. As matter of fact, a numerical method is obtained when the integrals (see (\ref{gamj1}) and (\ref{phij}))
$$\gamma_j(\sigma), \qquad \phi_j(\sigma), \qquad j=0,\dots,s-1,$$
are approximated by means of a suitable quadrature formula. In principle, they could be approximated by means of different quadrature formulae:
\begin{itemize}
\item one formula, based at the abscissae $0\le c_1<\dots<c_k\le 1$ and corresponding weights $\{b_i\}$, of order $q$, for approximating $\gamma_j(\sigma)$:
\begin{equation}\label{formulak}
\gamma_j(\sigma) = \sum_{\ell=1}^k b_i P_j(c_i) f(\sigma(c_ih)) -\Delta_j(h), \qquad i=0,\dots,s-1,
\end{equation}
with
\begin{equation}\label{Djh}
\Delta_j(h) = O(h^{q-j}),\qquad j=0,\dots,s-1;
\end{equation}

\item another formula, based at the abscissae $0\le \tau_1<\dots<\tau_r\le 1$ and corresponding weights $\{\bb_i\}$, of order $\hat{q}$, for approximating $\phi_j(\sigma)$:
\begin{equation}\label{formular}
\phi_j(\sigma) = \sum_{\ell=1}^{r} \bb_i P_j(\tau_i) \nabla L(\sigma(\tau_ih)) -\Psi_j(h), \qquad i=0,\dots,s-1,
\end{equation}
with
\begin{equation}\label{Psijh}
\Psi_j(h) = O(h^{\hat{q}-j}),\qquad j=0,\dots,s-1.
\end{equation}
\end{itemize}
In the following, for sake of simplicity, we shall consider the following choices of such abscissae:
\begin{eqnarray}\label{cik}
P_k(c_i)&=&0, \qquad i=1,\dots,k, \qquad \Rightarrow \qquad q = 2k,\\
P_r(\tau_j)&=&0, \qquad j=1,\dots,r, \qquad \Rightarrow \qquad \hat{q} = 2r.\label{taur}
\end{eqnarray}
\begin{defi}
We shall refer to such a method as {\em LIM$(r,k,s)$}, where LIM is the acronym for {\em Line Integral Method}.
\end{defi}
\begin{rem} We observe that:
\begin{itemize}
\item LIM$(0,s,s)$ is the $s$-stage Gauss method,
\item LIM$(0,k,s)$ is the HBVM$(k,s)$ method,
\end{itemize}
where $r=0$ means that no invariant conservation is seeked.
\end{rem}

After discretization, the polynomial $\sigma$ is obviously formally replaced by the polynomial $u\in\Pi_s$ such that:
\begin{equation}\label{polimet2}
u'(ch) = \sum_{j=0}^{s-1} \hat{\gamma}_j P_j(c) ~-~\hat{\phi}_0\hat\aa,\qquad c\in[0,1], \qquad u(0)=y_0,
\end{equation}
with, in general, (see (\ref{gamj1}), (\ref{phij}), (\ref{formulak})--(\ref{Psijh}))
\begin{eqnarray}\label{hgj}
\hat\gamma_j &=& \sum_{i=1}^k b_i P_j(c_i) f(u(c_ih)) ~\equiv~\gamma_j(u) + \Delta_j(h),\\
\hat\phi_j       &=& \sum_{\ell=1}^r \bb_\ell P_j(\tau_\ell)\nabla L(u(\tau_\ell h)) ~\equiv~\phi_j(u) + \Psi_j(h),
\label{hfj}
\end{eqnarray}
for $j=0,\dots,s-1$, and (see (\ref{cik})-(\ref{taur}))
\begin{eqnarray}\label{haa}\nonumber
\lefteqn{\left[\hat\phi_0^T\hat\phi_0\right]\hat\aa = \sum_{j=0}^{s-1}\hat\phi_j^T\hat\gamma_j
~=~ \sum_{j=0}^{s-1}\left[\phi_j(u)+\Psi_j(h)\right]^T\left[\gamma_j(u)+\Delta_j(h)\right]}\\ \nonumber
&=& \sum_{j=0}^{s-1}\left[\overbrace{\underbrace{\phi_j(u)^T}_{=\,O(h^j)}\underbrace{\gamma_j(u)}_{=\,O(h^j)}}^{=\,O(h^{2j})} +
\overbrace{\underbrace{\Psi_j(h)^T}_{=\,O(h^{2r-j})}\gamma_j(u)}^{=\,O(h^{2r})}  +\overbrace{\phi_j(u)^T\underbrace{\Delta_j(h)}_{=\,O(h^{2k-j})}}^{=\,O(h^{2k})}
+\overbrace{\Psi_j(h)^T\Delta_j(h)}^{=\,O(h^{2(r+k-j)})}\right]\\ \nonumber
&=& -\sum_{j\ge s}\phi_j(u)^T\gamma_j(u) + \sum_{j=0}^{s-1}\left[
\Psi_j(h)^T\gamma_j(u)  +\phi_j(u)^T\Delta_j(h)+\Psi_j(h)^T\Delta_j(h)\right]\\
&=& O(h^{2s}) + O(h^{2r}) + O(h^{2k}) + O(h^{2(r+k-s+1)}).\label{orddis}
\end{eqnarray}
\begin{rem}\label{costlim}
From (\ref{hgj})-(\ref{hfj}), and (\ref{orddis}), it follows that the (block) size of the discrete problem is $2s+1$:
\begin{itemize}
\item the $s$ coefficients $\hat\gamma_j\in\RR^m$, $j=0,\dots,s-1$,
\item the $s$ coefficients $\hat\phi_j\in\RR^{m\times \nu}$, $j=0,\dots,s-1$,
\item the vector $\hat\aa\in\RR^\nu$.
\end{itemize}
even though it must be stressed that $\hat\phi_j$ is a matrix with $\nu$ columns, whereas $\hat\gamma_j$ is a vector.
The efficient implementation of such methods is, however, still under investigation.
\end{rem}
The following result then holds true.

\begin{theo}\label{ordlim}
For all $r\ge s$ and $k\ge s$, for LIM$(r,k,s)$ one obtains: $$u(h)-y(h)=O(h^{2s+1}).$$
\end{theo}
\proof
From (\ref{orddis}) and the hypotheses $r\ge s$ and $k\ge s$, it follows that $$\hat\aa = O(h^{2s}).$$ Moreover, ones has, by repeating similar steps as in Theorem~\ref{ordpolimet1}:
\begin{eqnarray*}
\lefteqn{u(h)-y(h)~=~ y(h;h,u(h))-y(h;0,u(0)) ~=~ \int_0^h \frac{\dd}{\dd t} y(h;t,u(t))\dd t}\\
&=& \int_0^h \left(\frac{\partial}{\partial \theta} y(h;\theta,u(t))\Big|_{\theta=t} + \frac{\partial}{\partial \omega}y(h;t,\omega)\Big|_{\omega=u(t)}u'(t)\right)\dd t\\
&=& \int_0^h \Phi(h,t)[-f(u(t)) +u'(t)]\dd t \\
&=& h\int_0^1 \Phi(h,\tau h)[-f(u(\tau h)) +u'(\tau h)]\dd\tau \\
&=&-h\int_0^1 \Phi(h,\tau h)\left[ \sum_{j\ge0}P_j(\tau)\gamma_j(u) - \sum_{j=0}^{s-1}P_j(\tau)\hat\gamma_j + \hat\phi_0\hat\aa\right]\dd\tau\\
&=&-h\int_0^1 \Phi(h,\tau h)\left[ \sum_{j\ge0}P_j(\tau)\gamma_j(u) \right.\\
&&\qquad\qquad\qquad\left.-\sum_{j=0}^{s-1}P_j(\tau)\left[\gamma_j(u) -\overbrace{\Delta_j(h)}^{=\,O(h^{2k-j})}\right] + \left[\phi_0(u)-\overbrace{\Psi_0(h)}^{=\,O(h^{2r})}\right]\hat\aa\right]\dd\tau\\
&=&-h\int_0^1 \Phi(h,\tau h)\sum_{j\ge s}P_j(\tau)\gamma_j(u)\dd\tau -h\int_0^1\Phi(h,\tau h)\phi_0(u)\hat\aa\dd\tau \\
&& -h\int_0^1 \Phi(h,\tau h)\sum_{j=0}^{s-1}P_j(\tau)\Delta_j(h)\dd\tau +h\int_0^1\Phi(h,\tau h)\Psi_0(u)\hat\aa\dd\tau\\
&=&-h\sum_{j\ge s}\left[\underbrace{\int_0^1 \Phi(h,\tau h)P_j(\tau)\dd\tau}_{=\,O(h^j)}\right] \overbrace{\gamma_j(u)}^{=\,O(h^j)}-h\left[\underbrace{\int_0^1\Phi(h,\tau h)\dd\tau}_{=\,O(1)}\right]\underbrace{\phi_0(u)}_{=\,O(1)}\overbrace{\hat\aa}^{=\,O(h^{2s})}\\
&& -h\sum_{j=0}^{s-1}\left[\underbrace{\int_0^1 \Phi(h,\tau h)P_j(\tau)\dd\tau}_{=\,O(h^j)}\right]\overbrace{\Delta_j(h)}^{=\,O(h^{2k-j})} +h\left[\underbrace{\int_0^1\Phi(h,\tau h)\dd\tau}_{=\,O(1)}\right]\underbrace{\Psi_0(u)}_{=\,O(h^{2r})}\overbrace{\hat\aa}^{=\,O(h^{2s})}\\
&=& h\sum_{j\ge s} O(h^{2j})~+O(h^{2s+1})~+~O(h^{2k+1})~+O(h^{2(r+s)+1})~=~O(h^{2s+1}).\,\QED
\end{eqnarray*}

\bigskip
Consequently, the following statement easily follows.

\begin{cor}\label{ordlim1}
For all $r\ge s$ and $k\ge s$,  LIM$(r,k,s)$ has order $2s$.
\end{cor}

Concerning the conservations of the invariants, the following result holds true.

\begin{theo}\label{conslim} For given $r,k\ge s$, 
LIM$(r,k,s)$ exactly conserves polynomial invariants of degree
\begin{equation}\label{numax}
\nu\le\frac{2r}s.
\end{equation}
For all suitably regular invariants,
\begin{equation}\label{altricasi}
L(u(h))-L(y_0)=O(h^{2r+1}).
\end{equation}
\end{theo}
\proof One has, by virtue of (\ref{polimet2})--(\ref{hfj}):
\begin{eqnarray*}
L(u(h))-L(y_0) &=& L(u(h))-L(\sigma(0)) ~=~h\int_0^1 \nabla L(u(ch))^Tu'(ch)\dd c\\
   &=& h\sum_{j=0}^{s-1}\phi_j(u)^T\hat\gamma_j ~-~h\left(\phi_0(u)^T\hat\phi_0\right)\hat\aa\\
   &\equiv& E_L.
\end{eqnarray*}
In case $L$ is a polynomial of degree $\nu$ satisfying (\ref{numax}), the quadrature formula (\ref{hfj}) is exact, so that
$$\phi_j(u)=\hat\phi_j, \qquad j=0,\dots,s-1.$$
Consequently, $E_L=0$, by virtue of  (\ref{orddis}). This proves the first part of the thesis. In any other case, one has:
\begin{eqnarray*}
E_L&=&h\left[ \sum_{j=0}^{s-1}(\hat\phi_j-\Psi_j)^T\hat\gamma_j ~-~\left((\hat\phi_0-\Psi_0)^T\hat\phi_0\right)\hat\aa  \right]\\
&=&h \left[\underbrace{\sum_{j=0}^{s-1}\hat\phi_j^T\hat\gamma_j -\left(\hat\phi_0^T\hat\phi_0\right)\hat\aa}_{=0, \mbox{~\small see (\ref{orddis})}}~-~\sum_{j=0}^{s-1}\underbrace{\hat\Psi_j^T\hat\gamma_j}_{=\,O(h^{2r})} +\underbrace{(\Psi_0^T\hat\phi_0)}_{=\,O(h^{2r})}\overbrace{\hat\aa}^{=\,O(h^{2s})}
   \right]\\
   &=& O(h^{2r+1})~+~O(h^{2(r+s)+1})~=~O(h^{2r+1}),
\end{eqnarray*}
thus proving (\ref{altricasi}).\,\QED\bigskip

\begin{rem} From (\ref{altricasi}), one obtains that, for any suitably regular set of invariants, conservation  can {\em always} be practically obtained, provided that $r$ is large enough. Indeed, it is enough to obtain conservation up to roundoff errors.

Moreover, if some of the invariants are polynomials of low degree, then, in principle, a less accurate quadrature formula could be used for approximating the corresponding integrals (\ref{hfj}).
\end{rem}

The following result can be also proved, by using arguments similar to those used in Section~\ref{symmet}. 

\begin{theo} Provided that the abscissae (\ref{cik})-(\ref{taur}) are symmetrically distributed in the interval [0,1], LIM$(r,k,s)$ is a symmetric method.\end{theo}

We now provide a couple of straightforward fully-conserving generalizations of HBVMs and Gauss-Legendre Runge-Kutta methods.

\subsection{LIM$(k,k,s)$}
Such conserving methods can be regarded as a straightforward generalization of HBVM$(k,s)$ methods. In such a case, the same set of abscissae,
$$
0<c_1<\dots<c_k<1,\qquad P_k(c_i)=0, \qquad i=1,\dots,k,
$$
are used for the quadraures approximating $\gamma_j(u)$, $\phi_j(u)$, $j=0,\dots,s-1$. Consequently,
by setting as usual $Y_i=u(c_ih)$, one obtains:
\begin{eqnarray}\nonumber
Y_i &=& y_0 + h\sum_{j=0}^{s-1} \int_0^{c_i} P_j(x)\dd x\,\hat\gamma_j
~-~ \,h\hat\phi_0\hat\aa, \qquad i=1,\dots,k,\\  \nonumber
&&\hat\gamma_j ~=~ \sum_{\ell=1}^k b_\ell P_j(c_\ell) f(Y_\ell),\\
&&\qquad\qquad\qquad\qquad\qquad\qquad j=0,\dots,s-1,    \label{LIMkks}\\ \nonumber
&&\hat\phi_j       ~=~ \sum_{\ell=1}^k b_\ell P_j(c_\ell)\nabla L(u(Y_\ell)),\\ \nonumber
&&\left[\hat\phi_0^T\hat\phi_0\right]\hat\aa ~=~ \sum_{j=0}^{s-1}\hat\phi_j^T\hat\gamma_j,
\end{eqnarray}
with the new approximation given by
\begin{equation}\label{LIMkks1}
y_1=y_0+h\sum_{i=1}^k b_i f(Y_i) \,-h\hat\phi_0\hat\aa.\end{equation}

\subsection{LIM$(k,s,s)$} These methods turn out to be fully conservative variants of the $s$-stage Gauss methods. In such a case, two sets of abscissae are used, i.e.,
\begin{equation}\label{cis1}
0<c_1<\dots<c_s<1,\qquad P_s(c_i)=0, \qquad i=1,\dots,s,
\end{equation}
$$0<\tau_1<\dots<\tau_k<1, \qquad P_k(\tau_i)=1,\dots, s.$$
The resulting method can be easily seen to be, by setting $Y_i=u(c_ih)$:
\begin{eqnarray}\nonumber
Y_i &=& y_0 + h\sum_{j=0}^{s-1} \int_0^{c_i} P_j(x)\dd x\,\hat\gamma_j
~-~ \,h\hat\phi_0\hat\aa, \qquad i=1,\dots,s,\\  \nonumber
Z_\ell &=& \sum_{i=1}^s L_{is}(\tau_\ell)Y_i, \qquad \ell=1,\dots,k,\\ \nonumber
&&\hat\gamma_j ~=~ \sum_{\ell=1}^s b_\ell P_j(c_\ell) f(Y_\ell),\\
&&\qquad\qquad\qquad\qquad\qquad\qquad j=0,\dots,s-1,    \label{LIMkss}\\ \nonumber
&&\hat\phi_j       ~=~ \sum_{\ell=1}^k \bb_\ell P_j(c_\ell)\nabla L(Z_\ell),\\ \nonumber
&&\left[\hat\phi_0^T\hat\phi_0\right]\hat\aa ~=~ \sum_{j=0}^{s-1}\hat\phi_j^T\hat\gamma_j,
\end{eqnarray}
where the $L_{is}(\tau)$ are the Lagrange polynomials defined at the abscissae (\ref{cis1}), and
the new approximation is given by
\begin{equation}\label{LIMkss1}
y_1=y_0+h\sum_{i=1}^s b_i f(Y_i) \,-h\hat\phi_0\hat\aa.\end{equation}
In this case, the first equation in (\ref{LIMkss}) can be recast in the more usual form,
$$u'(c_ih) = f(u(c_ih)) - \hat\phi_0\hat\aa, \qquad i=1,\dots,k,$$
which emphasizes the connection with the collocation conditions of the original $s$-stage Gauss method.

\section{Numerical tests}

In this section, we report a few numerical tests on a couple of conservative problems, possessing multiple invariants.

\subsection*{The Kepler problem}
A noticeable example of Hamiltonian problem with multiple invariants, is the Kepler problem, defined by the (non-polynomial) Hamiltonian:
\begin{equation}\label{keplH}
H(\bfq,\bfp) = \frac{1}2\|\bfp\|_2^2 -\frac{1}{\|\bfq\|_2}, \qquad \bfq,\bfp\in\RR^2,
\end{equation}
It admits the following invariants of motions, besides the Hamiltonian:
\begin{itemize}
\item the {\em angular momentum}, \begin{equation}\label{angm}
L(\bfq,\bfp) = q_1p_2-q_2p_1
\end{equation}
which is a quadratic invariant;

\item the so called {\em Laplace-Runge-Lenz (LRL)} vector which, for the problem at hand, implies the conservation of the following quantity,
\begin{equation}\label{lrl}
F(\bfq,\bfp) = q_2p_1^2 - q_1p_1p_2 -\frac{q_2}{\|\bfq\|_2}.
\end{equation}
\end{itemize}
When starting at the initial point
\begin{equation}\label{kepl0}
q_1 = 1-\eps, \qquad q_2 = p_1 = 0, \qquad p_2=\sqrt{\frac{1+\eps}{1-\eps}},
\end{equation}
it has a periodic orbit of period $T=2\pi$, which is, in the $(q_1,q_2)$-plane, an ellipse of eccentricity $\eps$.
We first consider the integration of problem (\ref{kepl0}) with $\eps = 0.6$, by using the following fourth order methods with a  constant stepsize $h=10^{-2}\pi$:
\begin{itemize}
\item LIM(0,2,2), (i.e., the 2-stage Gauss method);

\item LIM(0,8,2), (i.e., HBVM(8,2));

\item LIM(8,2,2), (i.e., the ``fully conservative'' variant of the 2-stage Gauss method);

\item LIM(8,8,2), (i.,e., the ``fully conservative'' variant of the HBVM(8,2) method).
\end{itemize}
Indeed, for the ``fully conserving methods'', $r=8$ is enough to obtain a {\em practical} conservation of all invariants.

\begin{figure}
\centerline{\includegraphics[width=12cm,height=7cm]{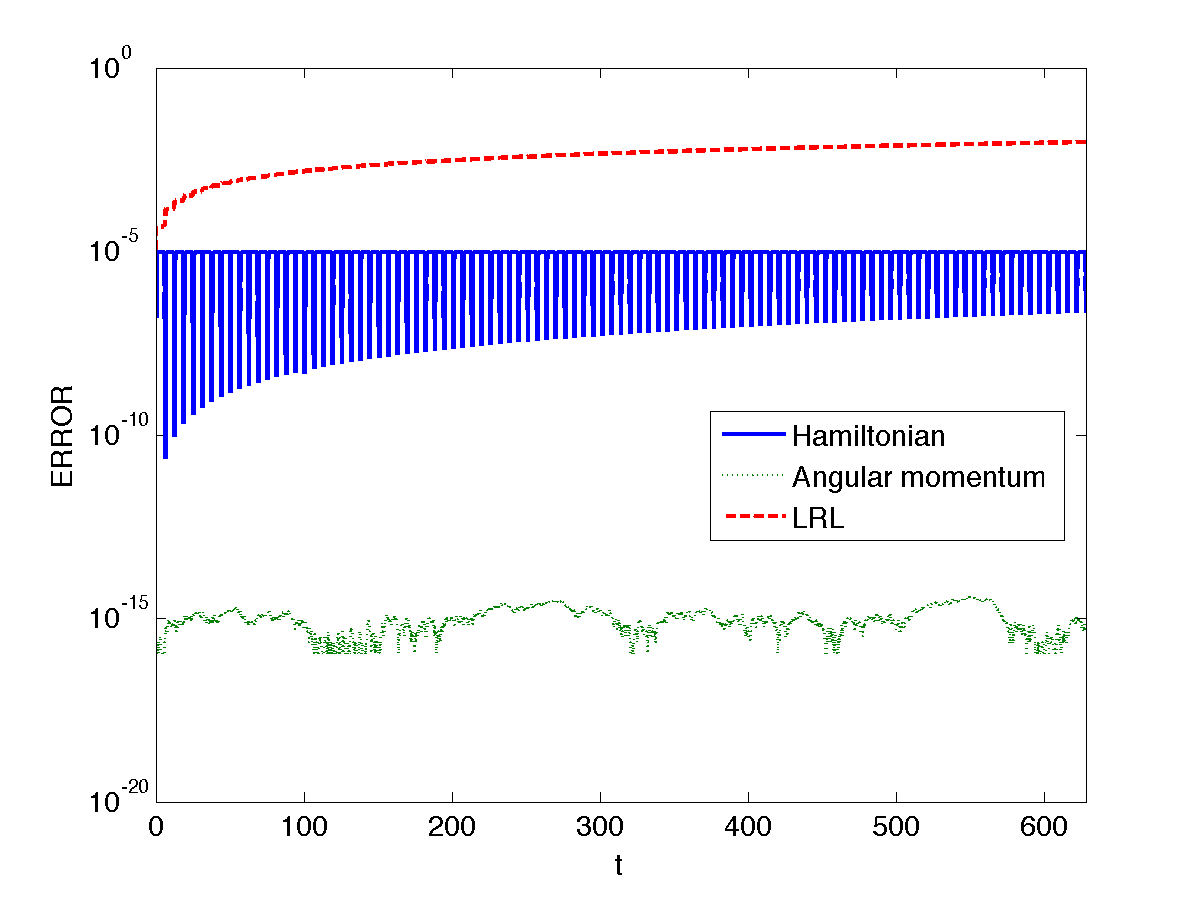}}
\caption{Kepler problem, $\eps=0.6$, 2-stage Gauss method, $h=10^{-2}\pi$, invariants errors.}
\label{kepl022}

\bigskip
\centerline{\includegraphics[width=12cm,height=7cm]{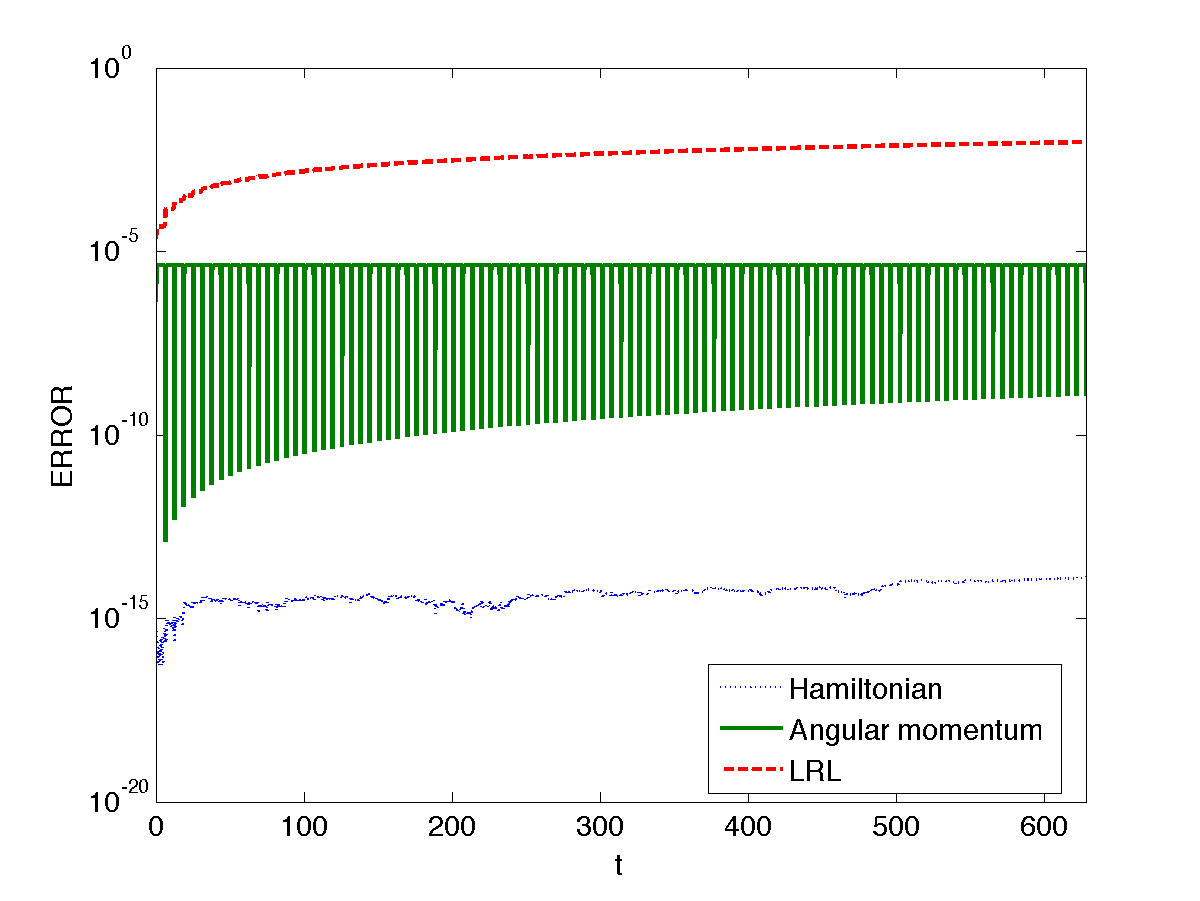}}
\caption{Kepler problem, $\eps=0.6$, HBVM(8,2) method, $h=10^{-2}\pi$, invariants errors.}
\label{kepl082}
\end{figure}

\begin{figure}
\bigskip
\centerline{\includegraphics[width=12cm,height=7cm]{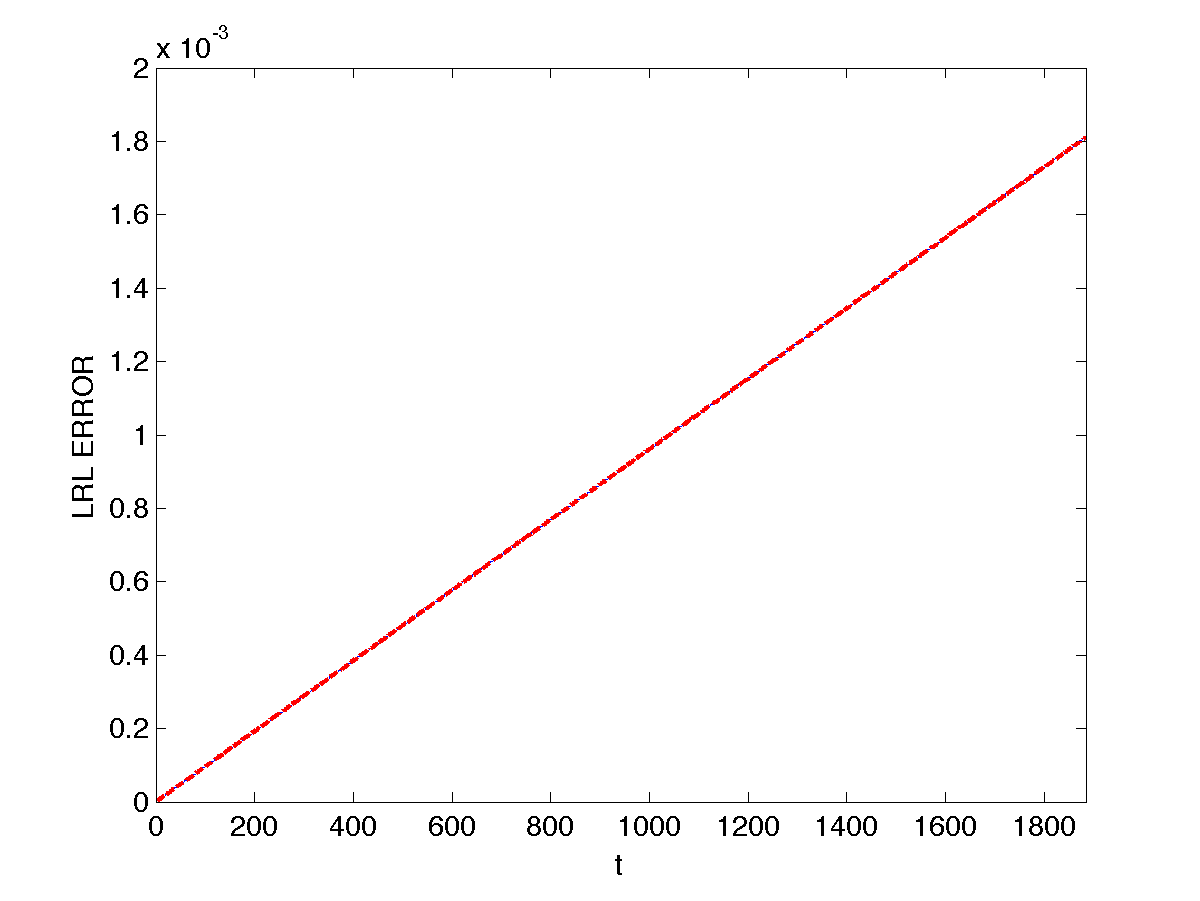}}
\caption{Kepler problem, $\eps=0.6$, drift in the LRL invariant for both the 2-stage Gauss method and the HBVM(8,2) method, $h=10^{-2}\pi$.}
\label{LRL_error}

\bigskip
\centerline{\includegraphics[width=12cm,height=7cm]{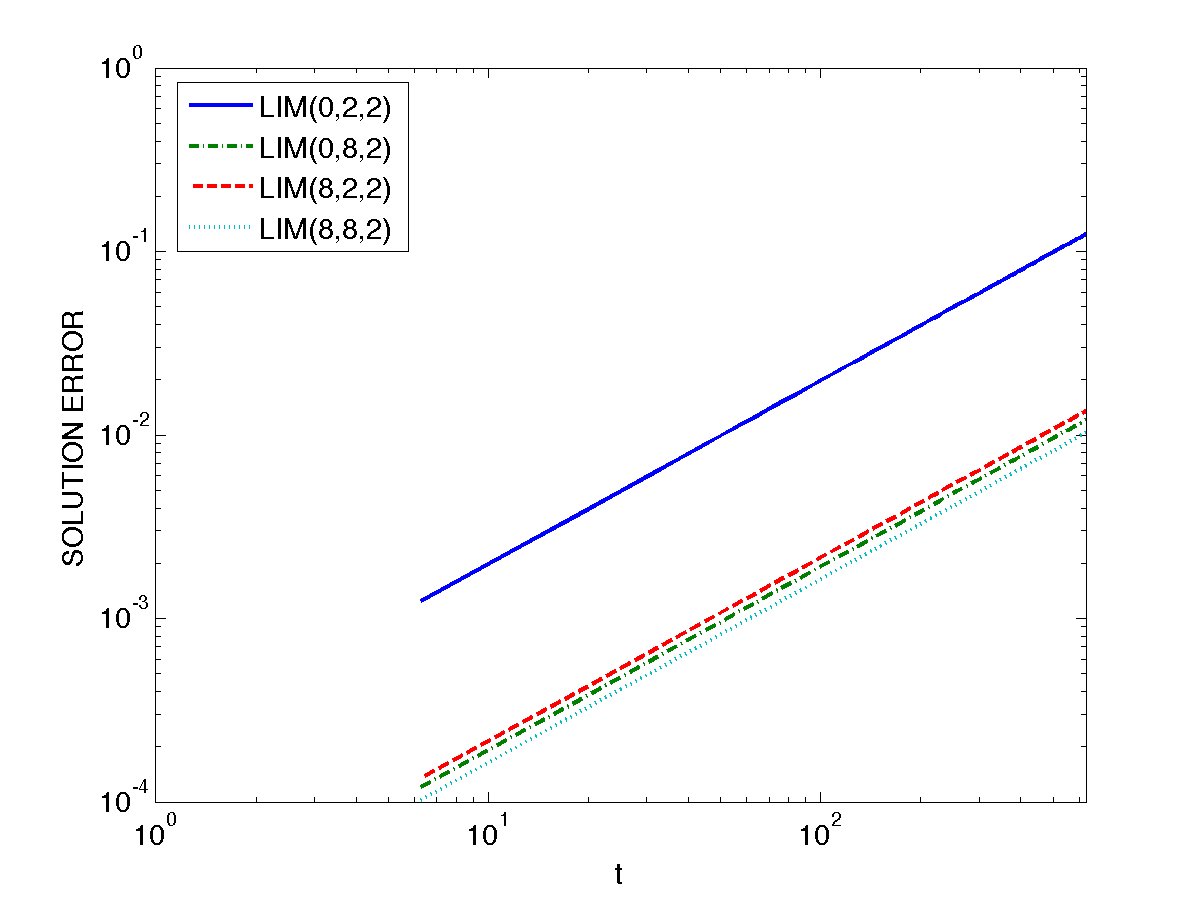}}
\caption{Kepler problem, $\eps=0.6$, error in the numerical solution, $h=10^{-2}\pi$.}
\label{keplerr}
\end{figure}

In Figures~\ref{kepl022} and \ref{kepl082} is the plot of the error in the invariants for LIM(0,2,2) ans LIM(0,8,2), respectively: for the first method, the error in the Hamiltonian is bounded, the angular momentum (which is a quadratic invariant) is conserved up to roundoff, and a drift apparently occurs in the LRL invariant; for the second method the situation is similar, with the roles of the Hamiltonian and of the angular momentum exchanged each other. Also in this case, a drift seems to occur in the LRL invariant: this is better evidenced in Figure~\ref{LRL_error}, where the drift of the LRL invariant (\ref{lrl}) is plotted over a longer interval, and it appears to be same for both LIM(0,2,2) and LIM(0,8,2). 

On the other hand, both LIM(8,2,2) and LIM(8,8,2) conserve all the invariants up to roundoff. In Figure~\ref{keplerr}, there is the plot of the error in the numerical solution (measured at each period) for all methods: it is evident that its growth is linear, and with the same order, even though LIM(0,2,2) (i.e., the two stage, symplectic Gauss method) is less accurate than the other methods.


\begin{figure}
\centerline{\includegraphics[width=12cm,height=7cm]{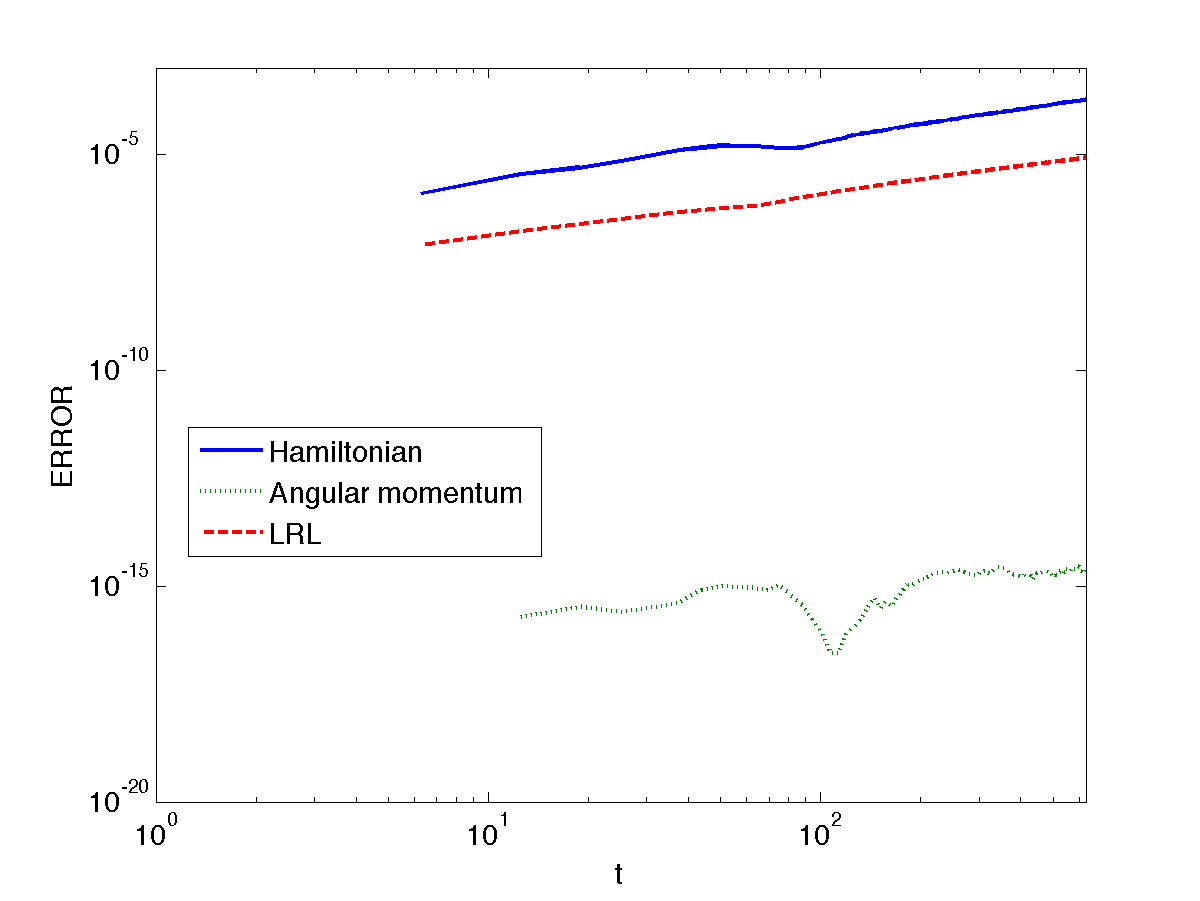}}
\caption{Kepler problem, $\eps=0.99$, 2-stage Gauss method, $tol=10^{-8}$, invariants errors.}
\label{kepl022v}

\bigskip
\centerline{\includegraphics[width=12cm,height=7cm]{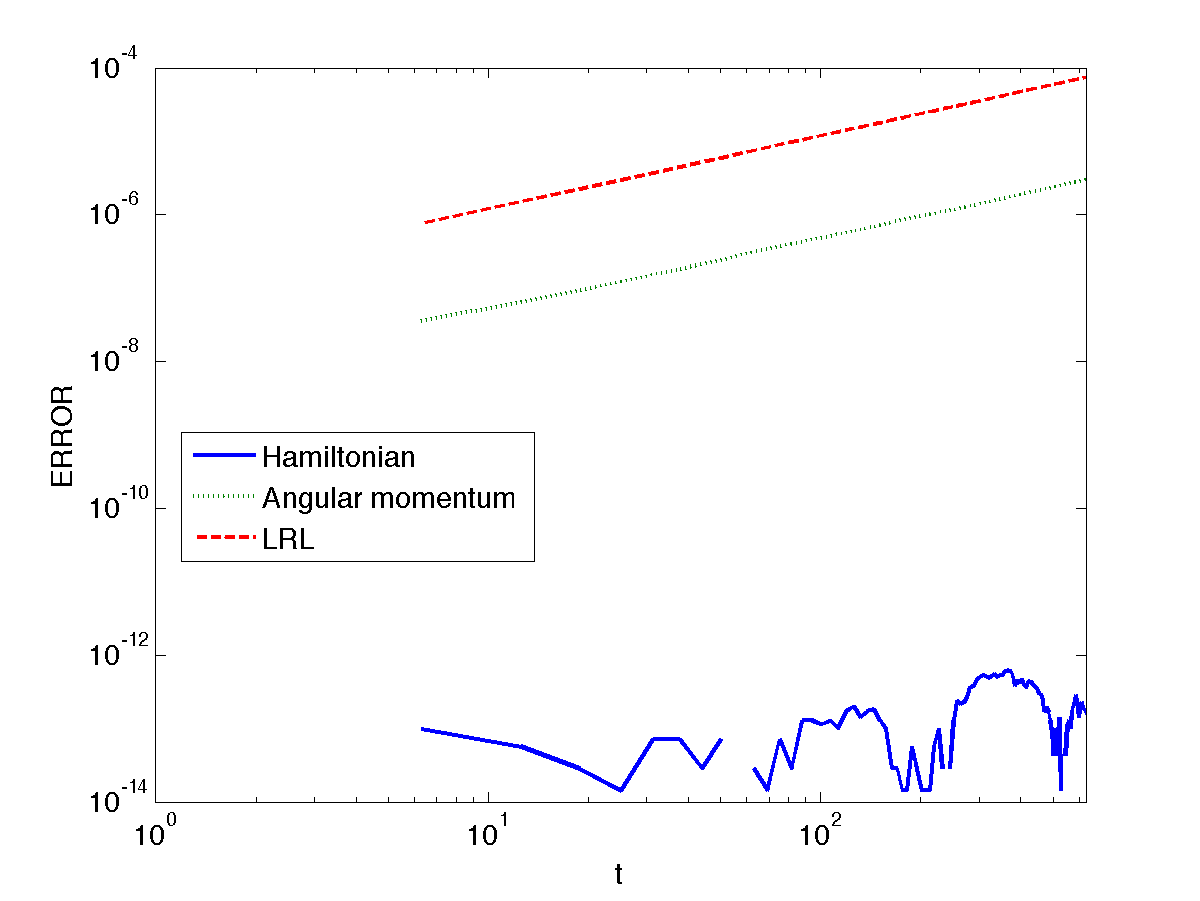}}
\caption{Kepler problem, $\eps=0.99$, HBVM(8,2) method, $tol=10^{-8}$, invariants errors.}
\label{kepl082v}

\bigskip
\centerline{\includegraphics[width=12cm,height=7cm]{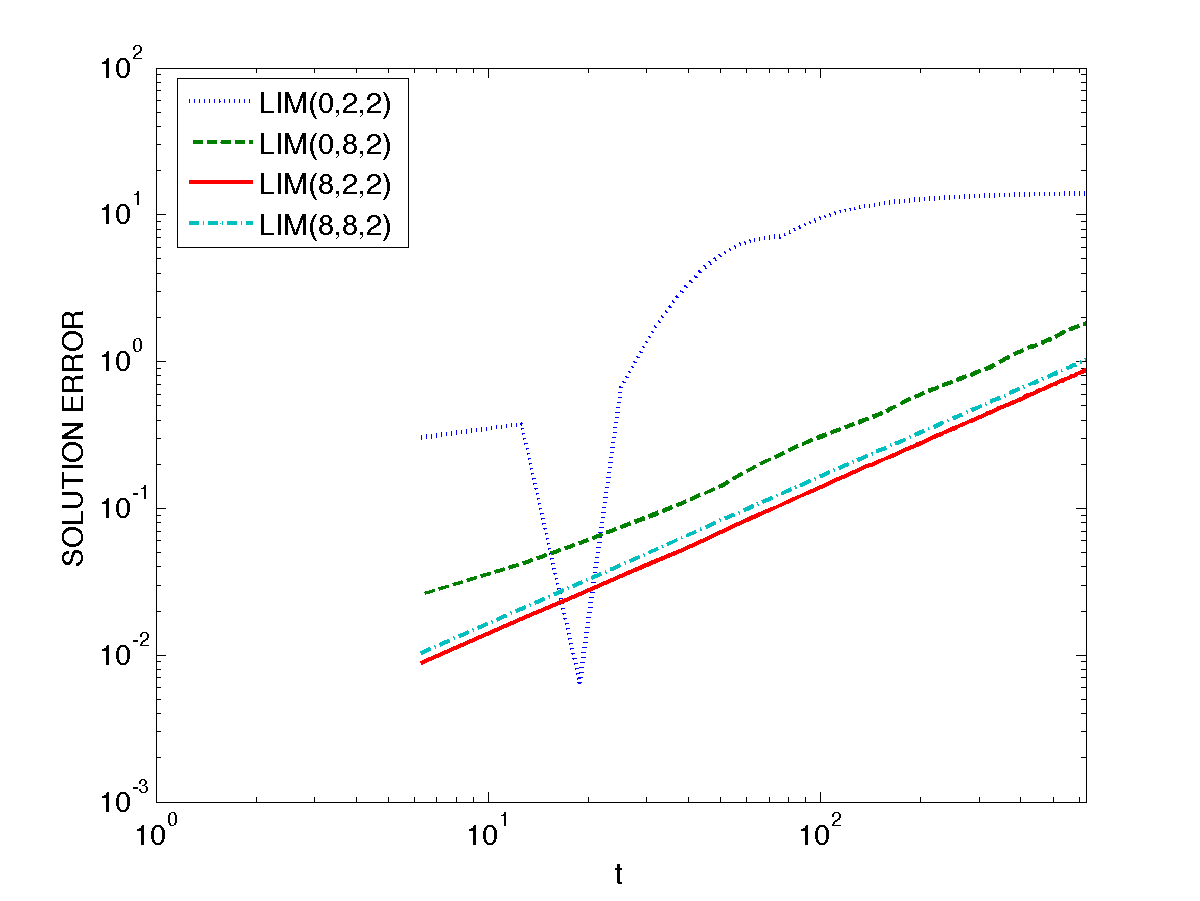}}
\caption{Kepler problem, $\eps=0.99$, error in the numerical solution, $tol=10^{-8}$.}
\label{keplerrv}
\end{figure}

This scenario changes when the eccentricity $\eps\approx1$: indeed, in such a case a constant stepsize is very inefficient and a variable stepsize would be preferable. By using  a standard mesh selection strategy, based on the control of the local error, e.g.,

$$h_{new} = 0.85\,h_{old}\left( \frac{tol}{\|\bfe\|}\right)^{\frac{1}{p+1}},$$
where:

\begin{itemize}
\item $h_{old}$ is the current stepsize;
\item $h_{new}$ is the new stepsize;
\item 0.85 is a ``safety'' factor;
\item $tol$ is the prescribed tolerance for the local error;
\item $\bfe$ is an esatimate of the latter error;
\item $p$ is the order of the method;
\end{itemize}
it is known that symplectic methods may suffer from a quadratic error growth (instead of a linear one, as shown, e.g., in \cite[pp.\,303--305]{HLW06}). Let us then see what happens when $\eps=0.99$, and a tolerance $tol=10^{-8}$ is used for all the above methods:
\begin{itemize}
\item for the 2-stage Gauss method, from the plot in Figure~\ref{kepl022v} one has that a drift in both the Hamiltonian and the LRL invariant is present, whereas the angular momentum is preserved up to roundoff;

\item for the HBVM(8,2) method, from the plot in Figure~\ref{kepl082v} one has that a drift in both the angular and the LRL invariant is present, whereas the Hamiltonian is preserved up to roundoff;

\item for both LIM(8,2,2) and LIM(8,8,2), all the invariants are conserved up to roundoff.
\end{itemize}
At last, in Figure~\ref{keplerrv} there is the plot of the error, measured at each period, for 100 periods: it is evident that the solution provided by the 2-stage Gauss method soon becomes meaningless, whereas for all the other methods a linear error growth is observed, the fully conserving methods being slightly more accurate than HBVM(8,2).

\section*{The Lotka-Volterra problem}
\begin{figure}
\centerline{\includegraphics[width=12cm,height=7cm]{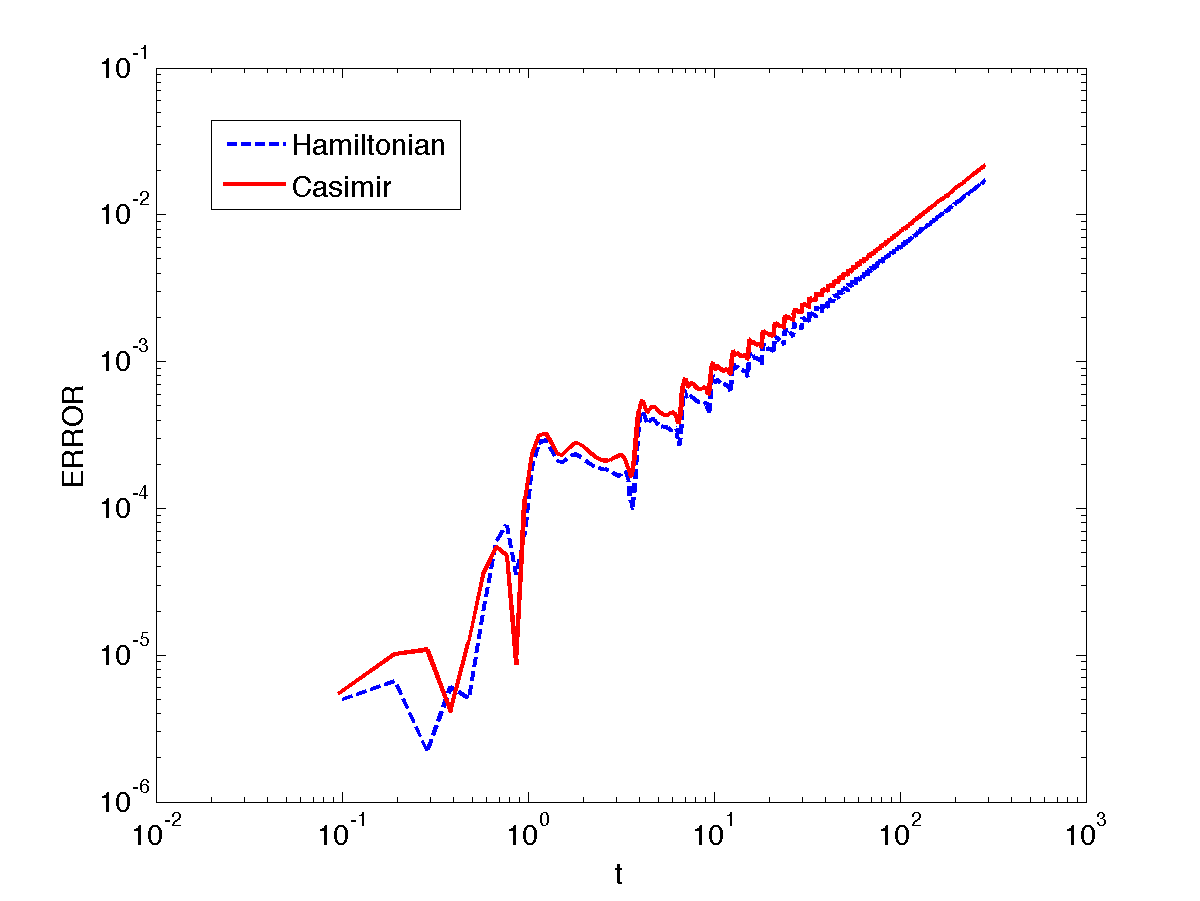}}
\caption{Lotka-Volterra problem, 2-stage Gauss method, $h = T/30$, invariants errors.}
\label{lotka022}

\bigskip
\centerline{\includegraphics[width=12cm,height=7cm]{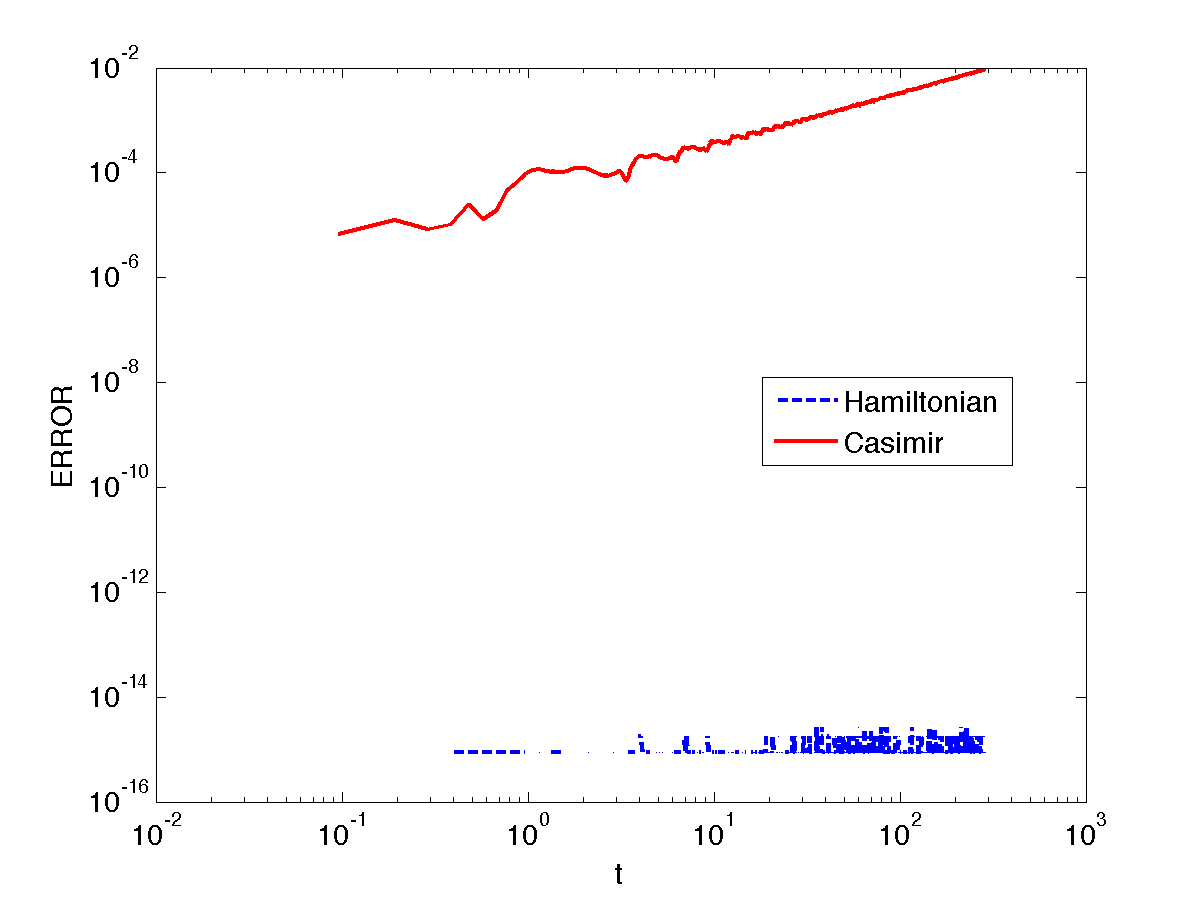}}
\caption{Lotka-Volterra problem, LIM(8,2,2) method with only the Hamiltonian preserved, $h = T/30$, invariants errors.}
\label{lotka822}

\bigskip
\centerline{\includegraphics[width=12cm,height=7cm]{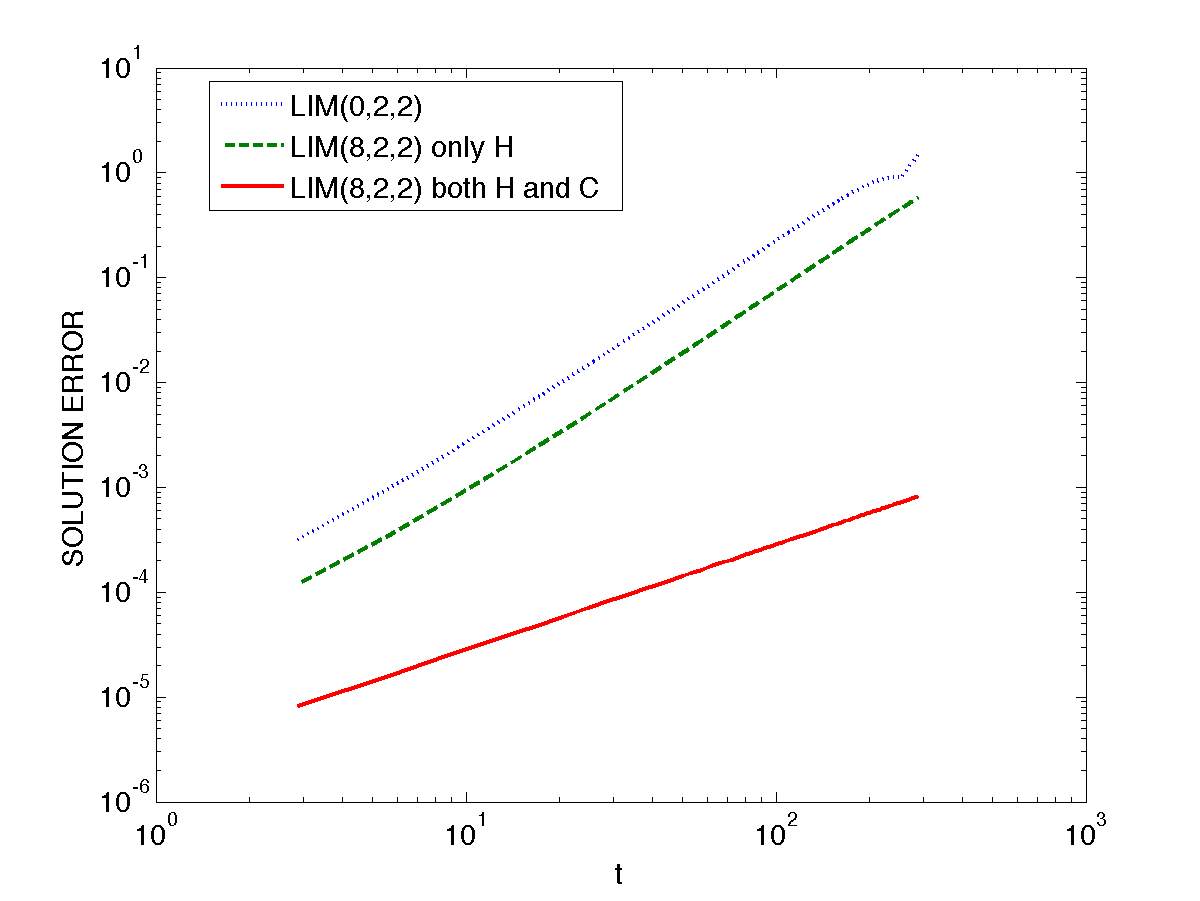}}
\caption{Lotka-Volterra problem, solution errors with stepsize $h = T/30$.}
\label{lotkaerr}
\end{figure}

The second test problem that we consider is the Lotka-Volterra problem, i.e., a problem in Poisson form,
\begin{equation}\label{poisson}
\bfy' = B(\bfy) \nabla H(\bfy), \qquad \bfy(0)=\bfy_0,
\end{equation}
with $$B(\bfy)^T=-B(\bfy), \qquad \forall \bfy,$$
and the scalar function $H(\bfy)$ is still called the Hamiltonian. Also in this case, the Hamiltonian is a constant of motion, since:
$$\frac{\dd}{\dd t}H(\bfy(t)) = \nabla H(\bfy(t))^T\bfy'(t) = \nabla H(\bfy)^TB(\bfy(t)\nabla H(\bfy(t)) = 0,$$
 $B(\bfy(t))$ being skew-symmetric. Moreover, each function $C(\bfy)$ such that
\begin{equation}\label{casimir}
\nabla C(\bfy)^TB(\bfy) = 0,
\end{equation}
is an invariant for the corresponding dynamical system. Indeed, one has:
$$\frac{\dd}{\dd t}C(\bfy(t)) = \nabla C(\bfy(t))^T\bfy'(t) = \nabla C(\bfy)^TB(\bfy(t))\nabla H(\bfy(t)) = 0,$$
because of (\ref{casimir}). A function $C(\bfy)$ which satisfies (\ref{casimir}) is said a {\em Casimir} for
(\ref{poisson}). Let us then consider the following problem \cite{CH11}, for which $\bfy=(y_1,y_2,y_3)^T$,
$$B(\bfy) = \pmatrix{ccc}
0 &c y_1y_2 &bcy_1y_3\\
-cy_1y_2 &0 & -y_2y_3\\
-bcy_1y_3 &y_2y_3 &0\endpmatrix, \qquad abc = -1,$$
the Hamiltonian is
$$H(\bfy) = aby_1 +y_2 -ay_3 +\nu\log y_2 - \mu\log y_3,$$
and, moreover, there is the following Casimir:
$$C(\bfy) = ab\log y_1 -b\log y_2 +\log y_3.$$
By using the following set of parameters,
$$a = -2,\qquad b = -1,\qquad c = -0.5,\qquad \nu = 1,\qquad \mu = 2,$$
and initial point:
$$\bfy_0 = \pmatrix{ccc} 1 & 1.9 & 0.5\endpmatrix^T,$$
the solution turns out to be periodic with period
$$T\approx 2.878130103817.$$
We now solve this problem with a constant stepsize
$$h = T/30 \approx 0.096,$$
so that we can check both the errors in the solution and in the invariants.
We solve, at first, the problem by using the 2-stage Gauss method. In Figure~\ref{lotka022} we plot the errors in the invariants along the numerical solution: as one can see, both of them exhibit a drift. Then, we use LIM(8,2,2), with the same stepsize, by imposing {\em only} the Hamiltonian conservation: indeed, in the case of the Kepler problem, this was sufficient to obtain a linear growth for the solution error. In Figure~\ref{lotka822} we plot the error in the numerical Hamiltonian and Casimir, thus showing a practical conservation of the former, and a linear drift for the latter. Finally, we use LIM(8,2,2), with the same stepsize, by imposing {\em both} the conservation of the Hamiltonian and of the Casimir, which are conserved up to roundoff. At last, in Figure~\ref{lotkaerr} we plot the measured error in the solution, measured over 100 periods: one then concludes that a linear error growth is attained only when preserving {\em both} the invariants; conversely, a quadratic error growth occurs.

\chapter{Further developments and references}
We end these lecture notes, by adding that further interesting developments, such as the possibility of getting, in a weakened sense, methods which are {\em both} symplectic and energy-conserving, have been considered in \cite{BIT12_3,BIT10_2}. We also mention that a noticeable extension of this approach, for PRK methods, has been recently devised in \cite{WXL13}. A further line of investigation deals with multistep energy-preserving method, as is sketched in \cite{BIT12_2}. Last, but not least, the efficient implementation of such methods deserves to be investigated as well.

\nopagebreak


\begin{thebibliography}{99}

\bibitem{BG94} G.\,Benettin, A.\,Giorgilli. On the Hamiltonian interpolation of near to the identity symplectic mappings with application to symplectic integration algorithms. {\em J. Statist. Phys.} {\bf  74} (1994) 1117--1143.

\bibitem{BS00}
P.\,Betsch, P.\,Steinmann. Inherently Energy Conserving Time Finite Elements for Classical Mechanics. {\em Journal of Computational Physics} {\bf 160} (2000) 88--116.

\bibitem{Bo97}
C.L.\,Bottasso. A new look at finite elements in time: a variational interpretation of Runge--Kutta methods. {\em Applied Numerical Mathematics} {\bf 25} (1997) 355--368.

\bibitem{Br00}
L.\,Brugnano. Blended Block BVMs (B3VMs): A Family of Economical Implicit Methods for ODEs. {\em Journal of Computational and Applied Mathematics} {\bf 116} (2000) 41--62.

\bibitem{BCMR12}
L.\,Brugnano, M.\,Calvo, J.I.\,Montijano, L.\,R\`andez.  Energy preserving methods for Poisson systems. {\em Journal of Computational and Applied Mathematics} {\bf 236} (2012) 3890--3904.

\bibitem{BI12}
L.\,Brugnano, F.\,Iavernaro.  Line Integral Methods which preserve all invariants of conservative problems. {\em Journal of Computational and Applied Mathematics} {\bf 236} (2012) 3905--3919.

\bibitem{BI12_1} L.\,Brugnano, F.\,Iavernaro. Recent Advances in the Numerical Solution of Conservative Problems. {\em AIP Conference Proc.} {\bf 1493} (2012) 175--182.

\bibitem{BI12_2}  L.\,Brugnano, F.\,Iavernaro. Geometric Integration by Playing with Matrices. {\em AIP Conference Proceedings} {\bf 1479} (2012) 16--19. 

\bibitem{BIS10}
 L.\,Brugnano, F.\,Iavernaro, T.\,Susca. Numerical comparisons between Gauss-Legendre methods and Hamiltonian BVMs defined over Gauss points. {\em Monografias de la Real Acedemia de Ciencias de Zaragoza} {\bf 33} (2010) 95--112.

\bibitem{BIT09}
L.\,Brugnano, F.\,Iavernaro, D.\,Trigiante. Analisys of Hamiltonian Boundary Value Methods (HBVMs) for the numerical solution of polynomial Hamiltonian dynamical systems. (2009) \url{arXiv:0909.5659v1}

 \bibitem{BIT09_1}
L.\,Brugnano, F.\,Iavernaro, D.\,Trigiante. Hamiltonian BVMs (HBVMs): a family of  "drift-free" methods for integrating polynomial Hamiltonian systems. {\em AIP Conf. Proc.} {\bf 1168} (2009) 715--718.

\bibitem{BIT10_0}L.\,Brugnano, F.\,Iavernaro, D.\,Trigiante. {\em The Hamiltonian BVMs (HBVMs) Homepage}, 2010. {\tt  arXiv:1002.2757}

\bibitem{BIT10}
L.\,Brugnano, F.\,Iavernaro, D.\,Trigiante. Hamiltonian Boundary Value Methods (Energy Preserving Discrete Line Methods). {\em Journal of Numerical Analysis, Industrial and Applied Mathematics} {\bf 5},1-2 (2010) 17--37.

\bibitem{BIT10_1}
 L.\,Brugnano, F.\,Iavernaro, D.\,Trigiante. Numerical Solution of ODEs and the Columbus' Egg: Three Simple Ideas for Three Difficult Problems. {\em Mathematics in Engineering, Science and Aerospace} {\bf 1},4 (2010) 407--426.
  
\bibitem{BIT10_2}  L.\,Brugnano, F.\,Iavernaro, D.\,Trigiante. Energy and quadratic invariants preserving integrators of Gaussian type. {\em AIP Conference Proceedings} {\bf 1281} (2010) 227--230.
 
\bibitem{BIT11}
 L.\,Brugnano, F.\,Iavernaro, D.\,Trigiante. A note on the efficient implementation of Hamiltonian BVMs. {\em Journal of Computational and Applied Mathematics} {\bf 236} (2011) 375--383.

\bibitem{BIT12}
  L.\,Brugnano, F.\,Iavernaro, D.\,Trigiante.  The Lack of Continuity and the Role of Infinite and Infinitesimal in Numerical Methods for ODEs: the Case of Symplecticity. {\em Applied Mathematics and Computation} {\bf 218} (2012) 8053--8063.

\bibitem{BIT12_1}
 L.\,Brugnano, F.\,Iavernaro, D.\,Trigiante.  A simple framework for the derivation and analysis of effective one-step methods for ODEs. {\em Applied Mathematics and Computation} {\bf 218} (2012) 8475--8485.

\bibitem{BIT12_2}
L.\,Brugnano, F.\,Iavernaro, D.\,Trigiante. A two-step, fourth-order method with energy preserving properties. {\em Computer Physics Communications} {\bf 183} (2012) 1860--1868.

\bibitem{BIT12_3}
L.\,Brugnano, F.\,Iavernaro, D.\,Trigiante. Energy and QUadratic Invariants Preserving integrators based upon Gauss collocation formulae. {\em SIAM Journal on Numerical Analysis} {\bf 50}, No. 6 (2012) 2897--2916.

 \bibitem{BM02} L.\,Brugnano, C.\,Magherini. Blended Implementation of Block Implicit Methods for ODEs. {\em Appl. Numer. Math.} {\bf 42} (2002) 29--45.

\bibitem{BM04} L.\,Brugnano, C.\,Magherini. The BiM Code for the Numerical Solution of ODEs. {\em Jour. Comput. Appl. Mathematics}  {\bf 164-165} (2004) 145--158.

\bibitem{BM09} L.\,Brugnano, C.\,Magherini. Recent Advances in Linear Analysis of Convergence for Splittings for Solving ODE problems. {\em Applied Numerical Mathematics} {\bf 59} (2009) 542--557.

\bibitem{BMM06} L.\,Brugnano, C.\,Magherini, F.\,Mugnai. Blended Implicit Methods for the Numerical Solution of DAE Problems. {\em Jour. Comput. Appl. Mathematics}  {\bf 189} (2006) 34--50.

\bibitem{BT98}
L.\,Brugnano, D.\,Trigiante. {\em Solving ODEs by Linear Multistep Initial and Boundary Value Methods},  Gordon and Breach, Amsterdam, 1998.

\bibitem{BB12} K.\,Burrage, P.M.\,Burrage. Low rank Runge-Kutta methods, symplecticity and stochastic Hamiltonian problems with additive noise. {\em Journal of Computational and Applied Mathematics}  {\bf 236} (2012) 3920--3930.

\bibitem{BB79} K.\,Burrage, J.C.\,Butcher. Stability criteria for implicit Runge--Kutta methods. {\em SIAM Journal on Numerical Analysis} {\bf 16} (1979) 46--57.

\bibitem{CLMR11}
M.\,Calvo, M.P.\,Laburta, J.I.\,Montijano, L.\,R\'andez. Error growth in the numerical integration of periodic orbits, {\em Math. Comput. Simulation} {\bf 81} (2011) 2646--2661.

\bibitem{CMcLMcLOQW09}
E.\,Celledoni, R.I.\,McLachlan, D.\,Mc\,Laren, B.\,Owren, G.R.W.\,Quispel, W.M.\,Wright. Energy preserving Runge--Kutta methods. {\em M2AN} {\bf 43} (2009) 645--649.

\bibitem{Bseries} E.\,Celledoni, R.I.\,McLachlan, B.\,Owren, G.R.W.\,Quispel. Energy-Preserving Integrators and the Structure of B-series. {\em Found. Comput. Math.} {\bf 10} (2010) 673--693.

\bibitem{CFM06}
P.\,Chartier, E.\,Faou, A.\,Murua. An algebraic approach to invariant preserving integrators: the case of quadratic and Hamiltonian invariants. {\em Numer. Math.} {\bf 103},\,4 (2006) 575--590.

\bibitem{CH11}D.\,Cohen, E.\,Hairer. Linear energy-preserving integrators for Poisson systems. {\em BIT Numer. Math.} {\bf 51} (1) (2011) 91--101.

\bibitem{C79} M.\,Crouzeix. Sur la B-stabilit\'e des m\'ethodes de Runge--Kutta. {\em Numerische Mathematik} {\bf 32} (1979) 75--82.

\bibitem{FK85} Feng Kang. On Difference Schemes and Symplectic Geometry. In {\em Proceedings of the 1984 Beijing symposium on differential geometry and differential equations}. Science Press, Beijing, 1985, pp.\,42--58.

\bibitem{GM88}
Z.\,Ge, J.E.\,Marsden. Lie-Poisson Hamilton-Jacobi theory and Lie-Poisson integrators. {\em Phys. Lett. A} {\bf 133} (1988) 134--139.

\bibitem{GPS01} H.\,Goldstein, C.P.\,Poole, J.L.\,Safko. {\em Classical Mechanics}. Addison Wesley, 2001.

\bibitem{Go96}
O.\,Gonzales. Time integration and discrete Hamiltonian systems. {\em J. Nonlinear Sci.} {\bf 6} (1996) 449--467.

\bibitem{G75}
W.\,Gr\"obner. {\em Gruppi, Anelli e Algebre di Lie}. Collana di Informazione Scientifica ``Poliedro'', Edizioni Cremonese, Rome, 1975.

\bibitem{Ha10}
E.\,Hairer. Energy preserving variant of collocation methods. {\em Journal of Numerical Analysis, Industrial and Applied Mathematics} {\bf 5},1-2 (2010) 73--84.

\bibitem{HLW06}  E.\,Hairer, C.\,Lubich, G.\,Wanner. {\em Geometric Numerical Integration. Structure-Preserving Algorithms for Ordinary Differential Equations}, Second ed., Springer, Berlin, 2006.

\bibitem{HW96} E.\,Hairer, G.\,Wanner. {\em Solving Ordinary Differential Equations II. Stiff and Differential-Algebraic Problems, 2nd edition}. Springer-Verlag, Berlin, 1996.

\bibitem{HZ12} E.\,Hairer, C.J.\,Zbinden. On conjugate symplecticity of B-series integrators. {\em IMA J. Numer. Anal.} (2012) 1--23.

\bibitem{HoS97} P.J.\,van der Houwen, J.J.B.\,de Swart. Triangularly implicit iteration methods for ODE-IVP solvers. {\em SIAM J.\,Sci. Comput.} {\bf 18} (1997) 41--55. 

\bibitem{HoS97a} P.J.\,van der Houwen, J.J.B.\,de Swart. Parallel linear system solvers for Runge-Kutta methods. {\em Adv. Comput. Math.} {\bf 7},\,1-2 (1997) 157--181.

\bibitem{Hu72}
B.L.\,Hulme. One-Step Piecewise Polynomial Galerkin Methods for Initial Value Problems. {\em Mathematics of Computation}, {\bf 26}, 118 (1972) 415--426.

\bibitem{IP07}
F.\,Iavernaro, B.\,Pace. $s$-Stage Trapezoidal Methods for the Conservation of Hamiltonian Functions of Polynomial Type. {\em AIP Conf. Proc.} {\bf 936} (2007) 603--606.

\bibitem{IP08}
F.\,Iavernaro, B.\,Pace. Conservative Block-Boundary Value Methods for the Solution of Polynomial Hamiltonian Systems. {\em AIP Conf. Proc.} {\bf 1048} (2008) 888--891.

\bibitem{IT09}
F.\,Iavernaro, D.\,Trigiante. High-order symmetric schemes for the energy conservation of polynomial Hamiltonian problems. {\em Journal of Numerical Analysis, Industrial and Applied Mathematics} {\bf 4},1-2 (2009) 87--101.

\bibitem{KMO99}
C.\,Kane, J.E.\,Marsden, M.\,Ortiz. Symplectic-energy-momentum preserving variational integrators. {\em Jour. Math. Phys.} {\bf 40},\,7 (1999) 3353--3371.

\bibitem{LR88} V.\,Lakshmikantham, D.\,Trigiante. {\em Theory of Difference Equations. Numerical Methods and Applications}. Academic Press, 1988.

\bibitem{McLQR99}
R.I.\,Mc\,Lachlan, G.R.W.\,Quispel, N.\,Robidoux. Geometric integration using discrete gradient. {\em Phil. Trans. R. Soc. Lond. A} {\bf 357} (1999) 1021--1045.

\bibitem{MW97} J.E.\,Marsden, J.M.\,Wendlandt.  Mechanical Systems with Symmetry,
Variational Principles, and Integration Algorithms. in {\em ``Current and Future Directions in Applied Mathematics''
M.\,Alber, B.\,Hu, and J.\,Rosenthal, Eds.}, Birkh\"auser, 1997, pp.\,219--261.

\bibitem{QMcL08}
G.R.W.\,Quispel, D.I.\,Mc\,Laren. A new class of energy-preserving numerical integration methods.
{\em J. Phys. A: Math. Theor.} {\bf 41} (2008) 045206 (7pp).

\bibitem{SS88} J.M.\,Sanz Serna. Runge-Kutta schemes for Hamiltonian systems. {\em BIT} {\bf 28} (1988) 877--883.

\bibitem{ST94} J.C.\,Simo, N.\,Tarnow. A new energy and momentum conserving algorithm for the non-linear dynamics of shells. {\em Internat. Jour. for Numerical Meth. in Engineering}  {\bf 37} (1994) 2527--2549.

\bibitem{STW92} J.C.\,Simo, N.\,Tarnow and K.K.\,Wong. Exact energy-momentum conserving algorithms and symplectic schemes for nonlinear dynamics. {\em Computer Methods in Applied Mechanics and Engineering} {\bf 100} (1992) 63--116.

\bibitem{Suris89} Y.B.\,Suris. On the canonicity of mappings that can be generated by methods of Runge--Kutta type for integrating systems $x''=\partial U/\partial x$. {\em U.S.S.R. Comput. Math. and Math. Phys.} {\bf 29},\,1 (1989) 138--144.

\bibitem{TC07}
Q.\,Tang, C.-m.\,Chen. Continuous finite element methods for Hamiltonian systems. {\em
Applied Mathematics and Mechanics}  {\bf 28},8  (2007) 1071--1080.

\bibitem{TS12}
W.\,Tang, Y.\,Sun.
Time finite element methods: a unified framework for numerical discretizations of ODEs.
{\em Applied Mathematics and Computation} {\bf 219},\,4 (2012)  2158--2179.

\bibitem{WXL13}
D.\,Wang, A.\,Xiao, X.\,Li. Parametric symplectic partitioned Runge-Kutta methods with energy-preserving properties for Hamiltonian systems. {\em Computer Physics Communications} {\bf 184} (2013) 303--310.


\end{thebibliography}
\end{document}